\documentclass[final,leqno,letterpaper]{etna}
\usepackage[utf8]{inputenc}
\usepackage[margin=1in]{geometry}
\usepackage{ragged2e}
\usepackage{amsmath}
\usepackage{amssymb}
\usepackage{bbm}
\usepackage{enumitem}
\usepackage{graphicx}
\usepackage{pdfpages}
\usepackage{tikz}
\usepackage{sgame}
\usepackage{varwidth}
\usepackage{mdframed}
\usepackage[nameinlink,capitalize]{cleveref}
\usepackage{hyperref}
\usepackage{microtype}
\usepackage{mathrsfs}
\usepackage[font=small,labelfont=bf,textfont=it]{caption}
\usepackage[font=small,labelfont=bf,textfont=it]{subcaption}
\crefname{equation}{}{}
\usepackage{url}
\renewcommand{\thefootnote}{\fnsymbol{footnote}}
\usepackage[space,noadjust,nocompress]{cite}

\usepackage{algorithm}
\usepackage{algpseudocode}

\allowdisplaybreaks

\usepackage{mathtools}
\usepackage{titlesec}

\titleformat{\paragraph}
{\normalfont\normalsize\bfseries}{\theparagraph}{1em}{}
\titlespacing*{\paragraph}
{0pt}{3.25ex plus 1ex minus .2ex}{1.5ex plus .2ex}

\newcommand{\mycomment}[1]{}

\usepackage{xcolor}

\newcommand{\ptxt}[1]{{\color{black} #1}}

\shorttitle{PARALLEL-IN-ITERATION OPTIMIZATION USING MGRIT} 
\shortauthor{G.~ARAUJO, O.~KRZYSIK, AND H.~DE STERCK}

\begin{document}
\title{Parallel-in-iteration optimization using multigrid reduction-in-time}

\author{G. H. M. Araújo\footnotemark[2] \and O. A. Krzysik \footnotemark[3] \and H. De Sterck\footnotemark[2]}

\footnotetext[2]{Department of Applied Mathematics, University of Waterloo, Waterloo, Ontario, Canada ({\tt \{ghmacieiradearaujo, hans.desterck\}@uwaterloo.ca}).}
\footnotetext[3]{Theoretical Division, Los Alamos Laboratory, Los Alamos, NM, USA ({\tt okrzysik@lanl.gov}, \url{https://orcid.org/0000-0001-7880-6512}). The work of the second author has been funded by the Los Alamos National Laboratory Advanced Simulation and Computation program through Contract No. 89233218CNA000001. The research was performed under the auspices of the National Nuclear Security Administration of the U.S. Department of Energy at Los Alamos National Laboratory, managed by Triad National Security, LLC under contract
89233218CNA000001. LA-UR-25-27245.}

    \date{\vspace{-5ex}}
	\maketitle
    
\renewcommand{\thefootnote}{\fnsymbol{footnote}}

    \begin{abstract}
    Standard gradient-based iteration algorithms for optimization, such as gradient descent and its various proximal-based extensions to nonsmooth problems, are known to converge slowly for ill-conditioned problems, sometimes requiring many tens of thousands of iterations in practice. 
    Since these iterations are computed sequentially, they may present a computational bottleneck in large-scale parallel simulations.
    In this work, we present a ``parallel-in-iteration'' framework that allows one to parallelize across these iterations using multiple processors with the objective of reducing the wall-clock time needed to solve the underlying optimization problem.
    Our methodology is based on re-purposing parallel time integration algorithms for time-dependent differential equations, motivated by the fact that optimization algorithms often have interpretations as discretizations of time-dependent differential equations (such as gradient flow). 
    Specifically in this work, we use the parallel-in-time method of multigrid reduction-in-time (MGRIT), but note that our approach permits in principle the use of any other parallel-in-time method.
    We numerically demonstrate the efficacy of our approach on  two different model problems, including a standard convex quadratic problem and the nonsmooth elastic obstacle problem in one and two spatial dimensions. 
    For our model problems, we observe fast MGRIT convergence analogous to its prototypical performance on partial differential equations of diffusion type.
    Some theory is presented to connect the convergence of MGRIT to the convergence of the underlying optimization algorithm.
    Theoretically predicted parallel speedup results are also provided.

    \end{abstract}
    
    

	\section{Introduction}\label{sec:intro}
    
	Consider the unconstrained optimization problem, 
	\begin{equation} \label{eq:min_f_intro}
			\underset{\mathbf{u} \in \mathbb{R}^N}{\min}\quad f(\mathbf{u}), 
	\end{equation}
	with $f: \mathbb{R}^N \rightarrow \mathbb{R}$ a convex and $L$-smooth function, i.e., a convex function with $L$-Lipschitz gradient:
    \begin{equation}\label{eq:lipschitz-smooth}
        \|\nabla f(\mathbf{u})-\nabla f(\mathbf{v})\| \leqslant L\|\mathbf{u}-\mathbf{v}\| ~\forall \mathbf{u},\mathbf{v} \in \mathbb{R}^N.
    \end{equation}
    Such optimization problems are ubiquitous throughout the sciences, with many theoretical properties established for them and applications in fields such as machine learning \cite{Suykens1999, Tai2021}, signal reconstruction \cite{Beck2009, Golub1999, Chambolle2016}, game theory \cite{Gao2017}, econometrics \cite{Angrist2008} and mathematical physics \cite{Rodrigues2014, Tran2015}.
    Various iterative gradient-based optimization algorithms for solving unconstrained minimization problems, such as gradient descent and Nesterov's accelerated gradient descent \cite{Nesterov1983} take the form of 1- or 2-step updates.
    Generally, the convergence rates of these algorithms for convex functions can be slow. For example, gradient descent has a convergence rate of $\mathcal{O}(\frac{1}{k})$ for problem \cref{eq:min_f_intro}\ptxt{, where $k$ denotes the number of iterations}; while Nesterov's improved rate of $\mathcal{O}(\frac{1}{k^2})$ was shown in \cite{Nesterov2004} to be optimal among all methods that only have information about the gradient of $f$ for its updates. 
    In practical terms, this means that many iterations are typically required to solve the problem; for example, tens of thousands of iterations is not uncommon \cite{Chambolle2016, ODonoghue2013}.  
    %
    %

    Related to \eqref{eq:min_f_intro} is the unconstrained problem of minimizing a sum of two functions
	\begin{equation}
		\min_{\mathbf{u} \in \mathbb{R}^{N}}\quad F(\mathbf{u}), \text{ where } F(\mathbf{u}):=f(\mathbf{u})+g(\mathbf{u}), \label{eq:min_f+g_intro}
	\end{equation}
	where $f: \mathbb{R}^N \rightarrow \mathbb{R}$ is convex and $L$-smooth, and $g: \mathbb{R}^{N} \rightarrow (-\infty,+\infty]$ is continuous, convex and possibly nondifferentiable.	
    The minimization problem \eqref{eq:min_f+g_intro} is typically solved by splitting algorithms \cite{Lions1979}, in which the differentiable $f$ and possibly nondifferentiable $g$ are handled differently. Many classical splitting algorithms also take the form of 1- or 2-step updates. Namely, the proximal gradient and FISTA algorithms \cite{Beck2009} are generalizations of gradient descent and Nesterov's accelerated gradient, respectively, for solving \cref{eq:min_f+g_intro}, maintaining the convergence rates of $\mathcal{O}(\frac{1}{k})$ and $\mathcal{O}(\frac{1}{k^2})$ as their counterparts do for \cref{eq:min_f_intro}. The convergence rates of these optimization algorithms are also slow, in practice requiring many iterations to get close to the true optimum \cite{Beck2009, Chambolle2016}. If gradient evaluations take some nontrivial amount of compute time, we may have to wait for a long wall-clock time to get an accurate approximation to the optimum because the algorithms are sequential in $k$.
	
	One can write 1-step algorithms in the form $\mathbf{u}_{k+1}=\Phi( \mathbf{u}_k) + \mathbf{w}_{k+1}$, where $\Phi(\cdot)$ is a linear or nonlinear operator and $\mathbf{w}_{k+1}$ is a generic term that does not dpeend on $\mathbf{u}_k$, for example, to encode the initial condition or a forcing term. For example, we can write the gradient descent update for \cref{eq:min_f_intro} given by
        \begin{equation}\label{eq:gd_intro}
		\mathbf{u}_{k+1}=\mathbf{u}_{k}-s \nabla f(\mathbf{u}_{k}), \quad k=0, 1, \ldots
	\end{equation}
	as 
	\begin{equation}
		\mathbf{u}_{k+1}=(I-s \nabla f)(\mathbf{u}_{k}).
	\end{equation}
	Suppose we need $K$ iterations to approximate the true optimum within the desired accuracy. Then, taking $\Phi:=(I-s \nabla f)$, computing the $K$ gradient descent iterations in parallel is equivalent to solving the below ``all-at-once'' system in parallel 
	\begin{equation}\label{eq:Au=g_intro}
		\mathcal{A}(\mathbf{u})=\left[\begin{array}{cccc}
			I & & & \\
			-\Phi(\cdot) & I & & \\
			& \ddots & \ddots & \\
			& & -\Phi(\cdot) & I
		\end{array}\right]\left[\begin{array}{c}
			\mathbf{u}_0 \\
			\mathbf{u}_1 \\
			\vdots \\
			\mathbf{u}_K
		\end{array}\right]=\left[\begin{array}{c}
			\mathbf{w}_0 \\
			\mathbf{w}_1 \\
			\vdots \\
			\mathbf{w}_K
		\end{array}\right]=\mathbf{w},
	\end{equation}
    where for this particular case $\mathbf{w}_0$ is the initial condition and $\mathbf{w}_1=\cdots=\mathbf{w}_K=\mathbf{0}$.
        We note that some 2-step methods, such as Nesterov's accelerated gradient descent \cite{Nesterov1983} for solving \cref{eq:min_f_intro}, can be rewritten as 1-step algorithms of the form $[\mathbf{u}_{k+1} \quad \mathbf{v}_{k+1}]^T=\Phi([\mathbf{u}_{k} \quad \mathbf{v}_{k}]^T)$ for block vector inputs $[\mathbf{u}_k \quad \mathbf{v}_k]^T$, such that an analogous all-at-once system to \cref{eq:Au=g_intro} holds in such cases too.
        
        The system \cref{eq:Au=g_intro} is exactly of the kind which parallel-in-time algorithms solve for time-dependent evolution equations \cite{Gander2015, Ong2020}. For example, the parareal method \cite{Lions2001} is a 2-level parallel-in-time method \cite{Gander2007} which solves systems of differential equations by employing sequential coarse-grid solves to speed up solving the original fine-grid problem in parallel. However, parareal is limited in concurrency, as the coarse-grid solve is still sequential. Thus, the multigrid-reduction-in-time algorithm (MGRIT) \cite{Falgout2014} was developed as a multilevel generalization of parareal based on multigrid reduction (MGR) techniques \cite{Ries1983}.

	In the linear case where $\Phi(\mathbf{u}_k)$ is computed by multiplying $\mathbf{u}_k$ by the matrix $\Phi$, MGRIT solves the system \cref{eq:Au=g_intro} via a coarse-grid system derived from an approximation to the Schur complement of \cref{eq:Au=g_intro} in terms of the coarse-level variables, and has a factor of $m \in \mathbb{N}$ fewer variables on each successively coarser grid. 
    Here, $m$ is called the coarsening factor. The biggest challenge in applying MGRIT is developing a suitable coarse-grid operator $\Phi_\Delta \approx \Phi^m$, where $\Phi^m$ is the operator appearing in the Schur complement.
    That is, $\Phi_\Delta $ should approximate $m$ successive applications with $\Phi$, but it should be computationally much less expensive than applying $\Phi$ $m$ times; this cost reduction combined with a quickly converging method is what allows for speedup on a parallel machine. 
    It has been well-documented that parallel-in-time methods such as parareal and MGRIT often work very well for diffusion-dominated (parabolic-type) partial differential equation (PDE) problems \cite{Falgout2014, Dobrev2017, Howse2019}. 
    Advection-dominated (hyperbolic-type) problems are more challenging, but careful construction of coarse-grid operators may facilitate fast convergence \cite{Howse2019, DeSterck2024, DeSterck2019}. 
    This provides motivation for the use of parallel-in-time methods for optimization problems such as \cref{eq:min_f_intro} and \cref{eq:min_f+g_intro}.
	
	For sufficiently small step lengths $s$ as in \cref{eq:gd_intro}, we can sometimes interpret optimization algorithms as discretizations of certain ordinary differential equations (ODEs), in which $k$ acts as a ``discrete time index,'' and the step length $s$ is related to the time-scale of the ODE. For example as $s \rightarrow 0$, gradient descent can be thought of as an explicit Euler discretization of the gradient flow ODE,
	\begin{equation}
		\frac{\mathrm{d}}{\mathrm{d}t}\mathbf{u}=-\nabla f(\mathbf{u}(t)), \quad  t>0,
	\end{equation}
	and as $s \rightarrow 0$, Nesterov's accelerated gradient descent can be thought of as an explicit discretization \cite{Su2016} of the ODE,
	\begin{equation}
		\frac{\mathrm{d}^2}{\mathrm{d}t^2}\mathbf{u}+\frac{3}{t} \frac{\mathrm{d}}{\mathrm{d}t}\mathbf{u}+\nabla f(\mathbf{u}(t))=0, \quad t>0.
	\end{equation}

	With these connections between solving optimization problems in parallel and solving systems of differential equations using parallel-in-time methods, we present a framework for implementing MGRIT for solving optimization problems in the classes \cref{eq:min_f_intro} and \cref{eq:min_f+g_intro} in parallel. We describe two model problems belonging to the former and latter classes, and discuss the choices of fine- and coarse-grid operators for solving these problems with MGRIT. These operators are based on standard choices of algorithms for solving these problems sequentially, such that the MGRIT implementation can be interpreted as a sped-up version of these sequential algorithms.
    Parallel-in-time methods have been previously used in the context of solving PDE-constrained optimization problems either as preconditioners or for solving time-dependent subproblems that arise in the context of time-dependent PDE-constrained optimization, such as in optimal control \cite{Ulbrich2015, Vuchkov2024, Gander2020, Maday2013} or for computing the geometry of induction motors \cite{Hahne2023}. These previous applications differ from our approach, which is to use parallel-in-time methods to parallelize across iterations of some underlying sequential optimization method.
	
	The rest of this paper is organized as follows. \cref{sec:mp} and \cref{sec:opt_algos} introduce the model problems and optimization algorithms that we use to solve them, respectively, with the latter elaborating on the interpretation of optimization methods as discretizations of time-dependent ODEs. 
    \cref{sec:p-in-t} reviews the MGRIT parallel-in-time method. 
    \cref{sec:mgrit_opt} discusses our methodology for solving optimization problems in parallel with MGRIT, which explores the aforementioned connections between ODEs for which parallel-in-time methods have originally been designed and optimization problems, providing the key insight behind this paper. 
    \cref{sec:num_res} provides numerical results, while \cref{sec:speedup} provides theoretical speedup estimates. \cref{sec:nesterov} explores the direct application of the proposed MGRIT approach to momentum-accelerated algorithms and discusses the difficulties that were encountered to obtain fast convergence of MGRIT iterations for this case. \cref{sec:concl} features concluding remarks and ideas for future work.
	
	\section{Model problems}\label{sec:mp}
	We now introduce the two model problems considered throughout this paper. The first model problem is relatively simple, and serves largely  as a stepping stone towards the more challenging second model problem, which presents a more realistic application of our proposed algorithm. The optimization problems we consider have a connection to discretizations of time-dependent heat or diffusion PDEs, and thus help to motivate our use of parallel-in-time methods.

    In what follows, we consider continuous optimization problems on $d$-dimensional domains $\Omega \subset \mathbb{R}^d$, $d \in \{1, 2\}$, with boundary denoted by $\partial \Omega$. 
    We let $\mathbf{x} \in \mathbb{R}^d$ denote a coordinate on $\Omega$, such that $\mathbf{x} = x$ when $d = 1$ and $\mathbf{x} = (x, y)$ when $d = 2$.
    The domain $\Omega$ is discretized with a mesh consisting of $n+2$ total points, and $n$ interior points, in each direction, equally separated by a distance of $h > 0$, with $h = \mathcal{O}(\frac{1}{n})$.
    Of particular relevance to our model problems are finite-difference discretizations of (negative) Laplacian operators, which, on the interior mesh points, are defined as 
    \begin{subequations}
    \begin{align}
    \label{eq:1d_laplacian}
        A_1 &= -
			\frac{1}{h^2}
            \begin{bmatrix}
                -2 & 1  \\
				1 & -2 & 1 \\
				& \ddots & \ddots & \ddots \\
				& & 1 & -2
            \end{bmatrix}
             \in \mathbb{R}^{n \times n}, 
             \quad \textrm{and} \quad
             \\
             \label{eq:2d_laplacian}
             A_2 &= -\frac{1}{h^2}\begin{bmatrix}
				-T & I &  & \\
				I & -T & I &    \\
                & \ddots& \ddots & \ddots \\
				& &   I & -T
			\end{bmatrix} \in \mathbb{R}^{n^2 \times n^2} 
            \end{align}
            \end{subequations}
            for $d = 1$ and $d = 2$, respectively.
            Here, $T \in \mathbb{R}^{n \times n}$ is a tridiagonal matrix with elements $-1$, $4$ and $-1$ on the sub-, main and super-diagonals, respectively, and $I \in \mathbb{R}^{n \times n}$ is the identity matrix. The 2-norm of $A_1$ is $\left\|A_1\right\|_2=\frac{1}{h^2} 4 \sin ^2\left(\frac{n \pi}{2(n+1)}\right) \approx \frac{4}{h^2}$ for $h \ll 1$, while the 2-norm of $A_2$ is $\left\|A_2\right\|_2=\frac{1}{h^2} 8 \sin ^2\left(\frac{n \pi}{2(n+1)}\right) \approx \frac{8}{h^2}$ for $h \ll 1$. Throughout the paper, we refer to $A_1$ and $A_2$ as 1D and 2D Laplacian matrices, respectively.

	\subsection{Quadratic minimization problem (MP1)}\label{subsec:mp1}
	For our first model problem, we consider the straightforward problem of quadratic minimization, formulated as
 	\begin{equation}
        \label{eq:mp1}
		\min_{\mathbf{u} \in \mathbb{R}^{N}} \quad \frac{1}{2} \langle A_1\mathbf{u}, \mathbf{u}  \rangle - \langle \mathbf{b}, \mathbf{u} \rangle, 
	\end{equation}
	where $A_1$ is the 1D Laplacian matrix \eqref{eq:1d_laplacian}, such that $N = n$, and $\mathbf{b} \in \mathbb{R}^{N}$ is some prescribed vector. Since $A_1$ is symmetric positive definite (SPD), the objective function here $f(\mathbf{u})=\langle A_1\mathbf{u}, \mathbf{u}  \rangle - \langle \mathbf{b}, \mathbf{u} \rangle$ is strongly convex and $L$-smooth, where $L$ denotes the smallest Lipschitz constant of $\nabla f$, which, in this case is $L = \Vert A_1 \Vert$. Since $A_1$ is SPD, and hence invertible, the unique solution of \cref{eq:mp1} is given by $\mathbf{u}=A_1^{-1}\mathbf{b}$, which is attained when $\nabla f(\mathbf{u})=\mathbf{0}$.
			
	\subsection{Elastic obstacle problem (MP2)}\label{subsec:mp2}
	
    The elastic obstacle problem (EOP) \cite{Caffarelli1998} consists of finding the equilibrium state of an elastic membrane under a non-penetration constraint. The EOP we consider was initially motivated by mathematical physics \cite{Rodrigues2014} such as steady state fluid interaction, thin-plate solid dynamics and elastostatics \cite{Tran2015}.
	Various stategies have been proposed for numerically solving the EOP, including multigrid approaches such as \cite{Brandt1983, Mandel1984, Ang2024, Wu2015}. The formulation we present here follows an  $L^1$ penalty approach \cite{Tran2015}.
	
    %
    Given an obstacle function $\phi(\mathbf{x}): \Omega \rightarrow \mathbb{R}$, the EOP seeks the location of a membrane ${\widehat{u}(\mathbf{x}): \Omega \rightarrow \mathbb{R}}$  with the lowest elastic potential energy. 
    Assuming the elastic potential energy of the membrane is proportional to its  surface area, the EOP can be written as
	\begin{equation}
	\min _{\widehat{u} \in \Omega} \quad \int_{\Omega} \sqrt{1+\|\nabla \widehat{u}\|_{2}^2} ~\mathrm{d}\mathbf{x} \quad \text {s.t. }  \widehat{u} \geq \phi \text { on } \Omega,~ \widehat{u}=0 \text { on } \partial \Omega.	\label{eq:EOP_cont}
	\end{equation}
	Here, $\nabla \widehat{u}: \Omega \rightarrow \mathbb{R}^d$ is the gradient field of $\widehat{u}$, and the norm $\|\cdot\|_2$ induced by the inner product $\langle\cdot, \cdot\rangle$ is the standard 2-norm for functions over $\Omega$. Note that we assume $\phi \leq 0$ on $\partial \Omega$, so the boundary condition $\widehat{u}=0$ on $\partial \Omega$ is well-posed.
	An approximation of \cref{eq:EOP_cont} can be obtained by linearization \cite{Wu2015}, based on the Maclaurin series $\sqrt{1+x^2} = 1+\frac{1}{2} x^2+o\left(x^4\right)$. Assuming that $\|\nabla \widehat{u}\|^2_2$ is sufficiently small such that we can ignore higher order terms, we have $\sqrt{1+\|\nabla \widehat{u}\|_{2}^2} \approx 1+\frac{1}{2}\|\nabla \widehat{u}\|_{2}^2$; further dropping the constant $1$ gives the linearized obstacle problem
	\begin{equation}
		\min_{\widehat{u} \in \Omega} \quad \int_{\Omega} \frac{1}{2}\|\nabla \widehat{u}\|_{2}^2 \mathrm{~d}\mathbf{x} \quad \text{s.t. } \widehat{u} \geq \phi,~ \widehat{u}=0 \text { on } \partial \Omega. \label{eq:EOP}
	\end{equation}
    In order to solve the linearized problem \cref{eq:EOP} we discretize it with finite differences on the aforementioned equispaced meshes.
    Consider first the 1D case. Let $x_i, i \in \{1, \ldots, n\}$, denote an interior mesh point, and let $\widehat{u}_i \approx \widehat{u}( x_i )$ be an approximation of $v$ at this mesh point. Then the integral is discretized as
    $\int_{\Omega} \frac{1}{2}\|\nabla \widehat{u}(x)\|^2_2 \mathrm{d} x \approx \sum_{i = 1}^n \frac{1}{2}\left\|\nabla \widehat{u}(x_i)\right\|^2_2 h \approx \frac{h}{2} \left\langle A_1 \mathbf{\widehat{u}}, \mathbf{\widehat{u}} \right\rangle$, where $\mathbf{\widehat{u}} = ( \widehat{u}_1, \ldots, \widehat{u}_{n} ) \in \mathbb{R}^n$, and $A_1 \in \mathbb{R}^{n \times n}$ is the 1D Laplacian matrix \eqref{eq:1d_laplacian}, which assumes zero boundary conditions.
    The approximation proceeds analogously in 2D, where now we let $(x_i, y_j)$, $i,j \in \{1, \ldots, n\}$, denote an interior mesh point, and $\widehat{u}_{i,j} \approx \widehat{u}(x_i, y_j)$.
    In this case, the integral is approximated as $\int_{\Omega} \frac{1}{2}\|\nabla \widehat{u}(x, y)\|^2_2 \mathrm{d} x \mathrm{d} y \approx \sum_{i, j = 1}^n \frac{1}{2}\left\|\nabla \widehat{u}(x_i, y_j) \right\|^2_2 h^2 \approx\frac{h^2}{2}\left\langle A_2 \mathbf{\widehat{u}}, \mathbf{\widehat{u}}\right\rangle$, where $\mathbf{\widehat{u}} = (\widehat{u}_{1,1}, \ldots, \widehat{u}_{n,n}) \in \mathbb{R}^{n^2}$ and $A_2 \in \mathbb{R}^{n^2 \times n^2}$ is the 2D Laplacian matrix \eqref{eq:2d_laplacian}. 

\begin{figure}[h!]
		\centering
		\begin{subfigure}{.5\textwidth}
			\centering
			\includegraphics[scale=0.42]{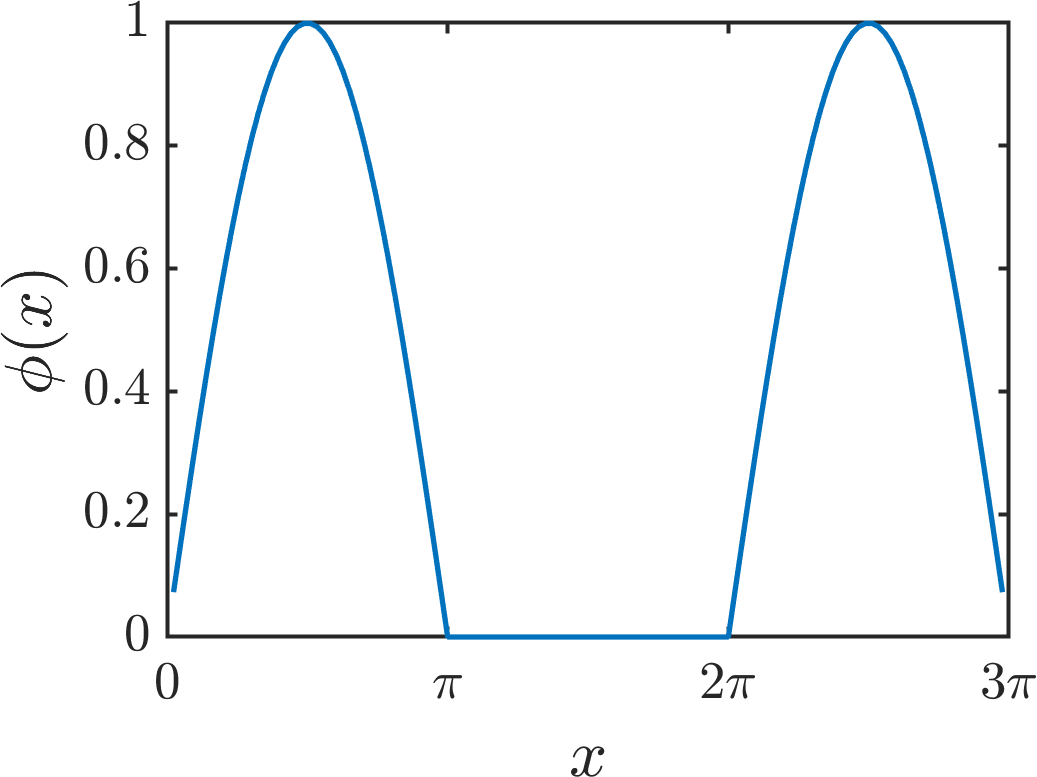}
			\caption{Obstacle $\phi$}
		\end{subfigure}%
		\begin{subfigure}{.5\textwidth}
			\centering
			\includegraphics[scale=0.42]{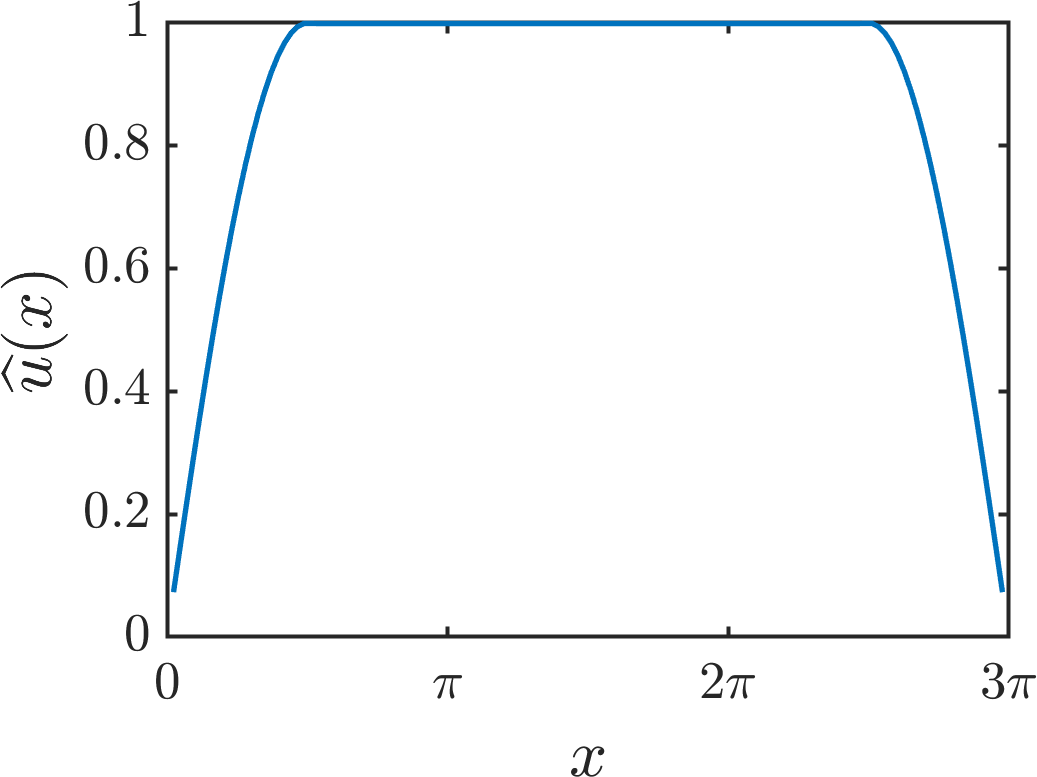}
			\caption{Solution $\widehat{u}$}
		\end{subfigure}
        \centering
		\begin{subfigure}{.5\textwidth}
			\centering
			\includegraphics[scale=0.42]{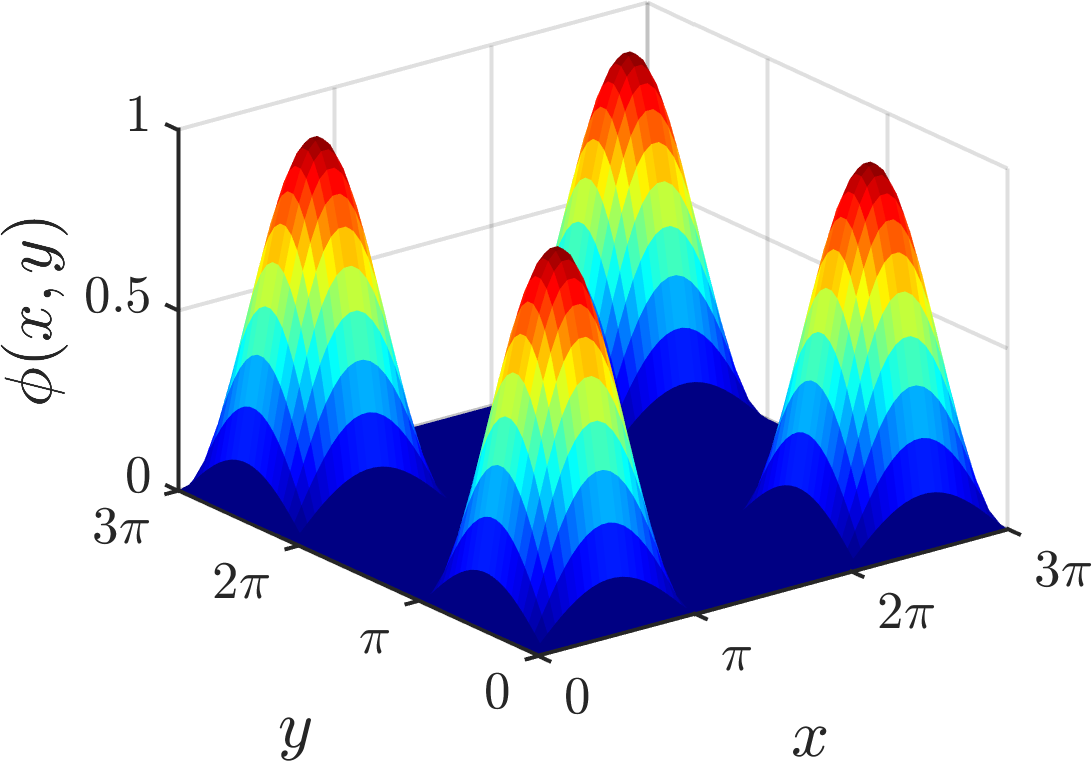}
			\caption{Obstacle $\phi$}
		\end{subfigure}%
		\begin{subfigure}{.5\textwidth}
			\centering
			\includegraphics[scale=0.42]{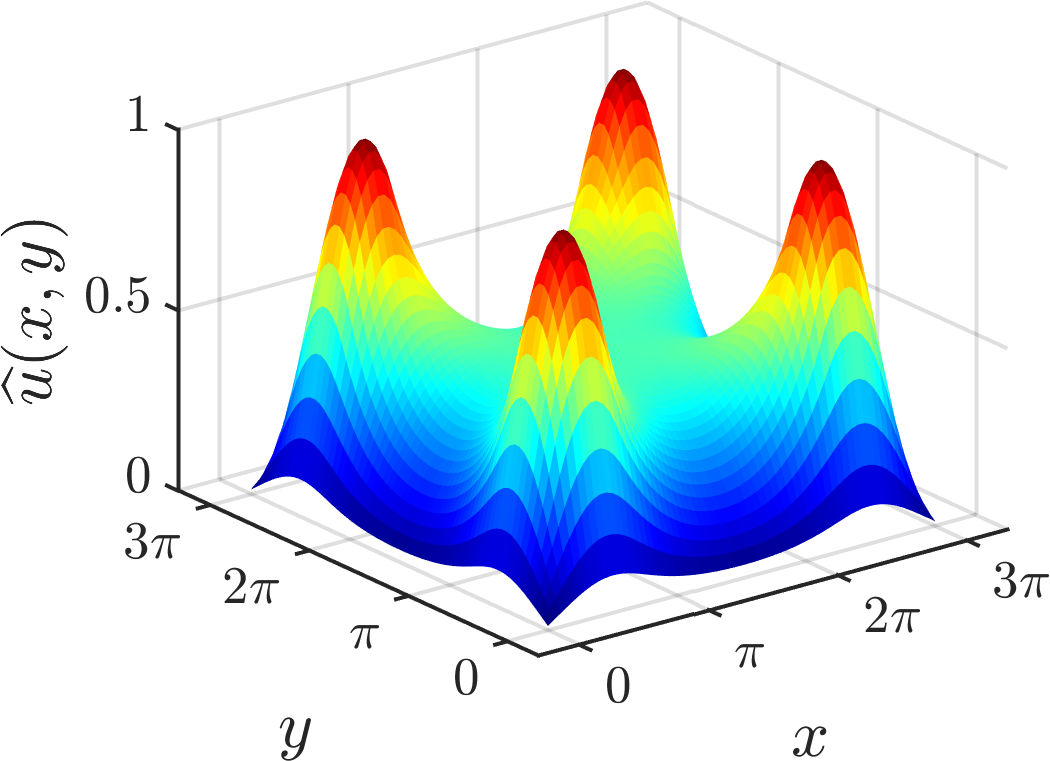}
			\caption{Solution $\widehat{u}$}
		\end{subfigure}
               \caption{(a) and (b): obstacle $\phi$ and solution $\widehat{u}$ of MP2-1D \cref{eq:mp2}. (c) and (d): obstacle $\phi$ and solution $\widehat{u}$ of   MP2-2D \cref{eq:mp2}, computed numerically using the proximal gradient algorithm \cref{eq:pg}.}
               \label{fig:eop_sol}
	\end{figure}

	Ignoring the leading $h$ or $h^2$ factors, since this does not change the solution of the minimization problem, our discretized versions of \cref{eq:EOP} take the form
	\begin{equation}\label{eq:aEOP}
		\underset{\mathbf{\widehat{u}} \in \mathbb{R}^N}{\min}\quad \frac{1}{2}\langle A_d \mathbf{\widehat{u}},\mathbf{\widehat{u}} \rangle \quad \text {s.t. } \mathbf{\widehat{u}} \geq \mathbf{\phi},
	\end{equation}
    where $N = n$ for $d = 1$, and $N = n^2$ for $d = 2$. Moreover, $\mathbf{\phi} \in \mathbb{R}^N$ denotes the discretized version of $\phi(\mathbf{x})$ on the underlying mesh, and $A_d$ the Laplacian matrix on this mesh. 
    By making a change of variable $\mathbf{u}:=\mathbf{\widehat{u}}-\mathbf{\phi} \geq \mathbf{0}$, we can write \cref{eq:aEOP} in its shifted form:
	\begin{equation}\label{eq:saEOP}
		\underset{\mathbf{u} \in \mathbb{R}^N}{\min}\quad \frac{1}{2}\langle A_d \mathbf{u}, \mathbf{u}\rangle +\langle A_d\mathbf{\phi}, \mathbf{u}\rangle + \frac{1}{2}\langle A_d\mathbf{\phi},\mathbf{\phi} \rangle \quad \text {s.t. } \mathbf{u} \geq \mathbf{0}. \end{equation}
    Dropping the $\frac{1}{2}\langle A_d \mathbf{\phi},\mathbf{\phi} \rangle$ term  (since it doesn't influence the minimizer), choosing $\mathbf{p}=-A_d\mathbf{\phi}$ and using an indicator function $i_{+} : \mathbb{R}^N \to \mathbb{R}$, where $i_{+}( \mathbf{a})=0$ if $(\mathbf{a})_q \geq 0$ for all $q$ and $i_{+}(\mathbf{a})=+\infty$ if $(\mathbf{a})_q<0$ for any $q$, allows \cref{eq:saEOP} to be written in unconstrained form as 
	\begin{equation}\label{eq:saEOP_i}
			\underset{\mathbf{u} \in \mathbb{R}^N}{\min}\quad \frac{1}{2}\langle A_d \mathbf{u}, \mathbf{u}\rangle-\langle \mathbf{p}, \mathbf{u} \rangle+i_{+}(\mathbf{u}).
	\end{equation}
    However, we choose to work with an alternative unconstrained formulation, given by
	\begin{equation}
		\underset{\mathbf{u} \in \mathbb{R}^N}{\min} \quad \frac{1}{2}\langle A_d \mathbf{u}, \mathbf{u}\rangle-\langle \mathbf{p}, \mathbf{u}\rangle+\lambda\left\|(-\mathbf{u})_{+}\right\|_1. \label{eq:mp2}
	\end{equation}
	where $\lambda>0$ is a predefined penalty parameter and $\sum_i\left(-(\mathbf{u})_i\right){_{+}}=\left\|(-\mathbf{u})_{+}\right\|_1$, where $(a)_{+}=\max \{0, a\}$ and $\|\cdot\|_1$ denotes the $L^1$ norm. For sufficiently large $\lambda$, the solution of \cref{eq:mp2} is the same as that of \cref{eq:saEOP_i} \cite{Mangasarian1985, Friedlander2008}, hence approximating the solution of \cref{eq:EOP} after reversing the change of variables $\mathbf{\widehat{u}}=\mathbf{u}+\mathbf{\phi}$.
    Throughout the rest of this paper, we refer to the shifted approximate obstacle problem in penalty form \cref{eq:mp2} as MP2-1D and MP2-2D, in the 1D and 2D cases, respectively. In our MP2-1D numerical experiments, we consider $x \in \Omega:=[0,3 \pi]$ and the obstacle $\phi(x)=\max \{0, \sin x\}$, for which \cref{eq:EOP} has the exact solution given by $\widehat{u}(x)=\sin x$ for $x \in\left[0, \frac{1}{2} \pi\right]$ or $x \in\left[\frac{5}{2} \pi, 3 \pi\right]$ and $\widehat{u}(x)=1$ for $x \in\left[\frac{1}{2} \pi, \frac{5}{2} \pi\right]$; see \cref{fig:eop_sol}. For MP2-2D, we consider $(x, y) \in \Omega:=[0,3 \pi]^2$ and $\phi(x, y)=\max \{0, \sin x\} \max \{0, \sin y\}$; see \cref{fig:eop_sol}.
    We take $\lambda=900$ in both the 1D and 2D cases, which is large enough such that solving \cref{eq:mp2} is equivalent to solving \cref{eq:saEOP_i}. 
    
    \cref{fig:eop_sol} presents visualizations of the obstacles and corresponding solutions of MP2-1D and MP2-2D.
    Note that the solutions pictured in panels (b) and (d) are those on the interior of the domain, and hence with the boundary points $\widehat{u} = 0$ eliminated.

    \section{Optimization algorithms}\label{sec:opt_algos}
	
	 \ptxt{The algorithms we employ are based on gradients and proximals of functions; in the case of nondifferentiable convex functions, gradients are replaced by subgradients. Let us recall the definition of the proximal operator $P_f$ of a function $f:\mathbb{R}^{N} \rightarrow \mathbb{R}$:
	\begin{equation}\label{eq:prox}
    P_f(\mathbf{v})=\left\{\mathbf{u} \in \mathbb{R}^N:\mathbf{u}=\underset{\mathbf{z} \in \mathbb{R}^{N}}{\arg\min}\left(f(\mathbf{z})+\frac{1}{2}\|\mathbf{z}-\mathbf{v}\|^2\right)\right\}.
	\end{equation}
    Under the assumption that $f$ is convex and continuous, which holds throughout this paper,  $\left(f(\cdot)+\frac{1}{2}\|\cdot-\mathbf{v}\|^2\right)(\mathbf{z})$ is strongly convex, such that $P_f(\mathbf{v})$ is single-valued.
    Recall the definition of $\partial f$, the subdifferential of $f$:
    \begin{equation}\label{eq:subdiff}
    \partial f(\mathbf{u})=\{\mathbf{z} \in \mathbb{R}^N: \forall \mathbf{v} \in \mathbb{R}^N,~ f(\mathbf{u})+\langle \mathbf{v}-\mathbf{u}, \mathbf{z}\rangle \leq f(\mathbf{v})\},
     \end{equation}
    where we call $\mathbf{z}\in \partial f(\mathbf{u})$ a subgradient of $f$ at $\mathbf{u}$; if $f$ is differentiable at $\mathbf{u}$, then $\partial f(\mathbf{u})=\{\nabla f(\mathbf{u})\}$. 
    We note that $\mathbf{u}=P_f(\mathbf{v})$ if and only if $\mathbf{v}-\mathbf{u} \in \partial f(\mathbf{u})$ or, equivalently, $\mathbf{v} \in \mathbf{u}+\partial f (\mathbf{u})=(I+\partial f)(\mathbf{u})$.
    Then, it follows that $\mathbf{u}=(I+\partial f)^{-1}(\mathbf{v})$, such that we can characterize the proximal operator as the resolvent of the subdifferential operator:
    \begin{equation}\label{eq:prox_def2}
     P_f=(I+\partial f)^{-1}.   
    \end{equation}
    
    We now discuss the optimization algorithms that we apply to the discrete model problems from \cref{sec:mp}. The main algorithms we consider are gradient descent for MP1, \cref{eq:mp1}, and proximal gradient descent for MP2, \cref{eq:mp2}, whose convergence we attempt to speed up by means of parallelization-in-time. In our method, we also make use of two auxiliary algorithms \ptxt{on the coarse grid}, namely proximal point for MP1 and alternating proximal mappings for MP2, which are closely related to gradient descent and proximal gradient descent, respectively.} We refer the reader to \cite[Chapters 26, 29-31]{Bauschke2023} for further details on the convergence properties of these algorithms.
		
	\subsection{Optimization algorithms for MP1}\label{subsec:mp1_algos}	
	Consider the problem of minimizing a convex and $L$-smooth function $f$,
    \begin{equation}\label{eq:min_f}
        \min_{\mathbf{u} \in \mathbb{R}^N} \quad f(\mathbf{u}).
    \end{equation}
	The first of two methods we consider for solving \cref{eq:min_f} is the classical gradient descent method. Given a starting point $\mathbf{u}_0 \in \mathbb{R}^{N}$ and a fixed positive scalar $s > 0$, the gradient descent iteration is given by
	\begin{equation}
		\mathbf{u}_{k+1} = (I-s \nabla f)(\mathbf{u}_{k}) =: G_{sf}(\mathbf{u}_k), \quad k=0, 1, \ldots.	 \label{eq:gd}
	\end{equation}
    The gradient descent sequence converges to a minimizer of function $f$ for $0 < s < \frac{2}{L}$ with a convergence rate of $\mathcal{O}(1/k)$; this means that if we want $|f(\mathbf{u}_k)-f(\mathbf{u}^*)| \leq \epsilon$, where $\mathbf{u}^*$ is a minimizer of $f$ and $\epsilon>0$, we need $\mathcal{O}(\frac{1}{\epsilon})$ iterations of gradient descent. Under strong convexity, linear convergence is attained \cite{Nesterov2004} and  $\mathcal{O}\left(\log\left(\frac{1}{\epsilon}\right)\right)$ are necessary asymptotically.

    We note that MP1 \cref{eq:mp1} is a case of \cref{eq:min_f} for the strongly convex function $f(\mathbf{u})=\langle A_1\mathbf{u}, \mathbf{u}  \rangle - \langle \mathbf{b}, \mathbf{u} \rangle$, with $A_1$ the Laplacian \eqref{eq:1d_laplacian}.
    Moreover, $\partial f(\mathbf{u})=\{\nabla f(\mathbf{u)}\}$, $\nabla f(\mathbf{u})=A_1\mathbf{u} - \mathbf{b}$ is $L$-Lipschitz with $L=\|A_1\|_2$, such that $f$ is $L$-smooth. Therefore, the objective function $f$ defined by \cref{eq:mp1} satisfies the conditions required for convergence of both \cref{eq:gd} and \cref{eq:ppm}.

    \ptxt{When it is necessary to take steps of size larger than the $\frac{2}{L}$ limit for gradient descent, we substitute gradient descent steps by proximal point steps, which is stable for all $s > 0$.} Given some starting point $\mathbf{u}_0$ and a positive scalar $s > 0$, the proximal point algorithm update for a continuous and convex $f$ is given by
	\begin{equation}\label{eq:ppm}
		\mathbf{u}_{k+1}=(I+s\partial f)^{-1}\left(\mathbf{u}_k\right):=P_{s f}\left(\mathbf{u}_k\right) \quad k=0, 1, \ldots,	 
	\end{equation}
	where we note that $P_{sf}$ is the proximal of $sf$ \cref{eq:prox}; since $s>0$, we have $s\partial f(\cdot)=\partial (sf)(\cdot)$.
    The proximal point iteration \cref{eq:ppm} converges to a minimizer of $f$; note that $L$-smoothness of $f$ is not required for this convergence, unlike in the case of gradient descent. \ptxt{We note that for simple problems, such as MP1, the proximal point method can converge in one step. However, our interest here is not in the proximal point method itself, but in using proximal steps to approximate a fixed number of gradient descent steps, as discussed in \cref{subsec:phi_delta_choice}.}
    
	\subsection{Optimization algorithms for MP2}\label{subsec:mp2_algos}
	
	Consider the possibly nonsmooth problem of minimizing the sum of functions 	\begin{equation}\label{eq:min_f+g}
			\min_{\mathbf{u} \in \mathbb{R}^{N}} \quad F(\mathbf{u}), \quad \text{where }F(\mathbf{u})=f(\mathbf{u})+g(\mathbf{u}),
	\end{equation}
	where $f$ is convex and $L$-smooth, and $g$ is continuous, convex and possibly nondifferentiable. We also assume that the set of minimizers of $F$ is nonempty. Minimizers of $F$ can be computed by splitting algorithms such as the proximal gradient method, which takes the form
	\begin{equation} \label{eq:pg}
		\mathbf{u}_{k+1}:=P_{sg}\circ G_{sf}(\mathbf{u}_{k})=(I+s \partial g)^{-1}\left((I-s \nabla f)(\mathbf{u}_{k})\right), \quad k=0, 1, \ldots ,	
	\end{equation}
	where $0 < s < \frac{2}{L}$ , where $P_{sg}$ is the proximal of $sg$ \cref{eq:prox} and $G_{sf}$ is the gradient descent step for $f$ of size $s$ \cref{eq:gd}. 
    The proximal gradient method is a composition of a gradient descent step over $f$ and a proximal point step over $g$; we note that proximal gradient inherits the stability constraint of gradient descent that $s<\frac{2}{L}$.  
    We also note that for $f=0$ in \eqref{eq:min_f+g}, iteration \eqref{eq:pg} reduces to the proximal point algorithm \cref{eq:ppm}, while for $g=0$ in \eqref{eq:min_f+g} the iteration reduces to gradient descent \cref{eq:gd}.

    For MP2 \cref{eq:mp2} written in the form of \eqref{eq:min_f+g}, we have $f=\frac{1}{2}\langle A_d \mathbf{u}, \mathbf{u}\rangle-\langle\mathbf{p}, \mathbf{u}\rangle$ and $g=\lambda\left\|(-\mathbf{u})_{+}\right\|_1$. 
    Here, $f$ is strongly convex and $L$-smooth, with $L=\|A_d\|_2$, where $A_d$ is given by \cref{eq:1d_laplacian} for $d=1$ and by \cref{eq:2d_laplacian} for $d=2$, and $g$ is convex and continuous. Therefore, $f$ and $g$ satisfy the requirements for convergence of \cref{eq:pg}, which converges to a solution of \cref{eq:min_f+g} as desired.

    \ptxt{When it is necessary to take steps of size $s>\frac{2}{L}$, we substitute proximal gradient descent by the alternating proximal mappings method, with iteration given by} 
	\begin{equation}
		\mathbf{u}_{k+1}:=P_{sg}\circ P_{sf}\left(\mathbf{u}_k\right)=(I+s\partial g)^{-1}\left((I+s \nabla f)^{-1}\left(\mathbf{u}_k\right)\right), \label{eq:apm}
	\end{equation}
    where $P_{sg}$ and $P_{sf}$ are proximals \cref{eq:prox} of $sg$ and $sf$, respectively, and we note that $\partial f=\{\nabla f\}$ since $f$ is differentiable by assumption. This method \eqref{eq:apm} is a composition of two proximal steps over $f$ and $g$. In particular, the sequence \eqref{eq:apm} converges to a minimizer of $\operatorname{env} (f)+g$ \cite[Proposition 31.1]{Bauschke2023}, where 
    \begin{equation}
             \operatorname{env} (f(\mathbf{u})) := \underset{\mathbf{u} \in \mathbb{R}^N}{\min}\quad \left(f\left(P_{f}(\mathbf{u})\right)+\frac{1}{2}\left\|\mathbf{u}-P_{f}(\mathbf{u})\right\|^2\right), \label{eq:min_envf+g}
    \end{equation}
      is the Moreau envelope of $f$ at $\mathbf{u}$. This means that \cref{eq:apm} does not generally converge to a minimizer of \cref{eq:min_f+g}.

    \subsection{Optimization methods as ODE discretizations}\label{subsec:opt_disc_ode}
	
	Here, we describe the connection between optimization algorithms and ODE discretizations, which provides key motivation for the parallel-in-iteration framework we outline in the next section. 
    The connection between these has been well established in the literature; see \cite{Parikh2014, Su2016,Attouch2016, Romero2022}.
    First, consider the gradient flow ODE of a convex and $L$-smooth function $f$, given by
	\begin{equation}
		\frac{\mathrm{d}}{\mathrm{d}t}\mathbf{u}(t)=-\nabla f(\mathbf{u}(t)), \quad \mathbf{u}(0)=\mathbf{u}_0. \label{eq:gf}
	\end{equation}
    Applying a forward finite-difference approximation of this equation yields the forward-Euler discretization 
    \begin{equation} \label{eq:fe_gf}
		\mathbf{u}_{k+1}=\mathbf{u}_{k}-s\nabla f(\mathbf{u}_{k}),
        \quad
        k = 0, 1, \ldots,
	\end{equation}
    where $\mathbf{u}_{k+1} \approx \mathbf{u}(t_k+\Delta t)$, $\mathbf{u}_{k} \approx \mathbf{u}(t_k)$, and $s=\Delta t$.
    Observe that this forward Euler discretization is identical to the gradient descent iteration \cref{eq:gd}. 
    
    Similarly, applying a backward finite-difference approximation to \eqref{eq:gf} yields the backward Euler discretization
      \begin{equation} \label{eq:be_gf}
          \mathbf{u}_{k+1}=\mathbf{u}_{k} -s\nabla f(\mathbf{u}_{k+1}),
          \quad k = 0, 1, \ldots
      \end{equation}
	Let us recall the proximal point method \cref{eq:ppm} and note that
	\begin{equation}
			\begin{aligned}
            \mathbf{u}_{k+1}= P_{sf}(\mathbf{u}) &\iff 
			\mathbf{u}_{k+1}=\underset{\mathbf{u}}{\operatorname{argmin} }\left(sf(\mathbf{u})+\frac{1}{2}\left\|\mathbf{u}-\mathbf{u}_{k}\right\|^2\right) \\
			& \iff \mathbf{u}_{k+1}=\underset{\mathbf{u}}{\operatorname{argmin} }\left(f(\mathbf{u})+\frac{1}{2s}\left\|\mathbf{u}-\mathbf{u}_{k}\right\|^2\right) \\
			& \iff \mathbf{0}=\nabla f\left(\mathbf{u}_{k+1}\right)+\frac{1}{s}\left(\mathbf{u}_{k+1}-\mathbf{u}_{k}\right) \\
			& \iff \mathbf{u}_{k+1}=\mathbf{u}_{k} -s\nabla f(\mathbf{u}_{k+1}).
		\end{aligned}
	\end{equation}
	Therefore, the proximal point algorithm \cref{eq:ppm} is equivalent to a backward Euler discretization of \cref{eq:gf}. Furthermore, the proximal point algorithm can be interpreted as an implicit variant of the gradient descent method \cref{eq:gd}. 
    In the context of ODEs and PDEs, the forward Euler method as in \eqref{eq:fe_gf} is known to have poor stability properties for stiff equations, putting strict requirements on the step size $\Delta t$. On the other hand, the backward Euler discretization as in \eqref{eq:be_gf} has relatively improved stability properties for stiff equations, typically permitting much larger time step sizes.
    These stability properties influence the design of our parallel-in-iteration methodology described in the next section.
    In particular, the methodology requires us to take large steps on coarse grids, and we chose to use the proximal point method \cref{eq:ppm} when these steps are larger than $s>\frac{2}{L}$ so as to avoid instability of gradient descent \eqref{eq:gd}.
	
	Now, consider the gradient flow of $f+g$, assuming, for simplicity, that both $f$ and $g$ are differentiable. The associated gradient flow ODE is
	\begin{equation}
		\frac{\mathrm{d}}{\mathrm{d}t}\mathbf{u}(t)=-\nabla f(\mathbf{u}(t))-\nabla g(\mathbf{u}(t)). \label{eq:gf2}
	\end{equation}
    Applying a forward Euler step on the $\nabla f$ term, and a backward Euler step on the $\nabla g$ term results in the implicit-explicit (IMEX) discretization
    \begin{equation}\label{eq:imex}
		\mathbf{u}_{k+1}=(I+s \nabla g)^{-1}\left((I-s \nabla f) (\mathbf{u}_{k})\right), \quad k = 0, 1, \ldots        
	\end{equation}
    Observe that \cref{eq:imex} is identical to \cref{eq:pg}, such that the proximal gradient method can be interpreted as an IMEX discretization of the corresponding gradient flow equation \cref{eq:gf2}. 
    Similarly, the alternating proximal mappings methods \cref{eq:apm} is equivalent to a split fully backward-Euler-based discretization of \cref{eq:gf2}.

	\section{Parallel-in-time methods}\label{sec:p-in-t}

    In this section we outline our parallel-in-iteration strategy for the solution of the discrete optimization problems outlined in the previous section.
    As outlined in \cref{sec:intro}, the motivation for parallelizing over iteration is that convergence rates of these methods tend to be slow, and thus, in practice, they typically require many iterations to reach convergence.
    By re-purposing techniques for the parallel-in-time solution of ODEs and PDEs to parallelize over iterations of optimization algorithms, we aim to reduce the wall-clock time required to solve discrete optimization problems.
    Specifically in this work we focus on the parallel-in-time method known as MGRIT (multigrid reduction in time) \cite{Falgout2014}. 
    However, we note that other suitable parallel-in-time methods can also be considered. 
    For example, the parareal method \cite{Lions2001} or the Parallel Full Approximation Scheme (FAS) in Space-Time (PFASST) method \cite{Minion2010, Emmett2012}.
    %


	\subsection{The MGRIT algorithm}\label{subsec:mgrit}

\begin{figure}[h!]
            \begin{center}
		\includegraphics[scale=1]{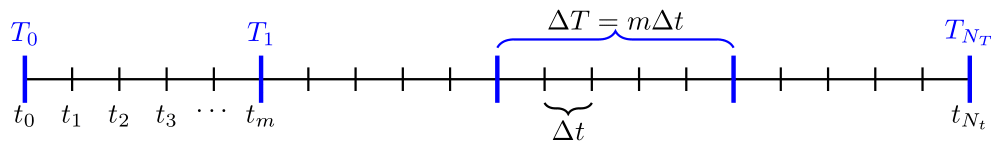}
		\caption{Fine- and coarse-grid discretization meshes with coarsening factor $m$. F-points (in black) are points $t_i$ for $i=1,\ldots,N_t$ such that $i \not = jm$ for $j=1,\ldots,N_T=\frac{N_t}{m}$, and C-points (in blue) are points $T_j=t_{jm}$ for $j=1,\ldots,N_T=\frac{N_t}{m}$.}
		\label{fig:fine_coarse_grids}
            \end{center}
	\end{figure}

    MGRIT \cite{Falgout2014} was developed as a multilevel generalization of the 2-level parareal method \cite{Lions2001} for solving time-dependent systems of differential equations in parallel.
    To outline this method, let us consider a general initial-value ODE problem of the form
	\begin{equation}\label{eq:cont_ode_system}
			\frac{\mathrm{d}}{\mathrm{d}t}\mathbf{u}(t)={\cal L}(t, \mathbf{u}(t)), \quad \mathbf{u}(0)=\mathbf{u}_0, \quad t \in (0, T].
	\end{equation}
    Discretizing this system of ODEs using a 1-step method (e.g., any Runge-Kutta method) results in a system of equations of the form
	\begin{equation}\label{eq:mgrit_problem}
		\begin{aligned}
			& \mathbf{u}_0=\mathbf{w}_0, \\
			& \mathbf{u}_i=\Phi_i(\mathbf{u}_{i-1})+\mathbf{w}_i, \quad i=1,2, \ldots, N_t,
		\end{aligned}
	\end{equation}
	where $\mathbf{w}_0 \equiv \mathbf{u}_0$, and $t_i=i \Delta t$, $i= 0,1, \ldots, N_t$ is a temporal mesh with spacing $\Delta t=\frac{T}{N_t}$, which we assume is constant for simplicity. 
    Furthermore, $\mathbf{u}_i \approx\mathbf{u}(t_i)$ is the discrete approximation to $\mathbf{u}(t)$ at time $t_i$, and $\Phi_i$ is the \emph{time-stepping operator}, responsible for advancing the solution at $t_{i-1}$ to $t_i$; for simplicity, we further assume $\Phi_i$ is a time-independent operator, such that $\Phi_i(\cdot)=\Phi(\cdot)$ for all $i=1,2,\ldots,N_t$.
    The vector $\mathbf{w}_i$ represents solution-independent forcing terms in the discretized problem. 
    We can then rewrite \cref{eq:mgrit_problem} as the following ``all-at-once'' system: 
    	\begin{equation} \label{eq:A(u)=g}
		\mathcal{A}(\mathbf{u})=\left[\begin{array}{cccc}
			I & & & \\
			-\Phi(\cdot) & I & & \\
			& \ddots & \ddots & \\
			& & -\Phi(\cdot) & I
		\end{array}\right]\left[\begin{array}{c}
			\mathbf{u}_0 \\
			\mathbf{u}_1 \\
			\vdots \\
			\mathbf{u}_{N_t}
		\end{array}\right]=\left[\begin{array}{c}
			\mathbf{w}_0 \\
			\mathbf{w}_1 \\
			\vdots \\
			\mathbf{w}_{N_t}
		\end{array}\right] = \mathbf{w}.
	\end{equation}
    For our purposes, it is also pertinent to introduce the associated residual vector, which, given some $\mathbf{v} \approx \mathbf{u}$, is defined as
    \begin{align}\label{eq:mgrit_res}
        \mathbf{r}(\mathbf{v})
        =
        \mathbf{w}
        -
        {\cal A}( \mathbf{v} )=
        \begin{bmatrix}
            \mathbf{w}_0 - \mathbf{v}_0 \\
            \mathbf{w}_1 - [ \mathbf{v}_1 - \Phi( \mathbf{v}_0 )] \\
            \vdots \\
            \mathbf{w}_{N_t} - [ \mathbf{v}_{N_t} - \Phi( \mathbf{v}_{N_t-1} )]
        \end{bmatrix}.
    \end{align}
	Instead of solving \cref{eq:mgrit_problem} sequentially with time-stepping, MGRIT aims to solve it iteratively and in parallel via the introduction of a coarser approximation of the problem. 
    To this end, given an integer coarsening factor $m>1$, we introduce a coarse temporal mesh $T_j=j \Delta T$, $j=0,1, \ldots, N_T$, where $N_T=\frac{N_t}{m}$, with constant spacing $\Delta T=m \Delta t$. This partitions the $N_t$ points of the fine grid into two sets: 1. \emph{C-points}, which are the $N_T$ points appearing exclusively on the coarse grid, and 2. \emph{F-points} which are the remaining $N_t - N_T$ points that do not appear on the coarse grid. 
    \cref{fig:fine_coarse_grids} illustrates the fine- and coarse-grid temporal meshes.
	
    As is standard for a multilevel method, an MGRIT iteration consists of a coarse-grid correction step sandwiched between pre- and post-relaxation steps.
    The standard relaxation, as used here, is known as FCF-relaxation \cite{Falgout2014}, and combines so-called F- and C-relaxations, both of which can be performed in parallel. 
    An F-relaxation is the process where, starting from a C-point, the time-stepping operator is used to update the values of the $m-1$ F-points that follow it (that is, until the F-point that immediately precedes the next C-point). 
	 C-relaxation is the update of each C-point by stepping from the preceding F-point using $\Phi$.
    %
    FCF-relaxation is simply the application of F-,C- and F-relaxations in this order. 
    The post relaxation is taken to be an F-relaxation.

    In between the pre- and post-relaxations, a coarse-grid correction problem is solved to provide a correction at C-points, as we now describe. 
    MGRIT has both linear and FAS variants \cite{Howse2019}, suitable for linear and nonlinear $\Phi$, respectively, and these variants differ in the way that coarse-grid correction is performed. 
    Let us first explain how this works in the linear case, and to this end let us assume that $\Phi_i$ is a time-independent, linear operator, such that $\Phi_i(\mathbf{u}_{i-1})=\Phi\mathbf{u}_{i-1}$ corresponds to a matrix-vector product with matrix $\Phi$.
    We also write the action of the space-time operator ${\cal A}( \mathbf{u} ) = {\cal A} \mathbf{u}$ as a matrix-vector product.
	From \eqref{eq:mgrit_problem}, we can then write, for coarse indices $j=1,2, \ldots, N_T$, 
	\begin{equation}
		\begin{aligned}
			\mathbf{u}_{j m}=\Phi \mathbf{u}_{j m-1}+\mathbf{w}_{j m}=\cdots=\Phi^m \mathbf{u}_{(j-1) m}+\widetilde{\mathbf{w}}_{j m},
		\end{aligned}
	\end{equation}
	where $
		\widetilde{\mathbf{w}}_{j m}=\mathbf{w}_{j m}+\Phi \mathbf{w}_{j m-1}+\cdots+\Phi^{m-1} \mathbf{w}_{(j-1) m+1}$.
	Therefore, the fine-grid system \eqref{eq:A(u)=g} is equivalent to the following system written only in terms of the C-points:
	\begin{equation} \label{eq:c_mgrit_problem}
			\begin{aligned}
			\mathbf{u}_0 & =\mathbf{w}_0, \\
			\mathbf{u}_{j m} & =\Phi^m \mathbf{u}_{(k-1) m}+\widetilde{\mathbf{w}}_{j m}, \quad j=1,2, \ldots, N_T,
		\end{aligned}
	\end{equation}
	where $\Phi^m$ plays the role of time-stepping operator for this coarse system.
    Or, in block-matrix form
    \begin{equation}\label{eq:A_Delta}
    \mathcal{A}_{\Delta} \mathbf{u}_{\Delta}=\left[\begin{array}{cccc}
			I & & & \\
			-\Phi^m & I & & \\
			& \ddots & \ddots & \\
			& & -\Phi^m & I
		\end{array}\right]\left[\begin{array}{c}
			\mathbf{u}_{0} \\
			\mathbf{u}_{m} \\
			\vdots \\
			\mathbf{u}_{N_T m}
		\end{array}\right]=\left[\begin{array}{c}
			\mathbf{w}_0 \\
			\widetilde{\mathbf{w}}_m \\
			\vdots \\
			\widetilde{\mathbf{w}}_{N_Tm}
		\end{array}\right]=\mathbf{w}_\Delta.
    \end{equation}
	While system \eqref{eq:A_Delta} has a factor of $m$ fewer equations than \eqref{eq:A(u)=g}, it is no less expensive to solve because (naively) computing the action of $\Phi^m$ is $m$ times as expensive as computing that of $\Phi$.
    So, instead, MGRIT solves an approximation to the \emph{ideal} coarse-grid system \eqref{eq:A_Delta} arising from  replacing the \emph{ideal} coarse-grid time-stepping operator $\Phi^m$ with an approximation:
	\begin{equation} \label{eq:phi_delta}
			\Phi_{\Delta} \approx \Phi^m. 
	\end{equation}
    That is, rather than inverting $\mathcal{A}_{\Delta}$ in \eqref{eq:A_Delta}, one inverts the approximation $\mathcal{B}_{\Delta}$ given by
	\begin{equation}\label{eq:B_Delta}
		\mathcal{B}_{\Delta}=\left[\begin{array}{cccc}
			I & & & \\
			-\Phi_{\Delta} & I & & \\
			& \ddots & \ddots & \\
			& & -\Phi_{\Delta} & I
		\end{array}\right].
	\end{equation}

    Let us now explain in more detail how MGRIT uses the approximation ${\cal B}_{\Delta} \approx {\cal A}_{\Delta}$.
	Consider the coarse system ${\cal A}_{\Delta} \mathbf{u}_{\Delta}=\mathbf{w}_{\Delta}$ \eqref{eq:A_Delta}. Let $\mathbf{v}_{\Delta} \approx \mathbf{u}_{\Delta}$ be an approximate solution of \eqref{eq:A_Delta} with error $\mathbf{e}_{\Delta}=\mathbf{u}_{\Delta}-\mathbf{v}_{\Delta}$. 
    Then, note that
	\begin{equation}\label{eq:c_residual}
		\begin{aligned}
			\mathcal{A}_{\Delta} \mathbf{e}_{\Delta} & =\mathcal{A}_{\Delta} \mathbf{u}_{\Delta}-\mathcal{A}_{\Delta} \mathbf{v}_{\Delta}=\mathbf{w}_{\Delta}-\mathcal{A} \mathbf{v}_{\Delta} =\mathbf{r}_{\Delta},
		\end{aligned}
	\end{equation}
	where $\mathbf{r}_{\Delta}$ denotes the coarse-grid residual. Since $\mathcal{A}_{\Delta}$ is invertible, it follows that $\mathbf{e}_{\Delta}=\mathcal{A}_{\Delta}^{-1} \mathbf{r}_{\Delta}$, such that
	\begin{equation}
			\begin{aligned}
			\mathbf{u}_{\Delta} =\mathbf{v}_{\Delta}+\mathbf{e}_{\Delta}=\mathbf{v}_{\Delta}+\mathcal{A}_{\Delta}^{-1} \mathbf{r}_{\Delta} .
		\end{aligned}
	\end{equation}
	Assuming that $\Phi_{\Delta} \approx \Phi^m$ such that $\mathcal{B}_{\Delta}^{-1} \approx \mathcal{A}_{\Delta}^{-1}$, then	
	\begin{equation}
			\mathbf{u}_{\Delta} \approx \mathbf{v}_{\Delta}+\mathcal{B}_{\Delta}^{-1} \mathbf{r}_{\Delta}.
	\end{equation}
	Let $\left\{\mathbf{u}_{\Delta}^{k}\right\}_{k=0,1,\ldots}$ be a sequence of approximate solutions generated by the $k=0,1,\ldots$ MGRIT iterations. Letting $\mathbf{u}_{\Delta}^{k+1}$ and $\mathbf{u}_{\Delta}^{k}$ play the roles of $\mathbf{u}_{\Delta}$ and $\mathbf{v}_{\Delta}$, respectively, where we think of $\mathbf{u}_{\Delta}^{k+1}$ as a better approximation of $\mathbf{u}_{\Delta}$ than $\mathbf{u}_{\Delta}^{k}$, we obtain the MGRIT residual correction scheme:
	\begin{equation}\mathbf{u}_{\Delta}^{k+1}=\mathbf{u}_{\Delta}^{k}+\mathcal{B}_{\Delta}^{-1}\left(\mathbf{w}_{\Delta}-\mathcal{A}_{\Delta} \mathbf{u}_{\Delta}^{k}\right),
    \quad k = 0, 1, \ldots
	\end{equation}
	Component-wise, this may be written with $\Phi$ and $\Phi_{\Delta}$ as
	\begin{equation}
		\mathbf{u}_{(j+1)m}^{k+1}=\Phi_{\Delta} \mathbf{u}_{jm}^{k+1}+\Phi^m \mathbf{u}_{jm}^{k}-\Phi_{\Delta} \mathbf{u}_{jm}^{k}+\mathbf{\widetilde{\mathbf{w}}}_{(j+1)m}, 
	\end{equation}
	where $\mathbf{u}_{0}^{k}=\mathbf{w}_{0}$.
    
    %

    Now let us consider the MGRIT coarse-grid correction in the more general FAS case, where $\Phi(\cdot)$ in \eqref{eq:A(u)=g} is not (necessarily) a linear operator. 
    Let $\mathbf{v} \approx \mathbf{u}$ be an approximate solution of $\mathcal{A}(\mathbf{u})=\mathbf{w}$ in \cref{eq:A(u)=g}.    
    After performing the FCF-relaxation step, the fine-grid approximation $\mathbf{v}$ and nonlinear residual $\mathbf{w}-\mathcal{A}(\mathbf{v})$ are injected to the coarse-grid:
	\begin{equation}
		\mathbf{v}_{\Delta}=R_{I}(\mathbf{v}), \quad \mathbf{r}_{\Delta}=R_{I}(\mathbf{w}-\mathcal{A}(\mathbf{v})),
	\end{equation}
    where $R_I$ stands for the injection operator.
    The \emph{ideal} FAS coarse-grid problem is then to solve the system
    \begin{equation}\label{eq:c_mgrit_problem_fas}		 \mathcal{A}_{\Delta}\left(\mathbf{u}_{\Delta}\right)=\mathcal{A}_{\Delta}\left(\mathbf{v}_{\Delta}\right)+\mathbf{r}_{\Delta}
    \end{equation}
    for the coarse-grid solution $\mathbf{u}_{\Delta}$, where similarly to \eqref{eq:A_Delta},
    \begin{equation} 
		\mathcal{A}_\Delta(\cdot)=\left[\begin{array}{cccc}
			I & & & \\
			-\Phi^m(\cdot) & I & & \\
			& \ddots & \ddots & \\
			& & -\Phi^m(\cdot) & I
		\end{array}\right].
	\end{equation}
    The coarse-grid system \eqref{eq:c_mgrit_problem_fas} results from the elimination of F-points variables from the residual equation ${\cal A}( \mathbf{u} ) - {\cal A}(\mathbf{v}) = \mathbf{r}( \mathbf{v} )$.
    Similarly to the linear case, the MGRIT iteration solves an approximation of \cref{eq:c_mgrit_problem_fas}, 
    \begin{equation} 
    \mathcal{B}_{\Delta}\left(\mathbf{u}_{\Delta}\right) = \mathcal{B}_{\Delta}\left(\mathbf{v}_{\Delta}\right)+\mathbf{r}_{\Delta},
    \end{equation} 
    where $\mathcal{B}_{\Delta}(\cdot) \approx \mathcal{A}_\Delta(\cdot)$ replaces the $\Phi^m(\cdot)$ terms in $\mathcal{A}_\Delta(\cdot)$ by a computationally cheaper $\Phi_\Delta(\cdot) \approx \Phi^m(\cdot)$. 
    Letting $\mathbf{u}^{k+1}_\Delta$ and $\mathbf{u}^{k}_\Delta$ play the roles of $\mathbf{u}_\Delta$ and $\mathbf{v}_\Delta$, respectively, the component-wise FAS residual correction scheme is then given by 
	\begin{equation}\mathbf{u}^{k+1}_{(j+1)m}=\Phi_{\Delta}\left(\mathbf{u}^{k+1}_{jm}\right)+\Phi\left(\mathbf{u}^{k}_{(j+1)m-1}\right)-\Phi_{\Delta}\left(\mathbf{u}^{k}_{jm}\right)+\mathbf{w}_{(j+1)m}, \quad j=1,2, \ldots, N_T.
	\end{equation}
    We note that when $\Phi(\cdot)$ and $\Phi_\Delta(\cdot)$ are linear operators $\Phi$ and $\Phi_\Delta$, the FAS variant simplifies to the linear variant of MGRIT. 

    While we have only described 2-level MGRIT methods, multilevel MGRIT schemes can be obtained by recursively applying the 2-level schemes. In our numerical tests, we apply the standard MGRIT scheme for MP1 \cref{eq:mp1} and the FAS-MGRIT scheme for MP2 \cref{eq:mp2}. For both MP1 and MP2, MGRIT is iterated until the 2-norm of the residual \cref{eq:mgrit_res} is reduced by 8 orders of magnitude from its initial value. We note that since the residual \cref{eq:mgrit_res} is measured after performing the FCF-relaxations, the residual at the F-points will be zero, such that it is sufficient to measure the residual at the C-points.
        	
	\section{MGRIT for solving optimization problems parallel-in-iteration}\label{sec:mgrit_opt}

    In this section, we discuss our methodology for using MGRIT to solve optimization problems in a novel parallel-in-iteration framework, using our first model problem, MP1 \cref{eq:mp1}, as an example. First, we show the equivalence between solving MP1 using the gradient descent algorithm \cref{eq:gd} and solving the 1D heat equation, in terms of the all-at-once system in \cref{eq:A(u)=g}. Inspired by well-known established practices for solving ODEs, we then present our heuristic for choosing a suitable coarse-grid operator for MGRIT applied to gradient descent optimization of MP1. Finally, we present theoretical results relating the convergence of approximate solutions obtained for the all-at-once system solved by MGRIT and the convergence of the original, sequential optimization algorithm.

	\subsection{Links between solving ODEs using MGRIT and parallel-in-iteration optimization}\label{subsec:mgrit_ode_link}

    Recall that the $i$th gradient descent iteration \cref{eq:gd} for MP1, \cref{eq:mp1}, is given by
	\begin{equation}\label{eq:mp1_fine_grid_system}
	\begin{aligned}
		\mathbf{u}_{i}&=(I-s\nabla f)(\mathbf{u}_{i-1})=(I-sA_1)\mathbf{u}_{i-1}+s\mathbf{b}.
	\end{aligned}
	\end{equation}
	Therefore, denoting $\Phi:=(I-sA_1)$ and $\mathbf{w}_i = s\mathbf{b}$ for $i=1, \ldots, N_t$, we can write $N_t$ gradient descent iterations as
	\begin{equation}\label{eq:Au=g_he}
		\mathcal{A} \mathbf{u}=\left[\begin{array}{cccc}
		I & & & \\
		-\Phi & I & & \\
		& \ddots & \ddots & \\
		& & -\Phi & I
	\end{array}\right]\left[\begin{array}{c}
		\mathbf{u}_0 \\
		\mathbf{u}_1 \\
		\vdots \\
		\mathbf{u}_{N_t}
	\end{array}\right]=\left[\begin{array}{c}
		\mathbf{w}_0 \\
		\mathbf{w}_1 \\
		\vdots \\
		\mathbf{w}_{N_t}
	\end{array}\right]=\mathbf{w}.
	\end{equation}
	where $\mathbf{u}_0=\mathbf{w}_0$ is the initial condition. Note that this is the same kind of system as in \cref{eq:A(u)=g}; in particular, it is the same system one would obtain from the forward Euler discretization of the 1D heat equation with Dirichlet boundary conditions and constant source term:
	\begin{equation}\label{eq:1d_heat_eq}		
		\begin{aligned}
			\frac{\partial}{\partial t}u(x,t) - \frac{\partial^2}{\partial x^2}u(x,t)& = b, & (x, t) & \in \Omega \times[0, T], \\
			u(x, 0) & =u_0(x), & x & \in \Omega,
		\end{aligned}
	\end{equation}
    where $\Omega \subset \mathbb{R}$.

    
    While previous research on solving \cref{eq:1d_heat_eq} with MGRIT, namely \cite{Falgout2014, Dobrev2017}, certainly informs our investigation, differences in objective and approach justify a new exploration of this model problem, regardless of the similarities.
    First, both \cite{Falgout2014, Dobrev2017} are strongly concerned with matching the accuracy of spatial and temporal discretizations of \cref{eq:1d_heat_eq}, with \cite{Falgout2014} testing different ratios between time step lengths and spatial mesh sizes and  \cite{Dobrev2017} employing higher-order time discretizations in order to match the spatial and temporal discretization orders. Secondly, with this accuracy concern in mind, \cite{Falgout2014, Dobrev2017} mostly employ implicit discretizations, which are unconditionally stable in regards to time step lengths, on all levels. 
    These objectives differ from ours significantly, as we are concerned with providing speedup over sequential solves for problems to which explicit optimization methods are typically applied, and in which the time step lengths follow a standard choice for the chosen method. 
    Furthermore, the aforementioned accuracy concerns in both space and time are not present here, as our accuracy requirement is that the MGRIT solution at the final time points  approximates the solution of the model problem to a similar tolerance accuracy as the sequential solve which we are trying to speed up. 
    Therefore, it is necessary for us to provide our own numerical experiments and speedup estimates for this model problem, despite the similarities with \cite{Falgout2014, Dobrev2017}. 
    
	
	\subsection{Choice of $\Phi_{\Delta}$ for MP1 using gradient descent}\label{subsec:phi_delta_choice}

	Inspired by the similarities between our model problem \cref{eq:mp1} and the 1D heat equation \cref{eq:1d_heat_eq}, in this paper we follow a simple heuristic inspired by \cite{Falgout2014} of simply replacing forward Euler steps on the finest level with backward Euler steps on the coarse s.
    Recall from \cref{eq:mp1_fine_grid_system} that $\Phi:=I-sA_1$. We further refine this by choosing the step size to be $s=\frac{1}{L}$, which is the classical choice for gradient descent. As discussed in \cref{subsec:opt_disc_ode}, the natural implicit candidate for the coarse-grid operator is the one in \cref{eq:ppm}:
	\begin{equation}\label{eq:mp1_phi_delta}
			\Phi_{\Delta}:=(I+msA_1)^{-1}.
	\end{equation}
	That is, we use the approximation
	\begin{equation}\label{eq:mp1_phi_delta_approx}
			\Phi_{\Delta}:=(I+msA_1)^{-1} \approx (I-sA_1)^m = \Phi^m.
	\end{equation}
	Intuitively, we are approximating  $m$ forward Euler time steps, each with step size $s$, with a single backward Euler time step of size $ms$. We recall that the proximal point method \cref{eq:ppm} converges for any $s>0$, such that $ms=\frac{m}{L}$ is a feasible step size choice for any coarsening factor $m$. In a multilevel MGRIT setting of $\widehat{\ell} \geq 2$ levels, we take $\ell=1$ to refer to the fine grid. For $2 \leq \ell \leq \widehat{\ell}$, the $\ell$th level employs the operator 
    \begin{equation}\label{eq:mp1_phi_delta_l}
     \Phi_{\Delta}^{\ell}:=(I+m^{\ell-1}sA_1)^{-1},   
    \end{equation}
    an implicit operator similar to \cref{eq:mp1_phi_delta} with the appropriate time step coarsening; when $\ell=2$, we omit the superscript. \ptxt{We note that the computation of the inverse matrices in \cref{eq:mp1_phi_delta} and \cref{eq:mp1_phi_delta_l} amounts to solving a linear system, which we do by means of LU-factorization; the L and U factors are stored, such that the factorization is computed only once per level.}
\ptxt{\subsection{Choice of $\Phi_{\Delta}$ for MP2 using proximal gradient descent}\label{subsec:phi_delta_choice2}

    The choice of coarse-grid operator for MP2, \cref{eq:mp2}, is analogous to that for MP1 in that we choose $\Phi$ to be the proximal gradient descent operator \cref{eq:pg} and $\Phi_\Delta$ to be its implicit analogue, the alternating proximals operator \cref{eq:apm}, with coarsened step size. Namely, we choose
    \begin{equation}\label{eq:mp2_phi_delta}
        \Phi:=P_{sg}(I-sA_d) \quad \text{ and } \quad \Phi_\Delta^\ell=P_{m^{\ell-1}sg}(I+m^{\ell-1}sA_d)^{-1}, \quad 2\leq \ell \leq \widehat{\ell},
    \end{equation}
     where $g=:\lambda\|(-\cdot)_+\|$, as in \eqref{eq:mp2}. The proximal operator $P_g$ in \cref{eq:mp2_phi_delta} is then given entry-wise by
    \begin{equation}
    \left(P_{\lambda\left\|(-\cdot)_{+}\right\|_1}(\mathbf{u})\right)_i= \begin{cases}(\mathbf{u})_i+\lambda, & (\mathbf{u})_i+\lambda<0, \\
    0, & (\mathbf{u})_i \leq 0 \leq (\mathbf{u})_i+\lambda  \\
    (\mathbf{u})_i, & (\mathbf{u})_i>0,\end{cases},
        \end{equation}
    where the $i$ subscript denotes the $i$th component of the corresponding vector.
    The coarse-grid operators in \cref{eq:mp2_phi_delta} are chosen under the assumption that, on the first coarse grid (i.e., $\ell=2$), 
    \begin{equation}\label{eq:mp2_phi_delta_approx}
        \Phi_{\Delta}:=P_{msg}(I+msA_d)^{-1} \approx \left(P_{sg}(I-sA_d)\right)^m = \Phi^m
    \end{equation}
    is a sufficiently accurate approximation, which is consistent with the assumption in \cref{eq:mp1_phi_delta_approx}.
    }
    
    \ptxt{We note that the heuristic presented in \cref{subsec:phi_delta_choice} and \cref{subsec:phi_delta_choice2}, by which we replace gradient descent terms by proximals, can be applied to any gradient-based method. For nonquadratic objective functions, we may need to invert a nonlinear system to compute these proximals by means of \cref{eq:prox_def2}. Alternatively, we can simply compute the proximals numerically by minimizing the function in the definition \cref{eq:prox} using some iterative method; this computational approach is specially useful when the proximal for the relevant function does not have a closed form. We also note that replacing the explicit gradient steps by implicit proximal steps is only required when taking steps larger than the stability limit of the explicit method. In some machine learning contexts, for example, parameter training stages are performed using arbitrarily small steps well within the stability limit, such that (proximal) gradient descent could be used on coarse levels without drawbacks. Finally, we note that our heuristic still applies in principle to more complex methods with varying step sizes and/or additional parameters. However, determining effective $\Phi_\Delta$ can become a quite complicated task depending on the complexity of the $\Phi$ algorithm.
    }
    
    \subsection{Gradient convergence}\label{subsec:mgrit_grad_conv}

    Given the exact fine-grid trajectory $\{\mathbf{u}_i\}$ and the MGRIT approximation $\{\mathbf{v}_i\}$ of $\{\mathbf{u}_i\}$, MGRIT convergence criteria in our numerical experiments are set in terms of the norm of the MGRIT coarse-grid residual $\mathbf{r}(\mathbf{v})=\mathbf{w}-\mathcal{A}(\mathbf{v})$ over the entire space-iteration domain, where $\mathbf{v} \in \mathbb{R}^{N\cdot N_t}$ (see \cref{eq:Au=g_he}), while the sequential convergence is given in terms of the gradient of the objective function $F(\mathbf{u}_i)$ at iteration $i$, where $\mathbf{u}_i \in \mathbb{R}^{N}$ . 
    The following lemma relates these two measures of convergence by establishing a bound for the (generalized) gradient of the MGRIT approximation $\mathbf{v}_i$ at iteration $i$ in terms of associated error $\mathbf{e}_i=\mathbf{u}_i-\mathbf{v}_i$ and the (generalized) gradient of the exact $\mathbf{u}_i$. 
    The setting for the lemma and subsequent remarks is the general problem of miniziming $F=f+g$ as in \cref{eq:min_f+g}, where  $f$ is a convex $L$-smooth function, and $g$ is a continuous, convex and possibly nondifferentiable function.
    More specifically, we consider the proximal-gradient iteration for minimizing $F=f+g$, with its convergence measured in terms of the  generalized gradient of $F$:
	\begin{equation} \label{eq:G_{sF}}
		G_{sF}\left(\mathbf{u}_i\right):=\frac{1}{s}\left[\mathbf{u}_i-P_{sg}\left(\mathbf{u}_i-s \nabla f\left(\mathbf{u}_i\right) \right)\right]=\frac{1}{s}(I-P_{sg} \circ G_{sf}(\mathbf{u}_i)). 
	\end{equation}
	Note that $P_{sg} \circ G_{sf}$ is the same as the proximal gradient operator \cref{eq:pg}; in particular, $P_{sg} \circ G_{sf}=I-sG_{sF}$. Furthermore, when $g=0$, $F=f$ and $G_{sF}=G_{sf}$. Therefore, the problem of minimizing $f$ \cref{eq:min_f} using the gradient descent method \cref{eq:gd} can be seen as a special case of the problem of minimizing $F=f+g$ with the proximal gradient method \cref{eq:pg}.
    We note that for both MP1 and MP2, we take $s=\frac{1}{L}$. For MP1 and MP2-1D, we have that $L=\|A_1\|_2$, where $A_1$ is the 1D Laplacian matrix \cref{eq:1d_laplacian}, while for MP2-2D, we have that $L=\|A_2\|_2$, where $A_2$ is the 2D Laplacian matrix \cref{eq:2d_laplacian}. 
    Moreover, we consider a general parallel-in-iteration method for an all-at-once system as in \cref{eq:A(u)=g}, such that the following results are not particular to the MGRIT method.

	\begin{lemma}\label{lem:mp2_bound_error}
		Let $f:\mathbb
		R^N\rightarrow\mathbb{R}$ be a convex and $L$-smooth function, $g:\mathbb{R}^N \rightarrow \mathbb{R}$ a continuous convex function. Consider a parallel-in-iteration method described by an all-at-once system of the form \eqref{eq:A(u)=g} with $\Phi:=I-sG_{sF}$ the proximal gradient operator of $f + g$ from \eqref{eq:pg}, with step size $s=\frac{1}{L}$. 
        Let $\mathbf{v}_i \approx \mathbf{u}_i$ be an approximation of the $i$th component of the solution of \eqref{eq:A(u)=g}, with associated error $\mathbf{e}_i=\mathbf{u}_i-\mathbf{v}_i$. 
        Then, the generalized gradient \eqref{eq:G_{sF}} of the approximate solution satisfies 
		\begin{equation}\label{eq:mp2_bound_error}
			\left\| G_{sF}\left(\mathbf{v}_i\right) \right\| \leq 
			L\sqrt{2}\left\|\mathbf{e}_i\right\| +\left\|G_{sF}\left(\mathbf{u}_i\right)\right\|.
		\end{equation}
		
		\begin{proof}
		Since $\mathbf{e}_i=\mathbf{u}_i-\mathbf{v}_i$, we can write $\mathbf{v}_i=\mathbf{u}_i-\mathbf{e}_i$. By \cite[Proposition 29.1]{Bauschke2023}, $\Phi:=I-sG_{sF}$ is $2/3$-averaged. Therefore,
			\begin{equation}\label{eq:2/3-avg}
				\begin{aligned}
					\frac{1}{2}\left\|sG_{sF}(\mathbf{v}_{i})-sG_{sF}(\mathbf{u}_{i})\right\|^2  &\leq  \left\|\mathbf{v}_i-\mathbf{u}_i\right\|^2-\left\|(I-sG_{sF})(\mathbf{v}_i)-(I-sG_{sF})(\mathbf{u}_i)\right\|^2 \\&\leq \left\|\mathbf{v}_i-\mathbf{u}_i\right\|^2= \left\|\mathbf{e}_i\right\|^2.    
				\end{aligned}
			\end{equation}
			Taking the square root on both sides, we get
			\begin{equation}
				\left\|sG_{sF}(\mathbf{v}_{i})-sG_{sF}(\mathbf{u}_{i})\right\|  \leq \sqrt{2}\left\|\mathbf{e}_i\right\|.
			\end{equation}
			Using the triangle inequality, we get
			\begin{equation}
				\begin{aligned}
					\left\| sG_{sF}\left(\mathbf{v}_i\right) \right\| 
                    &=
                    \left\|sG_{sF}\left(\mathbf{v}_i\right) - sG_{sF}\left(\mathbf{u}_i\right)+sG_{sF}\left(\mathbf{u}_i\right)\right\|\\
					&\leq \left\|sG_{sF}\left(\mathbf{v}_i\right) - sG_{sF}\left(\mathbf{u}_i\right)\right\|+\left\|sG_{sF}\left(\mathbf{u}_i\right)\right\| \\
					&\leq \sqrt{2}\left\|\mathbf{e}_i\right\| +s\left\|G_{sF}\left(\mathbf{u}_i\right)\right\|.
				\end{aligned}
			\end{equation}
			Therefore,
			\begin{equation}
				\begin{aligned}
					\left\| G_{sF}\left(\mathbf{v}_i\right) \right\| & \leq \frac{\sqrt{2}}{s}\left\|\mathbf{e}_i\right\| +\left\|G_{sF}\left(\mathbf{u}_i\right)\right\| \\
					&=L\sqrt{2}\left\|\mathbf{e}_i\right\| +\left\|G_{sF}\left(\mathbf{u}_i\right)\right\|.
				\end{aligned}
			\end{equation} 
		\end{proof}
	\end{lemma}

    Bound \cref{eq:mp2_bound_error} in lemma \ref{lem:mp2_bound_error} shows that, when $\mathbf{u}_i$ approximates a minimizer of $F(\mathbf{u})$ with accuracy $\|G_{sF}(\mathbf{u}_i)\|$, then our parallel-in-iteration method will approach the same accuracy $\|G_{sF}(\mathbf{v}_i)\|$ for $\mathbf{v}_i$ if we run sufficient iterations of the parallel-in-iteration to make $\|\mathbf{e}_i\|$ small relative to $\frac{\|G_{sF}(\mathbf{u_i})\|}{L \sqrt{2}}$.
    %
    This can also be restated in terms of residuals rather than errors; we present the result for the linear case in the following remark.
	
    \begin{remark}\label{rem:mp2_bound_res}
	Recall the all-at-once system \cref{eq:A(u)=g} and the related residual equation \cref{eq:mgrit_res}. If $\Phi$ is a linear operator such that lemma \ref{lem:mp2_bound_error} holds, then  $\mathcal{A}\mathbf{e}=\mathbf{r}$, such that $\mathbf{e}=\mathcal{A}^{-1}\mathbf{r}$. Therefore,  $\|\mathbf{e}_i\|\leq\|\mathbf{e}\|=\|\mathcal{A}^{-1}\mathbf{r}\|\leq\|\mathcal{A}^{-1}\|\|\mathbf{r}\|$, which together with \cref{eq:mp2_bound_error} implies
       \begin{equation}\label{eq:mp2_bound_res}
       \left\|G_{sF}\left(\mathbf{v}_i\right)\right\| \leq L\sqrt{2}\|\mathcal{A}^{-1}\|\|\mathbf{r}\|+\left\|G_{sF}\left(\mathbf{u}_i\right)\right\|.
       \end{equation}
        As such, considering large $i \approx N_t$, we find that a sufficiently small space-time residual implies accurate approximation of the minimizer.
        \end{remark}
        %





    
   In particular, let $\mathbf{u}^k_i$ be the approximate solution at time $i$ and iteration $k$, and let $\mathbf{e}^k_i=\mathbf{u}_i-\mathbf{u}^k_i$ be the error of the approximation. Similarly, let $\mathbf{u}^k$ be the approximate solution vector at iteration $k$. Next, we consider the assumption of linear convergence of the all-at-once system. For MGRIT, linear convergence has been proven for certain linear problems as in \cite{DeSterck2024, Southworth2019, Dobrev2017}, and widely observed numerically even for nonlinear problems \cite{DeSterck2023,Falgout2017}, including in the results we present in \cref{subsec:mgrit_mp2_num_res}. 

    \begin{remark}\label{rem:mp2_lin_conv}
		Suppose the parallel-in-iteration method described in lemma \ref{lem:mp2_bound_error} converges linearly in error for the problem of minimizing $F=f+g$ with convergence factor $\rho_{\mathbf{e}}$, that is,
		\begin{equation}
			\left\|\mathbf{e}^k\right\| \leq \rho_{\mathbf{e}}\left\|\mathbf{e}^{k-1}\right\| \leq (\rho_{\mathbf{e}})^k\left\|\mathbf{e}^0\right\|.
		\end{equation}
		for $\rho_{\mathbf{e}} \in(0,1)$ some constant independent of $k$. Then, from \cref{eq:mp2_bound_error} we get linear convergence of the gradient at any point $i$ in time:
		\begin{equation}\label{eq:mp1_lin_conv_error}
			\left\|G_{sF}\left(\mathbf{u}_i^{k}\right)\right\| \leq \left(L\sqrt{2}\left\|\mathbf{e}^{0}\right\|\right) (\rho_{\mathbf{e}})^k+\left\|G_{sF}\left(\mathbf{u}_i\right)\right\|.
		\end{equation}		
		In particular, suppose that after $k=k_*$ iterations $\left(L\sqrt{2}\left\|\mathbf{e}^{0}\right\|\right) (\rho_{\mathbf{e}})^{k_*}:=\epsilon_*$ is some sufficiently small number. Then,
		\begin{equation}
					\left\|G_{sF}\left(\mathbf{u}_i^{k}\right)\right\| \leq \epsilon_*+\left\|G_{sF}\left(\mathbf{u}_i\right)\right\|.
		\end{equation}				
		If we have linear convergence in the residual, that is,
		\begin{equation}
			\left\|\mathbf{r}^k\right\| \leq \rho_{\mathbf{r}} \left\|\mathbf{r}^{k-1}\right\| \leq (\rho_{\mathbf{r}})^k\left\|\mathbf{r}^0\right\|,
		\end{equation}
		for $\rho_{\mathbf{r}} \in(0,1)$ some constant independent of $k$, then from \cref{eq:mp2_bound_res} we obtain an analogous result for convergence of the gradient in terms of the residual for a linear $\Phi$:
		\begin{equation}\label{eq:mp2_lin_conv_res}
			\left\|G_{sF}\left(\mathbf{u}_i^{k}\right)\right\| \leq \left(L\sqrt{2}\left\|\mathcal{A}^{-1}\right\|\left\|\mathbf{r}^{0}\right\|\right) (\rho_{\mathbf{r}})^k+ \left\|G_{sF}\left(\mathbf{u}_i\right)\right\|.
		\end{equation}
	\end{remark}
       
        In particular, we conclude that when using MGRIT with the gradient descent fine-grid operator \cref{eq:gd} for MP1 \cref{eq:mp1}, and the proximal gradient descent fine-grid operator \cref{eq:pg} for MP2 \cref{eq:mp2}, minimizing the MGRIT residual implies minimizing the gradient of $\mathbf{u}_i^k$ to a similar accuracy as that of $\mathbf{u}_i$, with \cref{eq:mp2_lin_conv_res} providing a bound that can inform the number of iterations required to make $\|G_{sF}(\mathbf{u}_i^k)\|$ close to $\|G_{sF}(\mathbf{u}_i)\|$.
        For example, assuming a desired accuracy of the minimizer of $\|G_{sF}(\mathbf{u}_{N_t})\|=\tau$, one could guarantee accuracy $2\tau$ of the MGRIT approximation $\mathbf{u}_{N_t}^{k_*}$ with at most $k_*$ iterations, where $k_*$ is determined according to $\tau=L\sqrt{2}\|\mathcal{A}^{-1}\|\|\mathbf{r}^0\|(\rho_r)^{k_*}$.
    
    \section{Numerical results}\label{sec:num_res}
	Here, we present results of our parallel-in-iteration algorithm applied to the optimization problems MP1 \cref{eq:mp1} and MP2 \cref{eq:mp2} discretized with the  optimization methods  discussed in \cref{subsec:mp1_algos} and \cref{subsec:mp2_algos}, respectively. 
    It is important to note that we are not currently implementing MGRIT in parallel, which will be done in future work; the goal of the current paper is to present the method for optimization problems and investigate its convergence speed (the number of MGRIT iterations required to obtain an accurate solution) numerically, which is done using a serial implementation. 
    In a later section, \cref{sec:speedup}, we also provide estimates for speedup in a parallel setting. We note that while MGRIT by construction approximates the entire trajectory $\{\mathbf{u}_i\}$, $0\leq i \leq N_t$, of the sequential solve, we are for optimization problems really interested in calculating the final approximation $\mathbf{u}_{N_t}$, which is the MGRIT approximate solution at the final iteration point $N_t$, such that we give special importance to $\mathbf{u}_{N_t}$.
    
    First, we solve the problem sequentially using gradient descent or proximal gradient descent. 
    Next, we solve the same problem with MGRIT, i.e., we use MGRIT aiming to accelerate the computation of the fine-grid trajectory $\{\mathbf{u}_i\}$ in parallel. On coarse levels in the MGRIT hierarchy we use an implicit variant of the  optimization method used on the finest level.
    %
    We illustrate convergence results for a particular set of choices for the parameters $N$ (dimension of the Laplacian matrix), $m$ (coarsening factor) and $\widehat{\ell}$ (number of levels). \ptxt{Although we only present results for the fine-level step-size choice $s=\frac{1}{L}$, further tests (not shown here for brevity) indicate that MGRIT performs similarly well for larger step sizes up to the stability limit of $\frac{2}{L}$.} For our initial examples, we chose $\widehat{\ell} =2$ and $m=4$. Later, we present results for different values of $m$, $\widehat{\ell}$ and $N$ to analyze how MGRIT's performance is affected by the choices of these parameters. We choose $10^{-8}$ as the tolerance criteria for all iterative methods. For sequential solves, this applies to the reduction in the 2-norm of the gradient relative to its original value. For MGRIT solves, the tolerance criteria refers to the reduction on the 2-norm of the MGRIT residual, \cref{eq:mgrit_res}. For MP1, \cref{eq:mp1}, we measure the gradient of the quadratic function $f$, while for MP2, \cref{eq:mp2}, we measure the generalized gradient, \cref{eq:G_{sF}}, of the EOP function of the form $F:=f+g$. The initial MGRIT guess $\mathbf{u}^0 \in \mathbb{R}^{N}$ is chosen to be random with components in $[0,1]$, after which the condition $\mathbf{u}_0=\mathbf{g}_0$ is enforced. 

    \ptxt{We note that the sequential solve provides us with $N_t$, exactly the number of iterations required for convergence within the given tolerance threshold, which we then use to specify the MGRIT fine-grid time domain to be $\left[0,\frac{N_t-1}{L}\right]$, with $N_t$ fine-grid time points.
    However, solving the problem sequentially before solving it with MGRIT is of course not recommended in practice. Instead, the final time should be chosen adaptively, for example, using the following strategy.
    Given some tolerance threshold to which we want to reduce the gradient to, we take an educated guess at the final time required to achieve this reduction. This guess can possibly be informed by theoretical or practically observed convergence rates of the underlying optimization method.
    Setting the MGRIT time window with this guess for the final time, we run MGRIT until either the gradient at the final point reaches the prescribed tolerance, or the gradient at the final point stalls.
    
    In the case where the gradient reaches the prescribed tolerance, we either have precisely guessed the correct final time, such that the number of time points is exactly what we would get by solving the problem sequentially first, or we have overestimated the final time such that we have too many time points.
    Otherwise, if the gradient at the final time point stalls, then we have underestimated the final time and do not have enough time points to achieve the desired gradient reduction. Iterating MGRIT past this point is redundant, since MGRIT will only reduce the residual of the underlying algebraic system \cref{eq:mgrit_res}, but without reducing the gradient. 
    In that case, we must start a new MGRIT solve across a new time interval, starting from the final time point of the halted solve and setting the solution at the final time point of the halted solve as the initial guess for the new solve's all-at-once system \cref{eq:A(u)=g}; essentially, this is the same as extending the original trajectory. 
    Then, we take a new educated guess for the final time, which could be estimated based on the approximate convergence rate of the optimization method observed from the halted solve. Having set the new time horizon, we run MGRIT again and repeat this process until the desired gradient reduction is obtained. 
    
    Following this strategy, it may be that underestimating or overestimating the final time does not lead to performing any more work than when perfectly guessing the final time. However, we might lose some parallel speed-up efficiency, as in both cases we do more sequential work than when we work with the ideal number of time points, in which case the parallelization is perfectly optimized. For our discussion on parallel speed-up efficiency, see \cref{sec:speedup}.
    }
    
	\subsection{MGRIT for gradient descent applied to MP1}\label{subsec:mgrit_mp1_num_res}
	
\begin{figure}[h!]
	\centering
	\begin{subfigure}{.5\textwidth}
		\centering
		\includegraphics[scale=0.42]{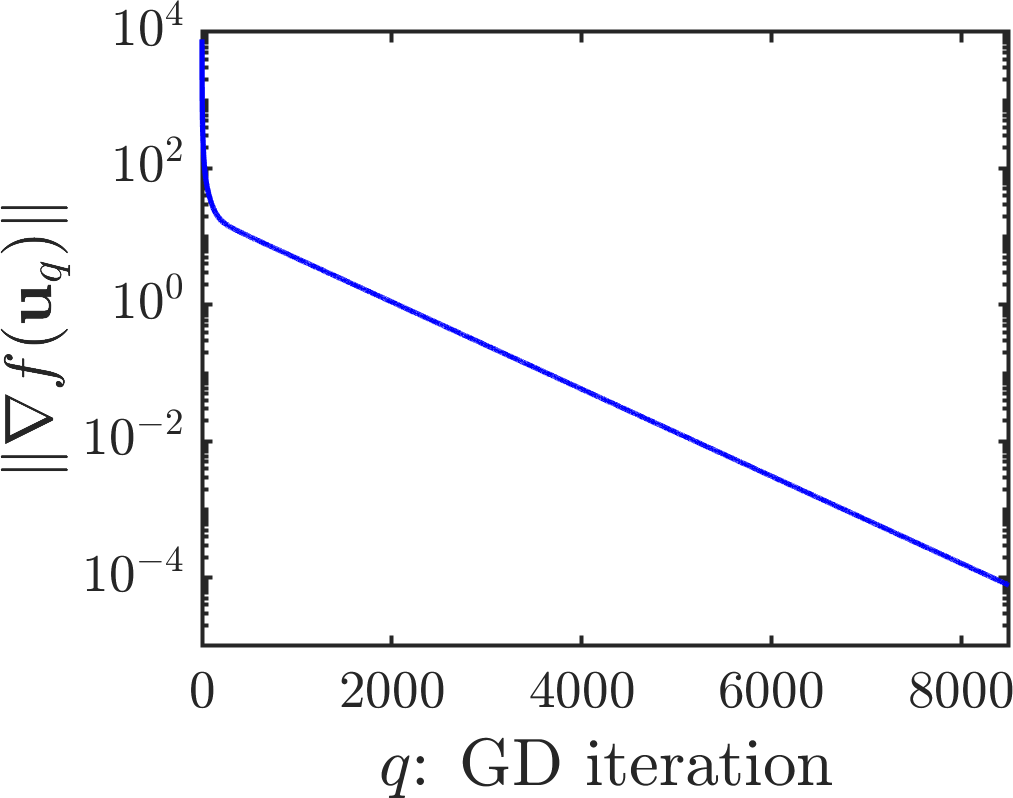}
		\caption{Sequential convergence}
        \label{fig:mp1_conv_seq}
	\end{subfigure}%
	\begin{subfigure}{.5\textwidth}
		\centering
		\includegraphics[scale=0.42]{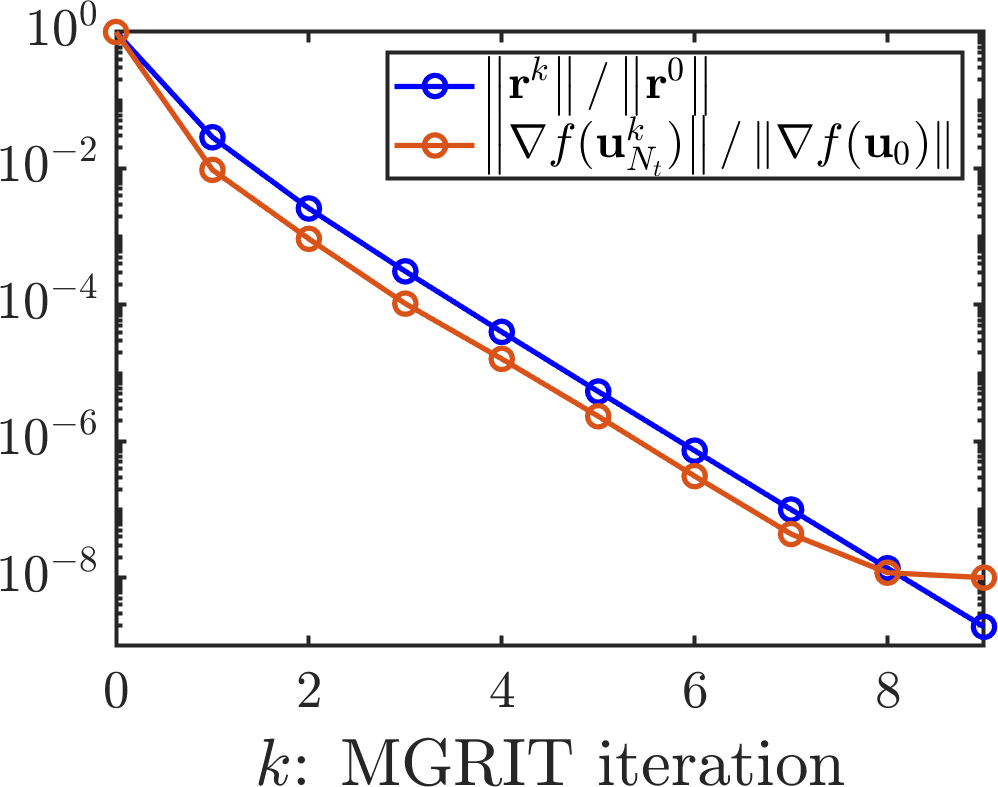}
		\caption{MGRIT convergence}
	\end{subfigure}
	\caption{(a): Convergence of MP1, \cref{eq:mp1}, with a standard sequential solve using gradient descent method \cref{eq:gd}; $\nabla f(\mathbf{u}_q)$ denotes the gradient of the $q$th sequential iterate. (b) Convergence of MP1, \cref{eq:mp1}, with MGRIT using gradient descent operator \cref{eq:gd} on the fine grid and proximal point operator \cref{eq:ppm} on the coarse grid; $\mathbf{r}^{k}$ denotes residual \cref{eq:mgrit_res} at iteration $k$, $\nabla f(\mathbf{u}_{N_t}^k)$ denotes the gradient of the approximate solution at the final time point $N_t$ and MGRIT iteration $k$, $\nabla f(\mathbf{u}_0)$ denotes the gradient of the initial guess $\mathbf{u}_0$, and $\mathbf{r}^0$ denotes the initial MGRIT residual.}
    \label{fig:mp1_conv}
	\end{figure}

	We solve \cref{eq:mp1} using gradient descent, \cref{eq:gd}. The 1D Laplacian matrix $A_1 \in \mathbb{R}^{n \times n}$ in \cref{eq:1d_laplacian} with 2-norm $L$ is taken with $n=40$, the coarsening factor is chosen to be $m=4$ and the number of levels to be $\widehat{\ell}=2$. We also choose $\mathbf{b} \in \mathbb{R}^n$ to be random with components in $[0,1]$. 
    For MGRIT, we choose the proximal-point operator \cref{eq:ppm} as the coarse-grid operator. First, we solve the problem sequentially, with $N_t=8,503$ iterations required for the convergence of the gradient descent algorithm \cref{eq:gd}. We then generate the MGRIT fine grid $t_i=i \Delta t$, $i=0, \ldots, N_t$ with $\Delta t=\frac{1}{L}$. Next, we generate the coarse-grid $T_j=j\Delta t$, $j=0,1,\ldots,N_T$, with $\Delta t=m \Delta t=\frac{4}{L}$, such that the coarse grid has $N_T=2,126$ points. MGRIT then solves \cref{eq:A(u)=g}. Convergence for the sequential and MGRIT solves is illustrated in \cref{fig:mp1_conv}.
    
    In \cref{fig:mp1_conv}, we see the MGRIT convergence both in terms of the MGRIT residual and the gradient at the final time point $N_t$. We recall that the MGRIT convergence criteria is given in terms of space-time residual \cref{eq:mgrit_res}, but using the gradient of solution at the final time point $\mathbf{u}_{N_t}^k$  of each MGRIT iteration as reference, we see that the gradient of the function is also reduced by a similar factor. Recalling remark \ref{rem:mp2_lin_conv}, we can interpret \cref{eq:mp2_lin_conv_res} as follows: if we look at a semi-log plot of $\left\|\nabla f\left(\mathbf{u}_i^{k}\right)\right\|$ as a function of $k$, as in \cref{fig:mp1_conv}(b), then it will be upper bounded by a straight line with slope equal to $\rho_{\mathbf{r}} \in (0,1)$, which is the slope of the residual curve, up until the point where $\left\|\nabla f\left(\mathbf{u}_i\right)\right\|  \approx \left(L\sqrt{2}\left\|\mathcal{A}^{-1}\right\|\left\|\mathbf{r}^{0}\right\|\right) (\rho_{\mathbf{r}})^k$, at which point the curve flattens out to a constant value of $\left\|\nabla f\left(\mathbf{u}_i\right)\right\|$. Therefore, the numerical results presented here are consistent with our theoretical predictions in \cref{subsec:mgrit_grad_conv}.
    
\begin{figure}[h!]		
    \centering
        \begin{subfigure}{.4\textwidth}
		\centering
	\includegraphics[scale=0.42]{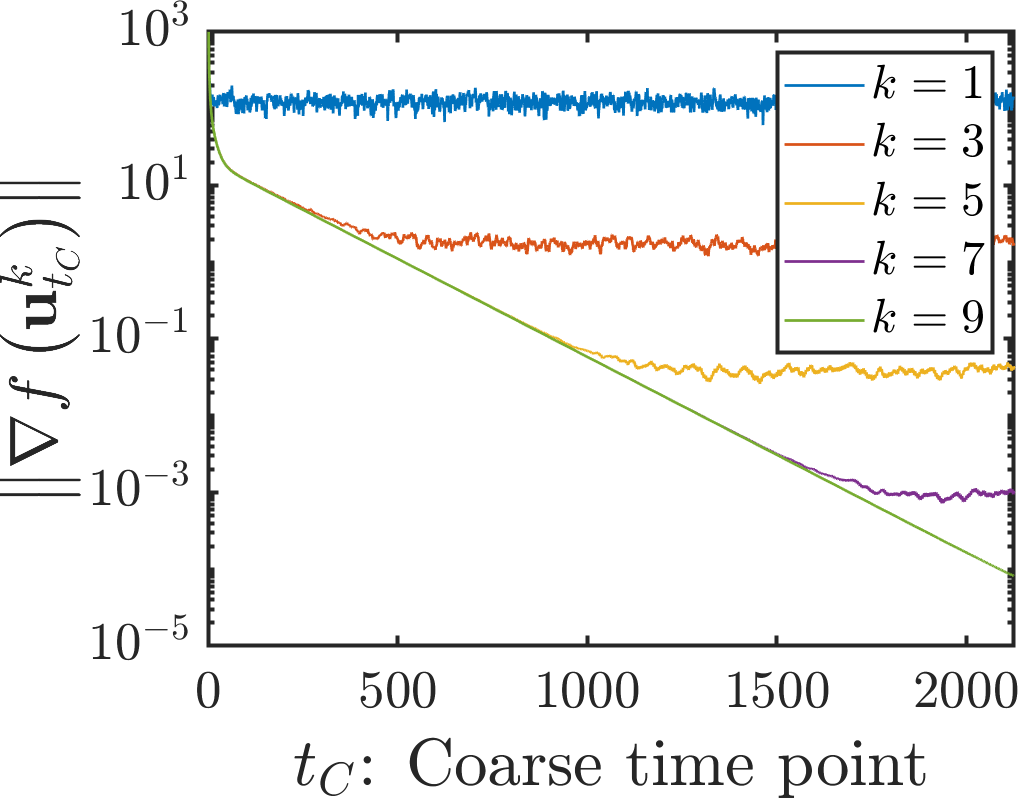}
		\label{fig:mp1_grad_iter}
            \caption{Gradient convergence by iteration}
            \end{subfigure}%
        \begin{subfigure}{.45\textwidth}
		\centering
		\includegraphics[scale=0.42]{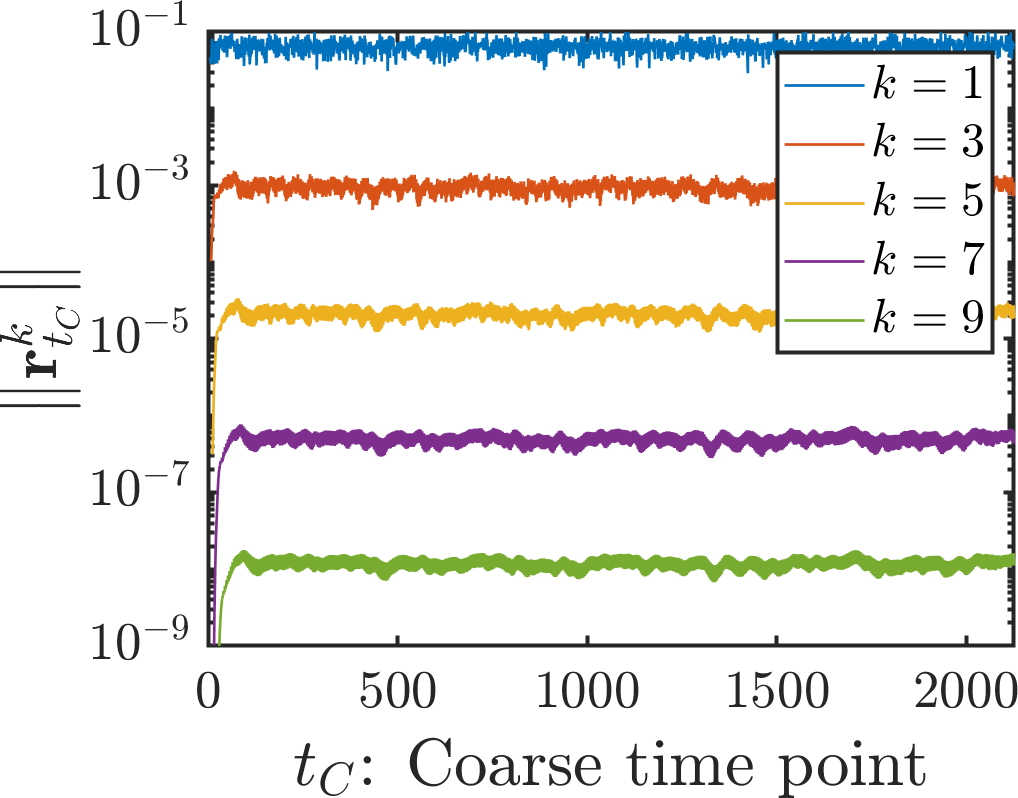}
		\label{fig:mp1_res_iter}
        \caption{Residual convergence by iteration}
         \end{subfigure}%
         \caption{(a): 2-norm of $\nabla f\left(\mathbf{u}^{k}_{t_C}\right)$, where $\nabla f$ is the gradient of $f$, and $\mathbf{u}^{k}_{t_C}$ is the approximate solution at coarse time point $t_C$ and iteration $k$. (b): 2-norm of $\mathbf{r}^{k}_{t_C}$, where $\mathbf{r}^{k}_{t_C}$ is the residual \cref{eq:mgrit_res} at coarse time point $t_C$ and iteration $k$.}
         \label{fig:mp1_conv_iter}
	\end{figure}

    In \cref{fig:mp1_conv_iter}, we further explore the behavior of the gradient and the residual \cref{eq:mgrit_res} as MGRIT converges. As expected, we see that as $k$ increases, the gradient curve for the $k$th iteration in \cref{fig:mp1_conv_iter}(a) approximates the curve in \cref{fig:mp1_conv}(a) increasingly closely. We note that in both \cref{fig:mp1_conv_iter}(a) and (b), the gradient and MGRIT residual are reduced by a consistent factor for each iteration. We also note that the behavior on the initial time points of each iteration in \cref{fig:mp1_conv_iter}(b) is consistent with the exactness property of MGRIT with FCF-relaxations; for $m=4$, at the $k$th MGRIT iteration, the first $2mk=8k$ fine points are solved exactly. \cref{fig:mp1_conv_iter}(b) shows that MGRIT converges in a small number of iterations, which is crucial for obtaining parallel speedup; see \cref{sec:speedup}.
    
	In \cref{table:1d_mp1_conv}, we present MGRIT convergence results for MP1 \cref{eq:mp1} with different coarsening factors $m$ and number of MGRIT levels $\widehat{\ell}$. We note that the larger the coarsening factor, the smaller the maximum number of allowed levels.
	\begin{table}[h!]
    \caption{MGRIT convergence results for different values of $m$ and $\widehat{\ell}$, using gradient descent operator \cref{eq:gd} on the fine grid and  proximal point operator \cref{eq:ppm} on the coarse grid.}\label{table:1d_mp1_conv}	
	\centering
	\begin{tabular}{l|l|l|l|l|l|l}
		 & $\widehat{\ell}=2$ & $\widehat{\ell}=3$ & $\widehat{\ell}=4$ & $\widehat{\ell}=5$ & $\widehat{\ell}=6$ & $\widehat{\ell}=7$ \\ 
		\hline
		$m=4$   & 9        & 9 		  & 9 		 & 9 		& 9		   & 9 \\
		\hline
		$m=16$  & 8 	   & 8        & 8        &          &          &   \\
		\hline
		$m=64$  & 8        & 8        &          &          &          &   \\
		\hline
		$m=256$ & 7        &          &          &          &          &   
	 \end{tabular}\subcaption{$n=40$, $N_t=8,503$}
\end{table} 

    From \cref{subsec:mgrit_ode_link} we recall MP1's \cref{eq:mp1} connection to the 1D heat equation, for which MGRIT convergence behavior is well understood \cite{Falgout2014, Dobrev2017}. For parabolic problems, MGRIT has been observed to have a desirable property of multigrid convergence known as $h$-independent convergence \cite{Trottenberg2000}, which means that the convergence speed does not depend on the size of the fine grid.	

   	\subsection{MGRIT for proximal gradient applied to MP2}\label{subsec:mgrit_mp2_num_res}
	
	We now consider solving the EOP problem, MP2 \cref{eq:mp2}, using proximal gradient method \cref{eq:pg}. For MGRIT, we choose the alternating proximals operator \cref{eq:apm} to be the coarse-grid operator. We noted in \cref{subsec:mp2_algos} how the sequence generated by alternating proximals does not necessarily converge to the same set of minimizers as that of proximal gradient. However, the similarity between these operators justifies the use of the alternating proximal method as the coarse-grid method (recall  discussions in \cref{subsec:opt_disc_ode} and \cref{sec:mgrit_opt} on choosing coarse-grid operators).
			
	\subsubsection{Numerical results for MP2-1D}
	
\begin{figure}[h!]
	\centering
	\begin{subfigure}{.5\textwidth}
		\centering
		\includegraphics[scale=0.42]{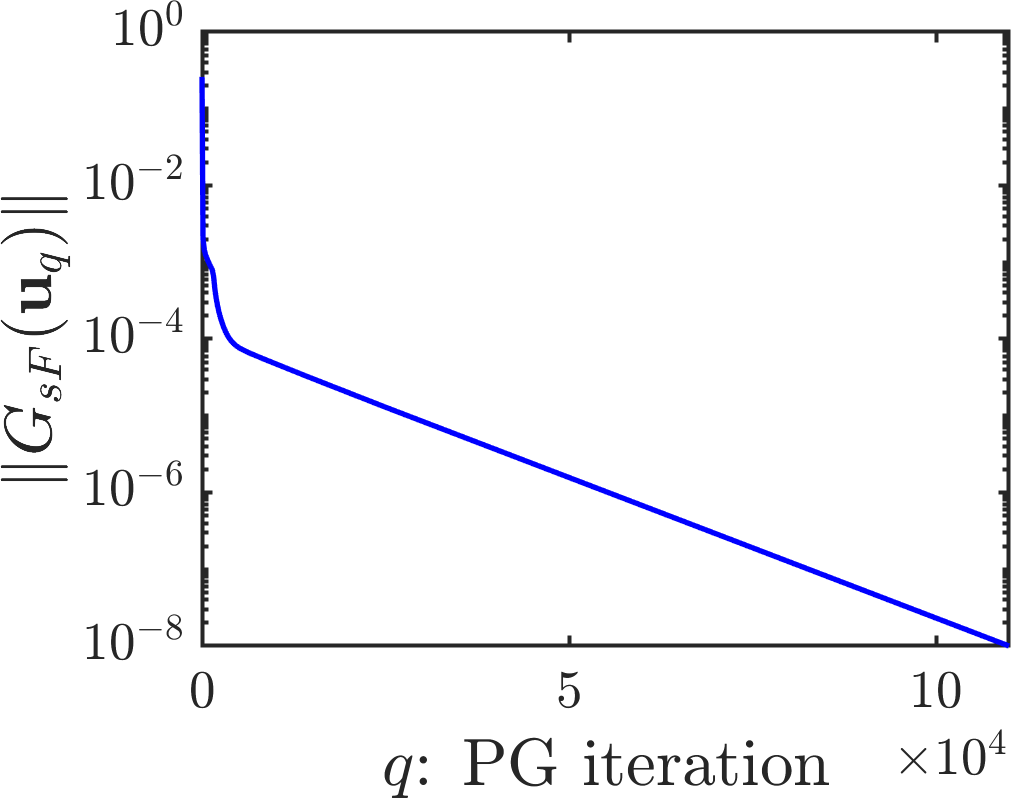}
		\caption{Sequential convergence}
	\end{subfigure}%
	\begin{subfigure}{.5\textwidth}
		\centering
		\includegraphics[scale=0.42]{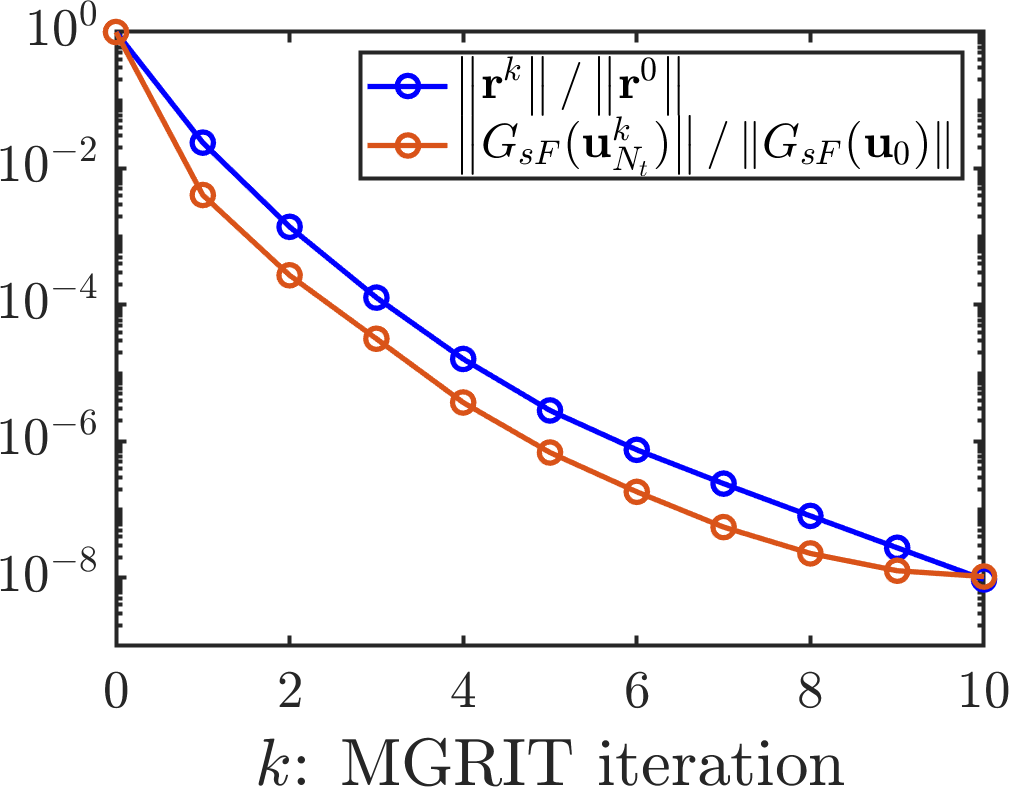}
		\caption{MGRIT convergence}
	\end{subfigure}
	\caption{(a) Convergence of MP2-1D, \cref{eq:mp2}, with a standard sequential solve using proximal gradient descent method \cref{eq:pg}; $G_{sf}(\mathbf{u}_q)$ denotes the generalized gradient \cref{eq:G_{sF}} of the $q$th sequential iterate. (b): Convergence of MP2-1D, \cref{eq:mp2}, with MGRIT using proximal gradient operator \cref{eq:pg} on the fine grid and alternating proximal mappings operator \cref{eq:apm} on the coarse grid; $\mathbf{r}^{k}$ denotes residual \cref{eq:mgrit_res} at iteration $k$, $G_{sF}(\mathbf{u}_{N_t}^k)$ denotes generalized gradient \cref{eq:G_{sF}} of the approximate solution at the final time point $N_t$ and MGRIT iteration $k$, $G_{sF}(\mathbf{u}_0)$ denotes the generalized gradient of the initial guess $\mathbf{u}_0$, and $\mathbf{r}^0$ denotes the initial MGRIT residual.}
	\label{fig:mp2_1d_conv}
	\end{figure}

	For the MP2-1D example, we consider MP2 in \cref{eq:mp2} with $d=1$ and $N=n=256$, with MGRIT parameters $m=4$ and $\widehat{\ell}=2$.  We again start by presenting results for the sequential and MGRIT convergence in \cref{fig:mp2_1d_conv}. 
    In this example, the proximal gradient method requires $N_t=109,888$ iterations to reach the prescribed halting tolerance, i.e., that the generalized gradient of $F$, $G_{sF}$ in \cref{eq:G_{sF}}, be reduced by 8 orders of magnitude. 
    MGRIT reaches the prescribed residual tolerance in 10 iterations.

    \begin{figure}[h!]
        \begin{subfigure}{.45\textwidth}
		\centering
         \end{subfigure}%
	\begin{subfigure}{.5\textwidth}
        \centering
		\includegraphics[scale=0.42]{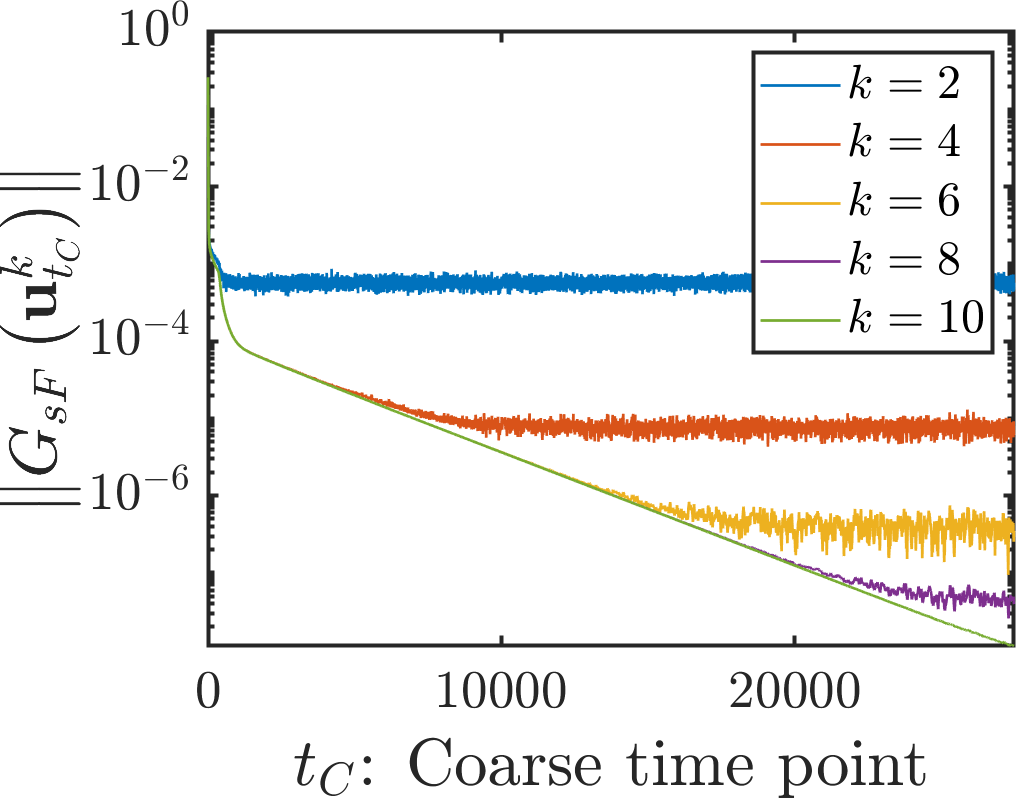}
		\label{fig:mp2_1d_grad}
        \caption{Gradient convergence by iteration}
	\end{subfigure}%
        \begin{subfigure}{.45\textwidth}
            		\centering
		\includegraphics[scale=0.42]{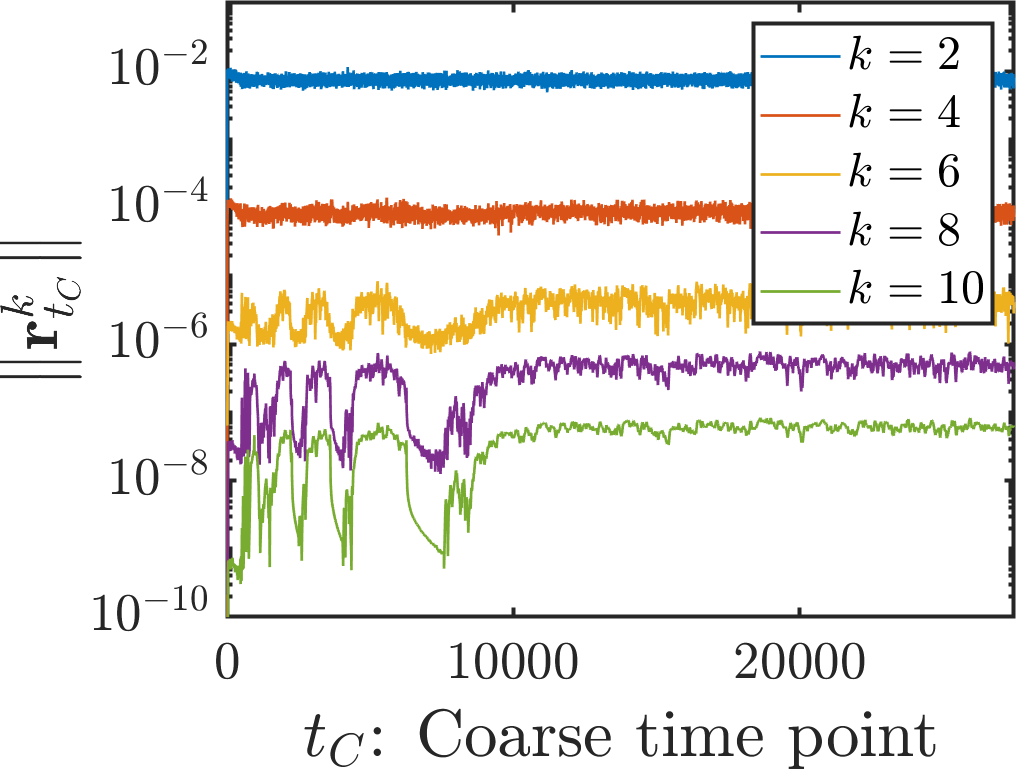}	
		\label{fig:mp2_1d_res}
        \caption{Residual convergence by iteration}
        \end{subfigure}
        	\caption{(a): 2-norm of $G_{sF}\left(\mathbf{u}^{k}_{t_C}\right)$, where $G_{sF}$ is the generalized gradient of $F$ \cref{eq:G_{sF}}, and $\mathbf{u}^{k}_{t_C}$ is the approximate solution at coarse time point $t_C$ and iteration $k$. (b) 2-norm of $\mathbf{r}^{k}_{t_C}$, where $\mathbf{r}^{k}_{t_C}$ is the residual \cref{eq:mgrit_res} at coarse time point $t_C$ and iteration $k$.}
            \label{fig:mp2_1d_conv_iter}
	\end{figure}

\begin{figure}[h!]		
		\centering
		\hspace*{1cm}\includegraphics[scale=0.42]{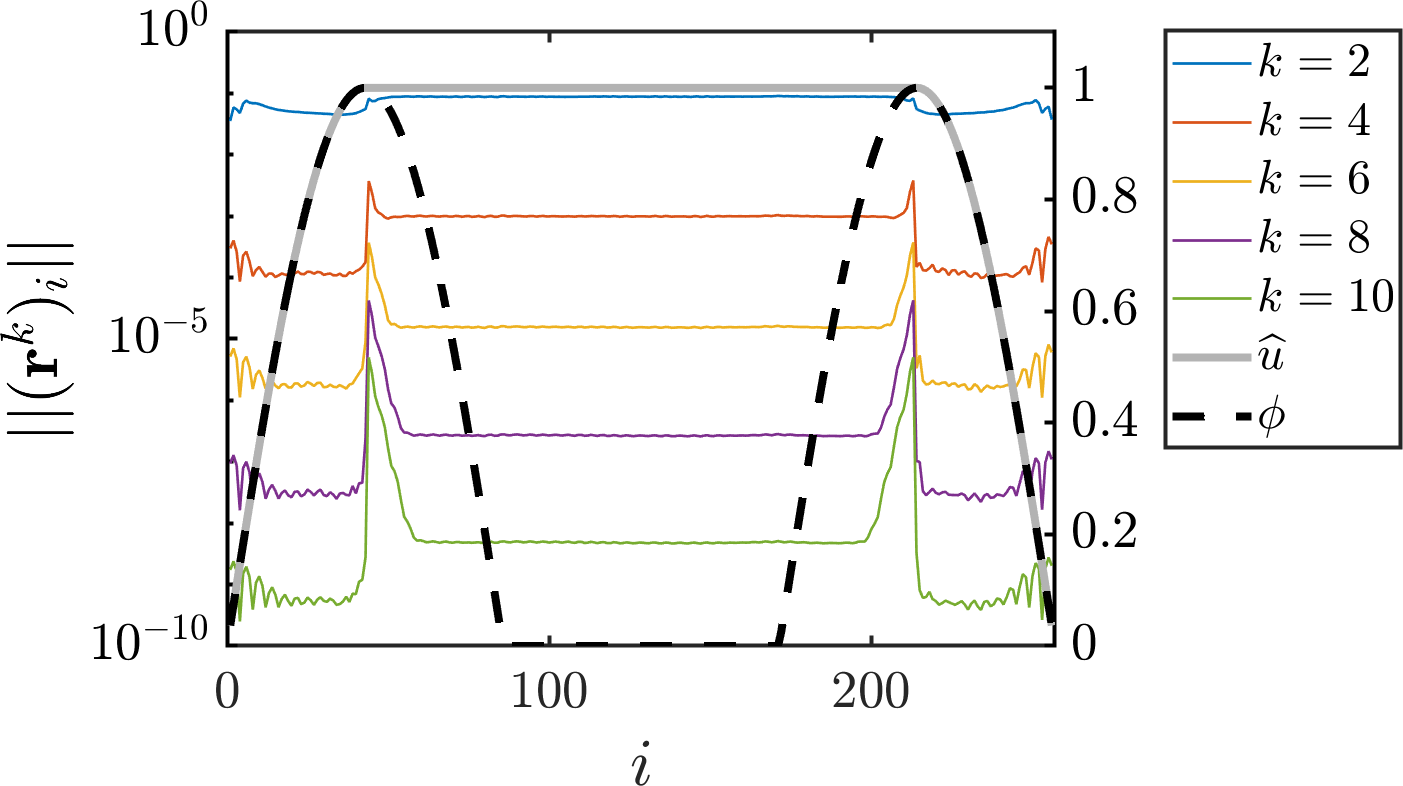}
		\caption{
        Left axis: 2-norm of $(\mathbf{r}^{k})_i$, where $(\mathbf{r}^{k})_{i}$ is the $i$th spatial component of residual \cref{eq:mgrit_res} at iteration $k$. Right axis: $\widehat{u}$ is the exact solution and $\phi$ is the obstacle of the 1D EOP, as in \cref{fig:eop_sol}. The $i$ points form a $256$-point discretization of the interior of the domain $\Omega=[0,3\pi]$ of both $\widehat{u}$ and $\phi$.}
		\label{fig:mp2_1d_obs}
	\end{figure}

	As in \cref{subsec:mgrit_mp1_num_res}, it is interesting to investigate the relation between the MGRIT residual \cref{eq:mgrit_res} and the generalized gradient \cref{eq:G_{sF}}; see \cref{fig:mp2_1d_conv_iter}. Much like in \cref{fig:mp1_conv_iter}(a), the gradient in \cref{fig:mp2_1d_conv_iter}(a) is consistently reduced across iterations. The space-time residual in \cref{fig:mp2_1d_conv_iter}(b), however, has an unusual behavior over the first coarse time points as the number of iterations increases. To further explore this, in \cref{fig:mp2_1d_obs} we look at the spatial components of MGRIT residual vector $\mathbf{r}^k$ \cref{eq:mgrit_res} at iteration $k$ in order to relate the residual to the exact solution $\widehat{u}$ and obstacle $\phi$ of the 1D EOP in \cref{eq:EOP}; see \cref{fig:eop_sol}.
	
	With the solution $\widehat{u}$ and obstacle $\phi$ in \cref{fig:eop_sol} overlaid on the residual \cref{eq:mgrit_res} in \cref{fig:mp2_1d_obs}, we see that the residual decreases more slowly and less smoothly around the spatial coordinates where the membrane attaches to the obstacle and the problem is nonsmooth.
    This nonsmoothness causes the slope of the convergence curves in \cref{fig:mp2_1d_conv}(b) to taper off, compared to \cref{fig:mp1_conv}(b).
    Looking back at \cref{fig:mp2_1d_conv_iter}(b), we observe that for large $k$, when MGRIT approaches the exact trajectory $\{\mathbf{u}_i\}$, the MGRIT residual is initially nonsmooth as a function of time point. This reflects that, at early time points, the set of spatial points where the nonsmooth constraint is active is not fully determined yet. Once the active set is fixed, the residual behaves similarly to how it does for MP1 \cref{eq:mp1}, as seen in \cref{fig:mp1_conv_iter}(b).

	\subsubsection{Numerical results for MP2-2D}
    
\begin{figure}[h!]
		\centering
		\begin{subfigure}{.5\textwidth}
			\centering
			\includegraphics[scale=0.42]{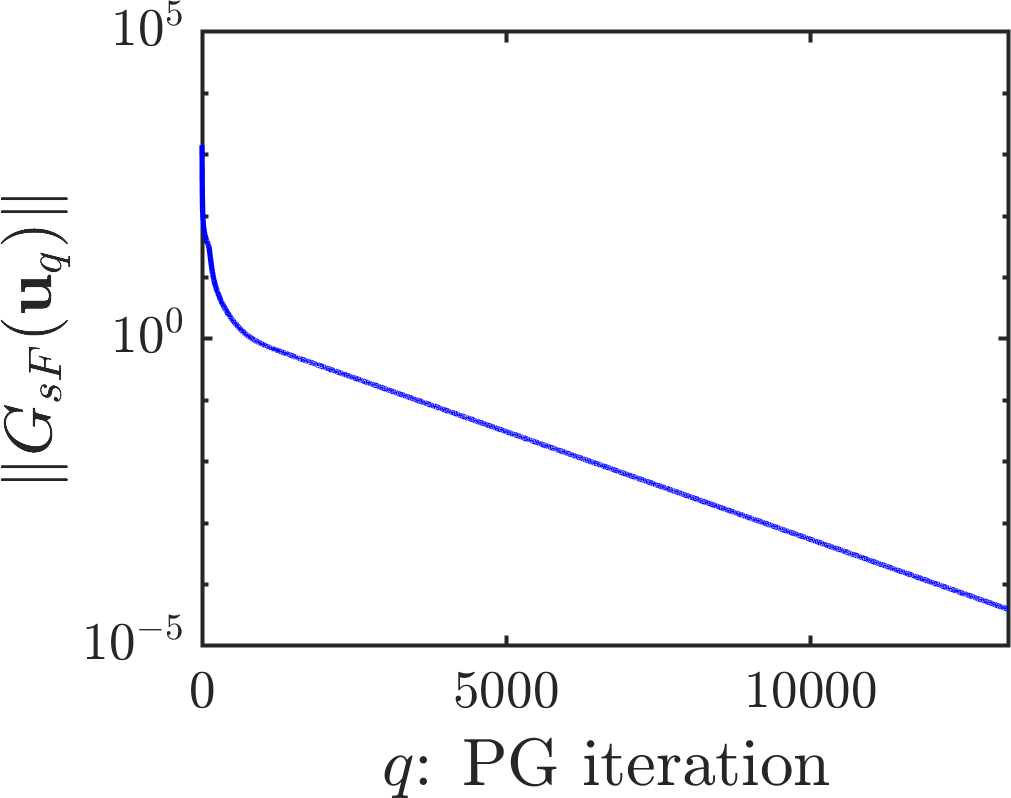}
			\caption{Sequential convergence}
		\end{subfigure}%
		\begin{subfigure}{.5\textwidth}
			\centering
			\includegraphics[scale=0.42]{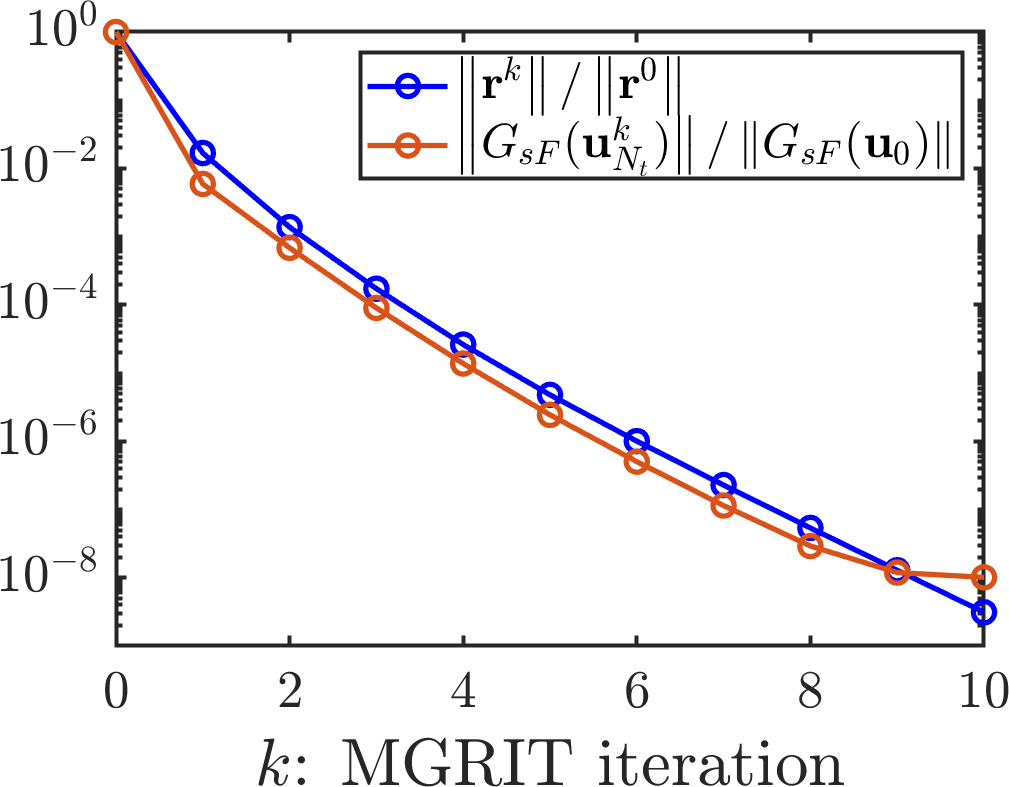}
			\caption{MGRIT convergence}
		\end{subfigure}
		\caption{(a) Convergence of MP2-2D, \cref{eq:mp2}, with a standard sequential solve using  proximal gradient descent method \cref{eq:pg}; $G_{sf}(\mathbf{u}_q)$ denotes  generalized gradient \cref{eq:G_{sF}} of the $q$th sequential iterate. (b): Convergence of MP2-2D, \cref{eq:mp2},  with MGRIT using proximal gradient operator \cref{eq:pg} on the fine grid and alternating proximal mappings operator \cref{eq:apm} on the coarse grid; $\mathbf{r}^{k}$ denotes the residual \cref{eq:mgrit_res} at iteration $k$, $G_{sF}(\mathbf{u}_{N_t}^k)$ denotes the generalized gradient \cref{eq:G_{sF}} of the approximate solution at the final time point $N_t$ and MGRIT iteration $k$, $G_{sF}(\mathbf{u}_0)$ denotes the generalized gradient of the initial guess $\mathbf{u}_0$, and $\mathbf{r}^0$ denotes the initial MGRIT residual.}
		\label{fig:mp2_2d_conv}
	\end{figure}

\begin{figure}[h!]
            \begin{subfigure}{.5\textwidth}
                \centering		
	        \includegraphics[scale=0.42]{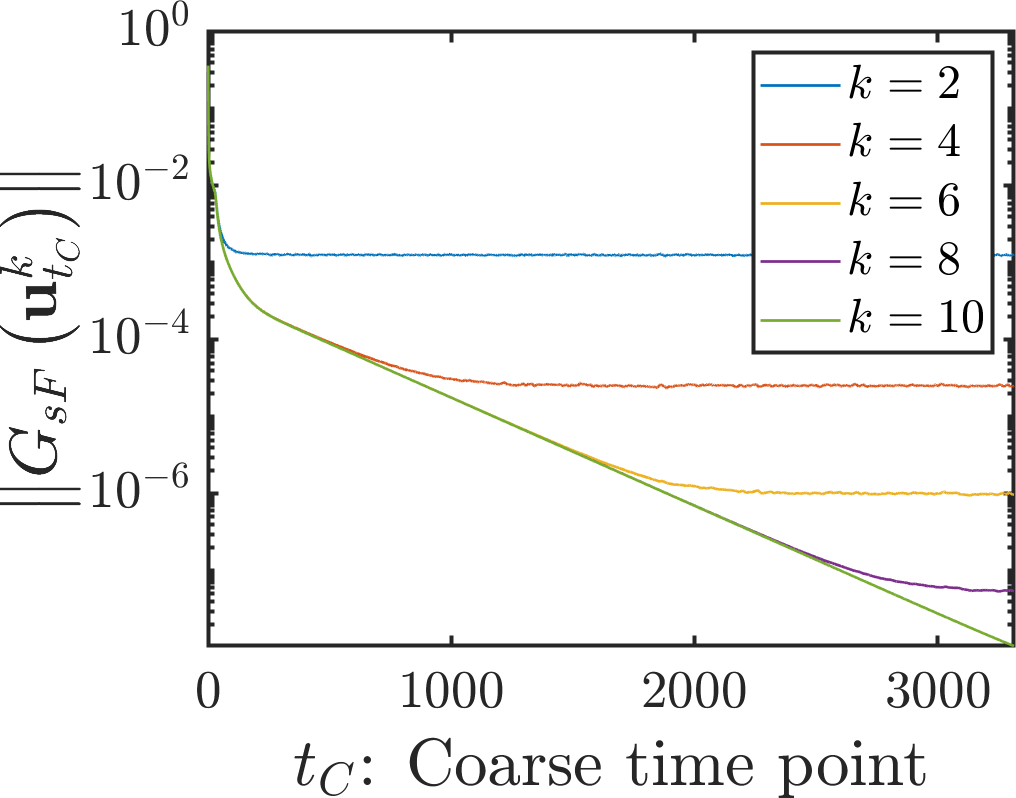}
	        \label{fig:mp2_2d_grad}
            \end{subfigure}
            \begin{subfigure}{.5\textwidth}
 	          \centering
	        \includegraphics[scale=0.42]{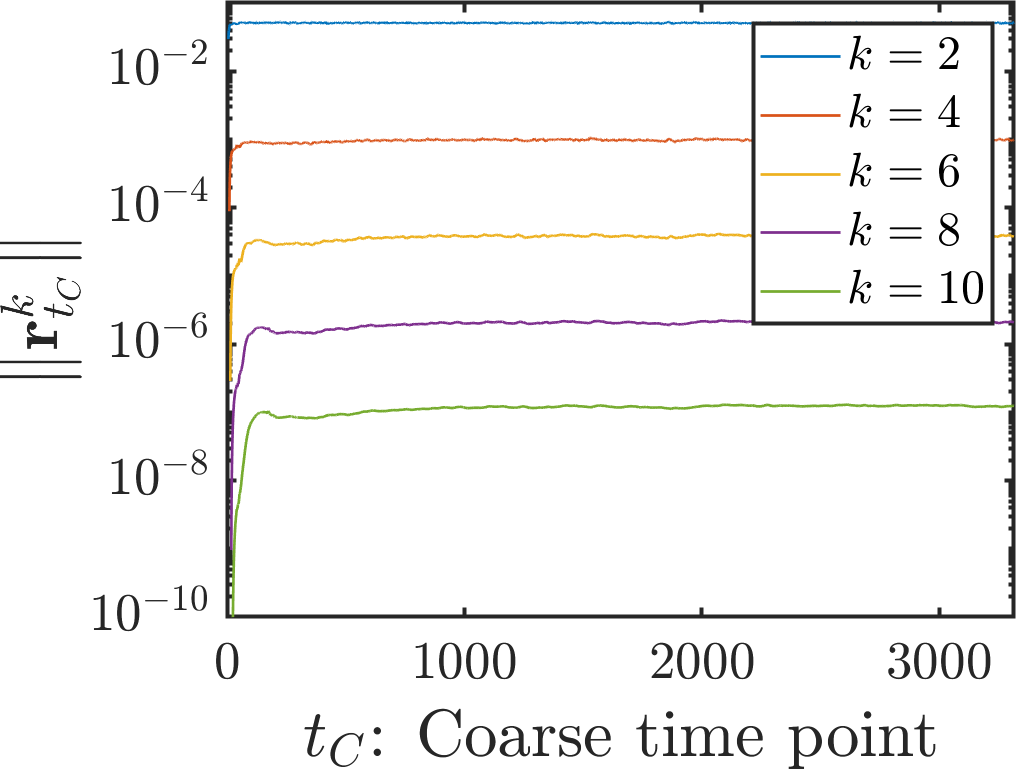}
            \end{subfigure}
            \caption{(a) 2-norm of $G_{sF}\left(\mathbf{u}^{k}_{t_C}\right)$, where $G_{sF}$ is the generalized gradient of $F$, \cref{eq:G_{sF}}, and $\mathbf{u}^{k}_{t_C}$ is the approximate solution vector at coarse time point $t_C$ and iteration $k$. (b) 2-norm of $\mathbf{r}^{k}_{t_C}$, where $\mathbf{r}^{k}_{t_C}$ is residual \cref{eq:mgrit_res} at coarse time point $t_C$ and iteration $k$.}
	    \label{fig:mp2_2d_conv_iter}
        \end{figure}

\begin{figure}[h!]		
		\includegraphics[scale=0.5]{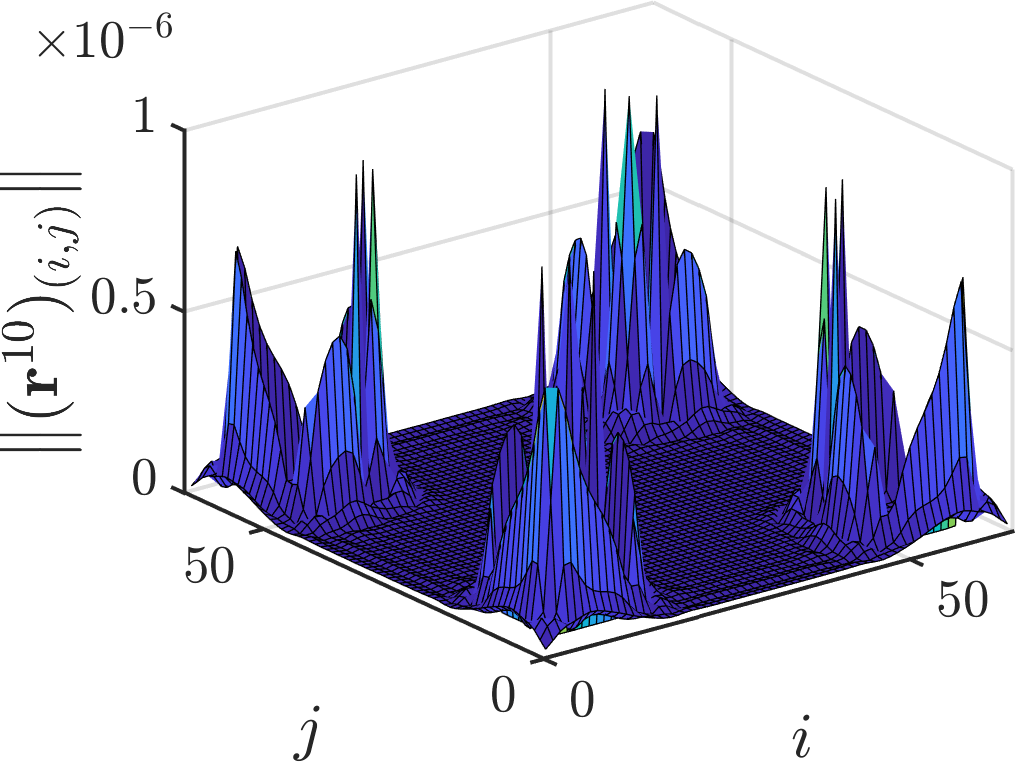}
		\centering
		\caption{2-norm of $(\mathbf{r}^{10})_{(i,j)}$, where $(\mathbf{r}^{10})_{(i,j)}$ is the $(i,j)$th spatial component of residual \cref{eq:mgrit_res} for the $10$th and final MGRIT iteration. 
        }
		\label{fig:mp2_2d_obs}
	\end{figure}

	For the MP2-2D example, we consider MP2 in \cref{eq:mp2} with $d=2$ and $N=n^2=64^2$, with MGRIT parameters $m=4$ and $\widehat{\ell}=2$. As in \cref{subsec:mgrit_mp1_num_res}, we start by presenting the sequential and MGRIT convergence results in \cref{fig:mp2_2d_conv}, noting that the sequential solver required $N_t=13,256$ iterations to converge. We then compare the generalized gradient \cref{eq:G_{sF}} and the MGRIT residual \cref{eq:mgrit_res} at the coarse time points in \cref{fig:mp2_2d_conv_iter}. As with \cref{fig:mp2_1d_obs} for the MP2-1D, we also look at the 2-norm of the spatial components of the residual \cref{eq:mgrit_res} in \cref{fig:mp2_2d_obs}. This time, in order to present a 3D plot, we only present the residual at the final iteration $k=10$.
    
    \begin{figure}[h!]		
		\centering
		\includegraphics[scale=0.42]{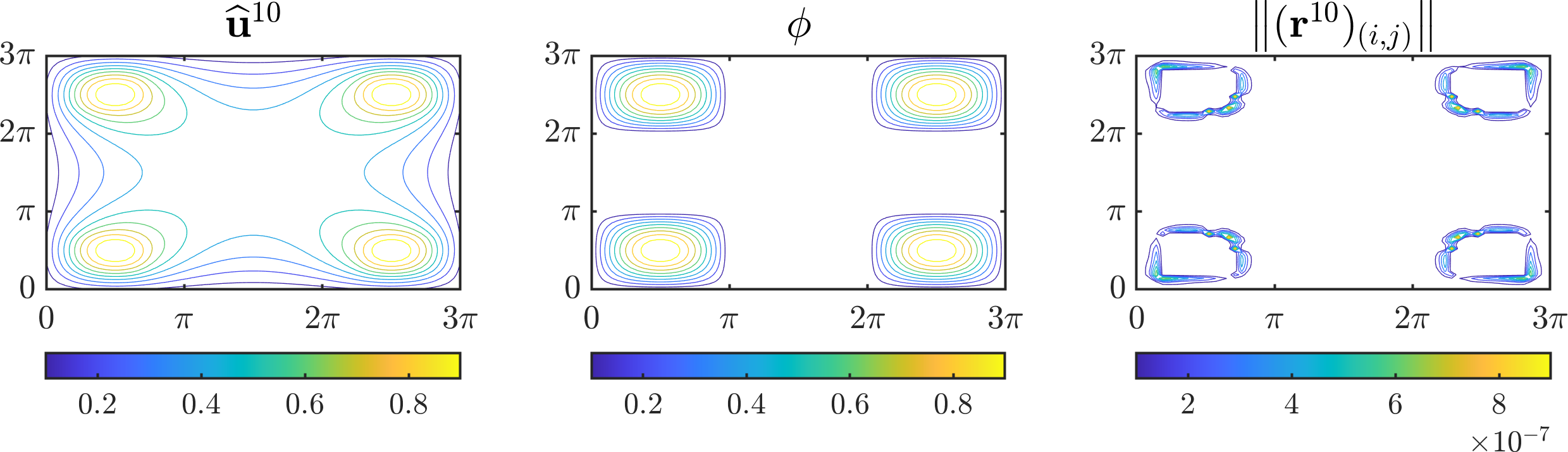}
		\caption{Contour plots of the approximate solution $\mathbf{\widehat{u}}$ and obstacle $\phi$ for the 2D EOP, as in \cref{fig:eop_sol}, and the residual of the 10th and final MGRIT iteration in \cref{fig:mp2_2d_obs}.}
		\label{fig:mp2_2d_cont}
	\end{figure}

    Comparing \cref{fig:mp2_2d_obs} to \cref{fig:eop_sol}, we see a similar relationship between the residual, the solution and the obstacle as we did with \cref{fig:mp2_1d_obs} and \cref{fig:eop_sol} for MP2-1D. To better illustrate this relationship, we present countours plots in \cref{fig:mp2_2d_cont}, from which we observe that the residual of the final MGRIT iteration is larger in the regions where the membrane attaches to the obstacle, as with the 1D case. Overall, we observe that MGRIT for MP2-2D behaves similar to that described for MP2-1D in \cref{subsec:mgrit_mp1_num_res}, and MGRIT convergence for both MP2-1D and MP2-2D behave according to the theory in \cref{subsec:mgrit_grad_conv}. 
    
	As with MP1 in \cref{table:1d_mp1_conv}, we investigate the convergence results for MP2, \cref{eq:mp2}, with different coarsening factors $m$ and number of MGRIT levels $\widehat{\ell}$. This time, we also consider different sizes for the matrices $A_1$ in \cref{eq:1d_laplacian} for MP2-1D and $A_2$ in \cref{eq:2d_laplacian} for MP2-2D. Overall, MGRIT consistently converges in a small number of iterations, as for MP1 in \cref{table:1d_mp1_conv}, and we obtain iteration numbers that are almost independent of problem size.
	
	\begin{table}[h!]
	\caption{MP2-1D convergence results for different values of $n$, $m$ and $\widehat{\ell}$, using  proximal gradient operator \cref{eq:pg} on the fine grid and  alternating proximals operator \cref{eq:apm} on the coarse grid.}\label{table:1d_mp2_conv}
		\centering
		\begin{tabular}{l|l|l|l|l|l|l|l}
			& $\widehat{\ell}=2$ & $\widehat{\ell}=3$ & $\widehat{\ell}=4$ & $\widehat{\ell}=5$ & $\widehat{\ell}=6$ & $\widehat{\ell}=7$ & $\widehat{\ell}=8$ \\ 
			\hline
			$m=4$   & 10        & 10 		 & 11 		& 12 	   & 12	   & 12       & 13 \\
			\hline
			$m=16$  & 8 	   & 9        & 10        &  10        &        &         & \\
			\hline
			$m=64$  & 10        & 10        &  10        &          &          &         & \\
			\hline
			$m=256$ & 10        & 10         &          &          &          &         & \\
			\hline
			$m=1024$ & 11        &          &          &          &          &         & \\
		\end{tabular}\subcaption{$n=256$, $N_t=109,888$.}
	\centering	
	\begin{tabular}{l|l|l|l|l|l|l|l|l|l}
		& $\widehat{\ell}=2$ & $\widehat{\ell}=3$ & $\widehat{\ell}=4$ & $\widehat{\ell}=5$ & $\widehat{\ell}=6$ & $\widehat{\ell}=7$ & $\widehat{\ell}=8$ & $\widehat{\ell}=9$ & $\widehat{\ell}=10$ \\ 
        \hline
		$m=4$   & 9   & 12 & 11 & 11  & 12	& 12 & 12 & 12 & 12 \\
		\hline
		$m=16$  & 10 & 10        & 11        &  11  & &      &        &         & \\
		\hline
		$m=64$  & 8        & 10        &  10        &     & &     &          &         & \\
		\hline
		$m=256$ & 9       & 9    & &     &          &          &          &         & \\
		\hline
		$m=1024$ & 10 & &        &          &          &          &          &         & \\
	\end{tabular}\subcaption{$n=512$, $N_t=369,612$.}
	\end{table}

	\begin{table}[h!]
		\caption{MP2-2D convergence results for different values of $n$, $m$ and $\widehat{\ell}$, using proximal gradient operator \cref{eq:pg} on the fine grid and  alternating proximals operator \cref{eq:apm} on the coarse grid.}\label{table:2d_mp2_conv}	
			\centering
		\begin{tabular}{l|l|l|l|l|l|l}
				& $\widehat{\ell}=2$ & $\widehat{\ell}=3$ & $\widehat{\ell}=4$ & $\widehat{\ell}=5$ & $\widehat{\ell}=6$  & $\widehat{\ell}=7$ \\ 
				\hline
				$m=4$   & 10        & 13 		  & 13 		 & 13 		& 13 &		13    \\
				\hline
				$m=16$  & 11 	   & 12        &  12       &          &       &   \\
				\hline
				$m=64$  & 10        &  10       &          &          &        &  \\
				\hline
                $m=256$  &   9     &         &          &          &         & 
		\end{tabular}\subcaption{$n^2=64^2$, $N_t=13,256$.}
        		\centering
			\begin{tabular}{l|l|l|l|l|l|l|l}
			& $\widehat{\ell}=2$ & $\widehat{\ell}=3$ & $\widehat{\ell}=4$ & $\widehat{\ell}=5$ & $\widehat{\ell}=6$ & $\widehat{\ell}=7$  & $\widehat{\ell}=8$\\ 
			\hline
			$m=4$   & 10         &  	10	  &  10		 &  10 		& 	11	 & 11 & 11  \\
			\hline
			$m=16$  &  9	   &       9  &       9  &          &    & &      \\
			\hline
			$m=64$  &   9      &      9   &          &          &     & &     \\
			\hline
                $m=256$  & 9        &         &          &          &      & &    
		\end{tabular}\subcaption{$n^2=128^2$, $N_t=46,431$.}
	\end{table}	

	\section{MGRIT parallel speedup estimates}\label{sec:speedup}
    

	In our numerical results, we verified that MGRIT was successful in solving both MP1, \cref{eq:mp1}, and MP2, \cref{eq:mp2}, in a small number of iterations. However, each MGRIT iteration does considerably more work than a full sequential solve. Since we are interested in the wall-clock time reduction in solving in parallel, we need to compute the relative cost of the MGRIT solve compared to the sequential solve, and estimate how much time would be saved by solving the problem in parallel. For that, we must first derive a formula for estimating the wall-clock time speedup as a function of the relevant parameters. Based on \cite{Ruprecht2012} and \cite{Minion2010}, which provide speedup estimates for the parareal method, we derive a formula for the MGRIT speedup estimate. 
    In the following estimates, we assume that parallel communication costs are negligible, as assumed for the parareal speedup estimates in \cite{Minion2010, Ruprecht2012}.
    %
    %
	\subsection{Estimating 2-level MGRIT wall-clock time speedup}\label{subsec:speedup_2}

    First, we consider a 2-level MGRIT scheme as in \cref{subsec:mgrit}. Let $N_f$ denote the number of time intervals on the fine-grid level $\ell=1$, and let $t_f$  denote the wall-clock time required to apply the fine-grid time-stepping operator $\Phi$ once. Similarly, let $N_c$ denote the total number of time intervals on the coarse-grid level $\ell=2$, and let $t_c$ denote the wall-clock time required to apply the coarse-grid time-stepping operator $\Phi_\Delta$ once. 
    The coarsening factor, $m$, is given by $m=\frac{N_f}{N_c}$. Let $N_p$ denote the number of processors; full parallelization is obtained when we take $N_p=N_c=\frac{N_f}{m}$, such that we allot each coarse-grid step to its own processor, along with $m$ fine-grid steps. Finally, let $N_{it}$ be the number of MGRIT iterations required for convergence to a set tolerance.
    
    The speedup, given by the ratio of the computation time for the sequential solve and the MGRIT solve, is given by
	\begin{equation}\label{eq:s2}
		S_2\left(N_p\right)=\frac{N_f t_f}{N_{i t}\left(N_c t_c+\frac{2 N_f t_f}{N_p}\right)}.
	\end{equation}
	The numerator $N_f t_f$ represents the total cost of the sequential solve. In the denominator, we have the cost of the MGRIT solve: for each of the $N_{it}$ MGRIT iterations, the $\frac{2 N_f t_f}{N_p}$ term accounts for the fine-grid cost of the parallel FCF-relaxation and correcting the coarse-grid solution, which in total amount to the cost of two full sequential solves that are solved in parallel by $N_p$ processors, while the $N_c t_c$ term accounts for the cost of computing the C-point residual. Recalling that $\frac{N_f}{N_c}=m$ and assuming full parallelization (i.e., $N_p=\frac{N_f}{m}$), we rearrange \cref{eq:s2} to find
	\begin{equation}\label{eq:s2_fp}
		S_2\left(m\right)=\frac{N_f t_f}{N_{i t}\left(N_c t_c+\frac{2 N_f t_f}{N_p}\right)}=\frac{1}{N_{i t}\left(\frac{1}{m} \frac{t_c}{t_f}+\frac{2}{N_p}\right)}=\frac{1}{N_{i t}\left(\frac{1}{m} \frac{t_c}{t_f}+\frac{2 m}{N_f}\right)}.
	\end{equation}
	In \cref{sec:num_res}, we observed that $N_{it}$ does not vary significantly with $m$, so in what follows we assume it is a constant independent of $m$. We further assume $t_c$ is independent of $m$, such that setting the derivative of $S_2$ with respect to $m$ equal to zero then yields
	\begin{equation}\label{eq:s2_m}
		S_2^{\prime}\left(m\right)=0 \iff -\frac{1}{m^2} \frac{t_c}{t_f}+\frac{2}{N_f}=0 \iff m=\sqrt{\frac{N_f}{2} \frac{t_c}{t_f}},
	\end{equation}    
    where for simplicity we assume the minimizer $\sqrt{\frac{N_f}{2} \frac{t_c}{t_f}}$ is an integer, such that we can choose $m=m^*:=\sqrt{\frac{N_f}{2} \frac{t_c}{t_f}}$ as the optimal coarsening factor; in our results in the tables below, we round $\sqrt{\frac{N_f}{2} \frac{t_c}{t_f}}$ to the nearest integer. Denote $\alpha:=\frac{t_c}{t_f}$, such that $m^*=\sqrt{\frac{\alpha N_f}{2}}$. Taking $m=m^*$ in \cref{eq:s2_fp}, and simplifying gives
	\begin{equation}\label{eq:s2_opt}
		S_2\left(m^*\right)=\frac{1}{N_{i t}\left(\frac{1}{\sqrt\frac{\alpha N_f}{2}}\alpha+2\frac{\sqrt{\frac{\alpha N_f}{2}}}{N_f}\right)}=\frac{1}{N_{it}\left(\sqrt{\frac{2\alpha}{N_f}}+\sqrt{\frac{2\alpha}{N_f}}\right)}=\frac{\sqrt{N_f}}{N_{i t} 2 \sqrt{2\alpha}}.
	\end{equation}
    We observe that the optimal speedup is a function of $N_f$, $N_{it}$ and $\alpha:=\frac{t_c}{t_f}$. 
    Since the optimal speedup is proportional to $\sqrt{N_f}$, it is clear that more speedup can be obtained and MGRIT will be more effective when the fine-level optimization method converges slowly and requires a large $N_f$. Similarly, to keep $N_{it}$ small it is important to choose a coarse-grid operator that closely matches the fine operator, see \cite{Dobrev2017,DeSterck2021}.
    Finally, $\alpha$ can be brought close to 1 by making $t_c$ as close as possible to $t_f$, for example, by optimizing the implicit computations done with the coarse-grid operator or by substituting it by a different one with smaller $t_c$, without increasing $N_t$. 
    
	
	\subsection{Estimating 3-level MGRIT wall-clock time speedup}\label{subsec:speedup_3}
    
	Next, we introduce a third level, $\ell=3$, which is obtained from coarsening the $\ell=2$ coarse grid level by the same coarsening factor $m$. Associated with $\ell=3$ is the time-stepping operator $\Phi_\Delta^3$, defined according to \cref{eq:mp1_phi_delta_l} (i.e., an operator similar to $\Phi_\Delta$ with an $m$ times larger step size), which takes $m^2 \Delta t$-sized steps.  Let $N_{cc}$ be the number of time intervals on the $\ell=3$ level, and let $t_{cc}$ denote the wall-clock time required to compute one $\Phi_\Delta^3$ step. We assume that $t_{cc}=t_c$, since we use the same implicit scheme on both coarse levels $\ell=2$ and $\ell=3$, and $N_{cc}=\frac{N_{c}}{m}=\frac{N_f}{m^2}$, assuming $N_{cc}$ is an integer.
    
    For $\widehat{\ell}=3$, we get a speedup estimate formula of
	\begin{equation}\label{eq:s3}
		\begin{aligned}
			S_3\left(N_p\right)&=\frac{N_f t_f}{N_{i t}\left(N_{cc} t_{cc}+\frac{3 N_{c} t_{c}}{\min\{N_p,N_c/m\}}+\frac{2 N_f t_f}{N_p}\right)}.
		\end{aligned}
	\end{equation}
	The reasoning is similar to that for $\widehat{\ell}=2$, with the difference that here we replace the original coarse-grid work in \cref{eq:s2}, with cost $N_c t_c$, with a new 2-grid MGRIT solve. Note that on the $\ell=2$ level, a post F-relaxation step must be done to distribute the correction from the coarse-grid solve in $\ell=3$, such that the total relaxation cost on the $\ell=2$ grid is roughly that of 3 sequential coarse solves. Since there are only $\frac{N_c}{m}$ intervals on the first coarse level, we can use at most $\frac{N_c}{m}$ processors on that level. Assuming full parallelization on the fine level $N_p=N_c$, we have that $\min\{N_p,\frac{N_c}{m}\}=\min\{N_c,\frac{N_c}{m}\}=\frac{N_c}{m}$. In that case,	
	\begin{equation}\label{eq:s3_fp}
		\begin{aligned}
			S_3\left(N_p\right)=\frac{N_f t_f}{N_{i t}\left(N_{cc} t_{cc}+\frac{3 N_c t_c}{N_c/m}+\frac{2 N_f t_f}{N_c}\right)}=\frac{1}{N_{i t}\left(\frac{1}{m^2}\frac{t_c}{t_f}+\frac{2m}{N_f}\left(1+\frac{3}{2}\frac{t_c}{t_f}\right)\right)}.
		\end{aligned}
	\end{equation}
	Following a similar reasoning to that for finding the optimal 2-level coarsening factor \cref{eq:s2_m}, we set the derivative of $S_3$ with respect to $m$ equal to zero and get 
	\begin{equation}\label{eq:s3_m}
		S_3^{\prime}\left(N_p\right)=0 	\iff -\frac{1}{m^3}\frac{t_c}{t_f} + \frac{2}{N_f}+\frac{2}{N_f}\frac{3}{2}\frac{t_c}{t_f}=0  \iff  m=\sqrt[3]{\frac{N_f}{2}\left(\frac{1}{t_f/t_c+3/2}\right)}.
	\end{equation}
        
    Assuming $\sqrt[3]{\frac{N_f}{2}\left(\frac{1}{t_f/t_c+3/2}\right)}$ is a integer, we take $m=m^*:=\sqrt[3]{\frac{N_f}{2}\left(\frac{1}{t_f/t_c+3/2}\right)}$ to be the optimal coarsening factor; we take the nearest integer, otherwise. We note that this analysis can be further extended for more levels, but here we restrict ourselves to the 2- and 3-level cases.


	\subsection{MGRIT parallel speedup estimates for MP1 and MP2}
	Using the formulas derived for 2 and 3 levels, we now present the computed MGRIT parallel speedup estimates for MP1 in \cref{eq:mp1} and MP2 in \cref{eq:mp2} in the following tables. First, we compute $m^*$ for $\ell=2$ and $\ell=3$ using \cref{eq:s2_m} and \cref{eq:s3_m}, respectively; if the resulting value is not an integer, we round it to the nearest integer and call that $m^*$ instead. We choose $m=m^*$ as the coarsening factor for the MGRIT solve, which gives us $N_{it}^*$ iterations in our numerical experiments; note that this is not the optimal number of MGRIT iterations, but simply the value of $N_{it}$ obtained experimentally when using $m=m^*$. We then use $m=m^*$ and $N_{it}=N_{it}^*$ in the speedup formulas \cref{eq:s2_fp} and \cref{eq:s3_fp}, where we also use an empirically estimated $\alpha=\frac{t_c}{t_f}$. 
    %
	
	Next, we consider a hypothetical \emph{ideal} case in which $\frac{t_c}{t_f}=1=\widehat{\alpha}$. Even if $\frac{t_c}{t_f}=1$ is unlikely, this provides a lower bound for the value of the optimal coarsening factor $m^*$, since $m^*$ in \cref{eq:s2_m} and \cref{eq:s3_m} is increasing with respect to $\frac{t_c}{t_f}$. Following the previously outlined process, we compute the optimal $\widehat{m}^*$, and then obtain $\widehat{N_{it}}^*$ numerically, which are used alongside $\widehat{\alpha}=1$ for the 2- and 3-level speedup formulas to provide an estimated upper bound for speedup. We note that since the fine grid, $t_i=i \Delta t$, $i= 0,1, \ldots, N_t$, has $N_t+1$ points, there are $N_f=N_t$ fine-grid intervals. Similarly, $N_c=N_T$.
	
	\begin{table}[h!]
		\caption{Speedup estimates for MP1, \cref{eq:mp1}, different values of $n$ and $\widehat{\ell}$, using gradient descent operator \cref{eq:gd} on the fine grid and proximal point operator \cref{eq:ppm} on the coarse grid. Here, $m^*$ and $S^*$ are computed based on the numerically estimated $\frac{t_c}{t_f}$ values given by $\alpha=2.9$ for $n=40$, $\alpha=3.1$ for $n=60$, and $\alpha=3.5$ for $n=80$. $N_p$ denotes the number of processors required to achieve optimal speedup.}\label{table:mp1_speedup}
       \begin{minipage}[c]{0.33\linewidth}
    \centering
		\begin{tabular}{l|l|l|l|l}
			$\widehat{\ell}$ & $m^*$ & $N_{it}^*$ & $S^*_\ell$ & $N_p$\\
			\hline
			$2$ & $112$ & $7$ & $2.73$ & $76$ \\
			\hline
			$3$ & $27$ & $8$  & $3.29$ & $315$
		\end{tabular}
		\subcaption{$n=40$, $N_t=8,503$.}
    \end{minipage}\begin{minipage}[c]{0.33\linewidth}
        \centering
        \begin{tabular}{l|l|l|l|l}
			$\widehat{\ell}$ & $m^*$ & $N_{it}^*$ & $S^*_\ell$ & $N_p$\\
			\hline
			$2$ & $167$ & $8$  & $3.37$ & $109$\\
			\hline
			$3$ & $35$ & $8$ & $5.11$ & $516$
		\end{tabular}\subcaption{$n=60$, $N_t=18,050$.}
        \end{minipage}
        \begin{minipage}[c]{0.33\linewidth}
	\centering	
        \begin{tabular}{l|l|l|l|l}
			$\widehat{\ell}$ & $m^*$ & $N_{it}^*$ & $S^*_\ell$ & $N_p$ \\
			\hline
			$2$ &  $229$ & $8$ & $4.09$ & $132$ \\
			\hline
			$3$ & $42$ & $8$ & $6.41$ & $715$
		\end{tabular}\subcaption{$n=80$, $N_t=30,018$.}
		\end{minipage}
	\end{table}

    \begin{table}[h!]
		\caption{Speedup estimates for MP1, \cref{eq:mp1}, for different values of $n$ and $\widehat{\ell}$, using gradient descent operator \cref{eq:gd} on the fine grid and proximal point operator \cref{eq:ppm} on the coarse grid. Here, $\widehat{m}^*$ and $\widehat{S}^{*}$ are computed based on a hypothetical lower bound for $\frac{t_c}{t_f}$ given by $\widehat{\alpha}=1$. $N_p$ denotes the number of processors required to achieve optimal speedup.}\label{table:mp1_speedup_ideal}
        \begin{minipage}[c]{0.33\linewidth}
		\centering
	   \begin{tabular}{l|l|l|l|l}
		$\widehat{\ell}$ & $\widehat{m}^{*}$ & $\widehat{N_{it}}^*$ & $\widehat{S}^{*}_\ell$ & $N_p$\\
		\hline
		$2$ & $65$ & $8$ & $4.07$ & $131$\\
		\hline
		$3$ & $22$ & $8$ & $8.33$ & $387$
		\end{tabular}
		\subcaption{$n=40$, $N_t=8,503$.}
        \end{minipage}\begin{minipage}[c]{0.33\linewidth}
	\centering   
       \begin{tabular}{l|l|l|l|l}
		$\widehat{\ell}$ & $\widehat{m}^{*}$ & $\widehat{N_{it}}^*$ & $\widehat{S}^{*}_\ell$ & $N_p$ \\
		\hline
		$2$ & $95$ & $8$ & $5.93$ & $190$ \\
		\hline
		$3$ & $28$ & $8$  & $13.84$ & $645$
 		\end{tabular}\subcaption{$n=60$, $N_t=18,050$.}
        \end{minipage}\begin{minipage}[c]{0.33\linewidth}
	    \centering   
           \begin{tabular}{l|l|l|l|l}
			$\ell$ & $\widehat{m}^{*}$ & $\widehat{N_{it}}^*$ & $\widehat{S}^{*}_\ell$ & $N_p$ \\
			\hline
			$2$ & $123$ & $8$ & $7.65$  & $245$\\
			\hline
			$3$ & $33$ & $8$ & $19.48$ & $910$
		\end{tabular}\subcaption{$n=80$, $N_t=30,018$.}
		\end{minipage}
	\end{table}

	\begin{table}[h!]
		\caption{Speedup estimates MP2-1D, \cref{eq:mp2}, for different values of $n$ and $\widehat{\ell}$, using proximal gradient operator \cref{eq:pg} on the fine grid and alternating proximals operator \cref{eq:apm} on the coarse grid. Here, $m^*$ and $S^*$ are computed based on the numerically estimated $\frac{t_c}{t_f}$ values given by $\alpha=2.4$ for $n=256$, and $\alpha=3.9$ for $n=512$. $N_p$ denotes the number of processors required to achieve optimal speedup.}\label{table:1d_mp2_speedup}
        \begin{minipage}[c]{0.5\linewidth}
		\centering
		\begin{tabular}{l|l|l|l|l}
			$\widehat{\ell}$ & $m^*$ & $N_{it}^*$ & $S^*_\ell$ & $N_p$ \\
			\hline
			$2$ & $363$ & $10$ & $7.56$ & $303$\\
			\hline
			$3$ & $60$ & $9$ & $19.52$ & $1832$
		\end{tabular}\subcaption{$n=256$, $N_t=109,888$.}
        \end{minipage}
        \begin{minipage}[c]{0.5\linewidth}
        \centering
		\begin{tabular}{l|l|l|l|l}
			$\widehat{\ell}$ & $m^*$ & $N_{it}^*$ & $S^*_\ell$ & $N_p$  \\
			\hline
			$2$ & $849$ & $10$ & $10.88$ & $436$\\
			\hline
			$3$ & $100$ & $9$ & $27.12$ & $3697$
		\end{tabular}\subcaption{$n=512$, $N_t=369,612$.}
        \end{minipage}
	\end{table}

    \begin{table}[h!]
		\caption{Speedup estimates for MP2-1D, \cref{eq:mp2}, for different values of $n$ and $\widehat{\ell}$, using the proximal gradient operator \cref{eq:pg} on the fine grid and the alternating proximals operator \cref{eq:apm} on the coarse grid. Here, $\widehat{m}^*$ and $\widehat{S}^{*}$ are computed based on an hypothetical lower bound for $\frac{t_c}{t_f}$ given by $\widehat{\alpha}=1$. $N_p$ denotes the number of processors required to achieve optimal speedup.}\label{table:1d_mp2_speedup_ideal}
        \begin{minipage}[c]{0.5\linewidth}
		\centering
	   \begin{tabular}{l|l|l|l|l}
		$\widehat{\ell}$ & $\widehat{m}^{*}$ & $\widehat{N_{it}}^*$ & $\widehat{S}^{*}_\ell$ & $N_p$ \\
		\hline
		$2$ & $234$  & $10$ & $11.72$ & $470$\\
		\hline
		$3$ & $52$  & $9$  & $40.61$ & $2114$
		\end{tabular}\subcaption{$n=256$, $N_t=109,888$.}
        \end{minipage}
        \begin{minipage}[c]{0.5\linewidth}
        \centering
		\begin{tabular}{l|l|l|l|l}
		$\widehat{\ell}$ &  $\widehat{m}^{*}$ & $\widehat{N_{it}}^*$ & $\widehat{S}^{*}_\ell$ & $N_p$ \\
		\hline
		$2$ & $430$ & $10$  & $21.49$ & $860$\\
		\hline
		$3$ & $77$ & $10$  & $82.62$ & $4801$
		\end{tabular}\subcaption{$n=512$, $N_t=369,612$.}
		\end{minipage}
	\end{table}
	
		\begin{table}[h!]
		\caption{Speedup estimates for MP2-2D, \cref{eq:mp2}, for different values of $n$ and $\widehat{\ell}$, using proximal gradient operator \cref{eq:pg} on the fine grid and alternating proximals operator \cref{eq:apm} on the coarse grid. Here, $m^*$ and $S^*$ are computed based on the numerically estimated  $\frac{t_c}{t_f}$ values given by $\alpha=13$ for $n^2=64^2$, and $\alpha=7$ for $n^2=128^2$. $N_p$ denotes the number of processors required to achieve optimal speedup.}\label{table:2d_mp2_speedup}
        \begin{minipage}[c]{0.5\linewidth}
        \centering
		\begin{tabular}{l|l|l|l|l}
			$\widehat{\ell}$ & $m^*$ & $N_{it}^*$ & $S^*_\ell$ & $N_p$ \\
			\hline
			$2$ & $294$  & $9$ &  $1.25$ & $46$ \\
			\hline
			$3$ & $46$ & $10$ &  $0.67$ & $289$
		\end{tabular}\subcaption{$n^2=64^2$, $N_t=13,256$.}
		\end{minipage}
                \begin{minipage}[c]{0.5\linewidth}
        \centering
		\begin{tabular}{l|l|l|l|l}
			$\widehat{\ell}$ & $m^*$ & $N_{it}^*$ & $S^*_\ell$ & $N_p$ \\
			\hline
			$2$ & $403$  & $9$ &  $3.19$ & $116$\\
			\hline
			$3$ & $58$ & $10$ & $3.24$ & $801$
		\end{tabular}\subcaption{$n^2=128^2$, $N_t=46,431$.}
		\end{minipage}
	\end{table}

    \begin{table}[h!]
		\caption{Speedup estimates for MP2-2D, \cref{eq:mp2}, for different values of $n$ and $\widehat{\ell}$, using the proximal gradient operator \cref{eq:pg} on the fine grid and the alternating proximals operator \cref{eq:apm} on the coarse grid. Here, $\widehat{m}^*$ and $\widehat{S}^{*}$ are computed based on an hypothetical lower bound for $\frac{t_c}{t_f}$ given by $\widehat{\alpha}=1$. $N_p$ denotes the number of processors required to achieve optimal speedup.}\label{table:2d_mp2_speedup_ideal}
        \begin{minipage}[c]{0.5\linewidth}
        \centering
		\begin{tabular}{l|l|l|l|l}
		$\widehat{\ell}$ & $\widehat{m}^{*}$ & $\widehat{N_{it}}^*$ & $\widehat{S}^{*}_\ell$ & $N_p$\\
		\hline
		$2$ & $81$  & $10$ & $4.07$ & $164$ \\
		\hline
		$3$ & $25$ & $9$ &  $10.07$ & $531$
		\end{tabular}\subcaption{$n^2=64^2$, $N_t=13,256$.}
		\end{minipage}
        \begin{minipage}[c]{0.5\linewidth}
        \centering
		\begin{tabular}{l|l|l|l|l}
		$\widehat{\ell}$ & $\widehat{m}^{*}$ & $\widehat{N_{it}}^*$ & $\widehat{S}^{*}_\ell$ & $N_p$\\
		\hline
		$2$ & $152$ & $10$  & $7.61$ & $306$ \\
		\hline
		$3$ & $39$ & $9$ &  $22.87$ & $1191$
		\end{tabular}\subcaption{$n^2=128^2$, $N_t=46,431$.}
		\end{minipage}
	\end{table}

	First, we note that 2-level MGRIT is predicted to provide speedup over sequential solves in all of our estimates. For most cases, the estimated speedups for $\ell=3$ were considerable improvements over the $\ell=2$ speedup.
    For all cases, increasing $N$ led to improved speedup. As the sizes $N$ of the Laplacian matrices \cref{eq:1d_laplacian,eq:2d_laplacian} increase, more sequential iterations are needed until convergence is achieved, which in turn increases the size of the MGRIT time grids. Optimal speedup formula \cref{eq:s2_opt} shows that larger grids often benefit from more parallelization. On the other hand, increasing $N$ is also correlated with an increased $\frac{t_c}{t_f}$ ratio, due to the cost of the implicit coarse-grid time-stepping operations increasing relatively more than that of the explicit fine-grid ones. However, in our experiments the larger $\frac{t_c}{t_f}$ ratios were offset by the increase in $N_t$, such that speedup was improved with larger $N$.
    
    The results in \cref{table:2d_mp2_speedup} are exceptions to the above. For the $n^2=64^2$ experiment, the estimated $\frac{t_f}{t_c}$ ratio of $\alpha=13$ was very large. This is expected, since we are using 2D Laplacians \cref{eq:2d_laplacian} instead of 1D ones \cref{eq:1d_laplacian}, thus increasing the relative cost of implicit steps. In fact, $\alpha$ was so large that the 2-level solve is predicted to perform better than the 3-level solve, likely since using 3 levels requires an additional post-relaxation step on the $\ell=2$ coarse grid. However, when increasing the problem size to $n^2=128^2$, the $\frac{t_c}{t_f}$ ratio was reduced to $\alpha=7$. To understand why this happens, recall that the fine-grid proximal gradient operator \cref{eq:pg} is $\Phi=P_{sg}\circ G_{sf}$, while the coarse-grid alternating proximals operator \cref{eq:apm} is $\Phi_\Delta=P_{sg}\circ P_{sf}$, with $s=\frac{1}{L}$ for $\Phi$ and $s=\frac{m}{L}$ for $\Phi_\Delta$. As $N$ increases, the cost of the $P_{sg}$ step in both $\Phi$ and $\Phi_\Delta$ accounts for a larger share of both $t_c$ and $t_f$, such that the $\frac{t_c}{t_f}$ ratio is reduced. While still large, $\frac{t_c}{t_f}$ was reduced enough that the 3-level solve is estimated to provide a slight speedup improvement over the 2-level solve.

    In particular, \cref{table:1d_mp2_speedup} and \cref{table:2d_mp2_speedup} exemplify the conditions for which our method performs best and worse, respectively; in the former, we have large $N_t$ and small $\alpha$, while the opposite is true for the latter. The implicit steps in our MGRIT implementation were calculated using LU-factorization with saved $L$ and $U$ factors to avoid having to repeatedly invert large matrices directly; improving this calculation can bring the $\frac{t_c}{t_f}$ ratio to lower, more beneficial levels, hence bringing speedup estimates closer to the hypothetical $\frac{t_c}{t_f}=1$ results in \cref{table:mp1_speedup_ideal}, \cref{table:1d_mp2_speedup_ideal} and \cref{table:2d_mp2_speedup_ideal}. While not explored here, we expect speedups to potentially improve with more levels; this is supported by \cref{table:1d_mp1_conv}, \cref{table:1d_mp2_conv} and \cref{table:2d_mp2_conv}, in which we observe that the $N_{it}$ remains fairly consistent across for all possible values $\ell$.
    Finally, increasing the resolution of MP2-2D to $n^2=256^2$ or $n^2=512^2$ to match the tests for MP2-1D, would dramatically increase $N_t$, which would also be expected to substantially increase the speedup.

	\section{Issues with momentum-accelerated algorithms}\label{sec:nesterov}
	
	So far, we have only worked with standalone versions of gradient-based methods such as gradient descent, \cref{eq:gd}, and proximal gradient descent, \cref{eq:pg}. However, in practice they are often replaced by accelerated variants; for example, gradient descent is often replaced by Nesterov's accelerated gradient descent \cite{Nesterov1983}. Starting from $\mathbf{u}_0 \in \mathbb{R}^N$ and letting $0 < s < \frac{2}{L}$, Nesterov's accelerated gradient sequence  $\{\mathbf{u}_k,\mathbf{v}_k\}$, given by
	\begin{equation}\label{eq:nest_iter}
		\begin{aligned}
			& \mathbf{u}_k=\mathbf{v}_{k-1}-s \nabla f\left(\mathbf{v}_{k-1}\right), \\
			& \mathbf{v}_k=\mathbf{u}_k+\beta_k\left(\mathbf{u}_k-\mathbf{u}_{k-1}\right),
		\end{aligned}
	\end{equation}
	with momentum parameter $\beta_k=\frac{k-1}{k+2}$, converges to a minimizer of $f$ with a convergence rate of $\mathcal{O}(\frac{1}{k^2})$ for convex $f$ with Lipschitz gradient \cite{Nesterov2004}. The connection between optimization method and ODE still holds, although not as trivially as in \cref{subsec:opt_disc_ode}. In \cite{Su2016}, a second-order ODE was derived which is the exact limit of Nesterov's scheme as the time step $s$ approaches $0$:
	\begin{equation}\label{eq:nest_ode}
			\frac{\mathrm{d}^2}{\mathrm{d}t^2}\mathbf{u}(t)+\frac{3}{t} \frac{\mathrm{d}}{\mathrm{d}t}\mathbf{u}(t)+\nabla f(\mathbf{u}(t))=0,
	\end{equation}
	for $t>0$ with initial conditions $\mathbf{u}(0)=\mathbf{u}_0$ and $\frac{\mathrm{d}}{\mathrm{d}t}\mathbf{u}(0)=0$, where $\mathbf{u}_0$ is the initial guess. The time parameter in this ODE is related to the step size $s$ in Nesterov's scheme via $t \approx k \sqrt{s}$.
	
	Initially, we expected that MGRIT in standard form with implicit discretization on the coarse grid, as explained above, would work well for MP1, \cref{eq:mp1}, with Nesterov's accelerated gradient operator \cref{eq:nest_ode} as the fine-grid operator and a similar momentum-accelerated proximal point operator \cref{eq:ppm} as the coarse-grid operator, given the positive results obtained with gradient descent \cref{eq:gd}. However, in our preliminary work, MGRIT has failed to converge except for when using very small constant values of $\beta=\epsilon>0$, in which case Nesterov's scheme behaves almost exactly like gradient descent. Unlike gradient descent, which we have observed to behave similarly to a parabolic PDE in \cref{subsec:mgrit_ode_link}, the same does not apply for Nesterov's accelerated gradient descent as explained in the following remark.
    
	\begin{remark}[Nesterov's accelerated gradient ODE and hyperbolic PDE connection]\label{eq:nest_hyper}
		Consider Nesterov's accelerated gradient ODE \cref{eq:nest_ode} applied to MP1, \cref{eq:mp1}, such that $\nabla f(\mathbf{u}(t))=A_1\mathbf{u}$, where $A_1$ is the 1D Laplacian matrix in \cref{eq:1d_laplacian} with $(A_1\mathbf{u})_i \approx - \frac{\partial^2}{\partial x^2}u(x_i,t)$. Then, we have
		\begin{equation}
			\frac{\partial^2}{\partial t^2}u(x,t)+\frac{3}{t}\frac{\partial}{\partial t}u(x,t)-\frac{\partial^2}{\partial x^2}u(x,t)=0,
		\end{equation}
		which is a hyperbolic PDE. 
	\end{remark}

	
	This behavior similar to a hyperbolic PDE, which MGRIT in standard form is known to struggle with, might be a reason for why Nesterov's accelerated gradient did not work with MGRIT using the standard algorithmic approach operators described in \cref{sec:p-in-t}.
    On the other hand, it has recently been shown that MGRIT can be made to work for hyperbolic PDEs by carefully considering non-standard specialized coarse operators:
    research on on modified semi-Lagrangian coarse-grid operators for advection-dominated and hyperbolic PDEs \cite{DeSterck2019,DeSterck2021,DeSterck2023,DeSterck2023_2,Krzysik2025} has shown progress in adapting MGRIT to solve hyperbolic PDEs, which indicates a promising direction for further research on applying MGRIT to momentum-accelerated optimization algorithms. 
		
	\section{Conclusions}\label{sec:concl}
	
	We have presented a framework for solving optimization problems all-at-once, in the sense that sequential solves are replaced by a parallel solver in which the sequential iterations are computed simultaneously.
    Using the MGRIT algorithm to parallelize over iteration of the optimization method instead of over time, we demonstrated a multilevel parallel-in-iteration method that can converge in a small number of iterations, and also provides considerable estimated speedup for both smooth and nonsmooth problems. Further improvements to this framework and the consideration of other model problems can be investigated by the continued exploration of the connection between optimization algorithms and ODEs, especially those related to diffusion-dominated problems for which MGRIT performs well. \ptxt{Among the problems of our interest, we list constrained optimization problems using ADMM, minimizations of nonquadratic functions (such as the cross-entropy loss function from logistic regression) and applications of stochastic approximations such as stochastic gradient descent, to name a few.}
    We established that in order to extend this framework to momentum-accelerated algorithms such as Nesterov's accelerated gradient descent and FISTA, we will likely need to address MGRIT's difficulties with handling hyperbolic PDEs. Adapting MGRIT with this intent is an intricate problem-specific endeavor, which we leave as a suggestion for future work.
	

    \bibliographystyle{etna}
    \bibliography{refs.bib}

@article{Ang2024,
	title = {{MGProx}: A Nonsmooth Multigrid Proximal Gradient Method with Adaptive Restriction for Strongly Convex Optimization},
	volume = {34},
	ISSN = {1095-7189},
	url = {http://dx.doi.org/10.1137/23M1552140},
	DOI = {10.1137/23m1552140},
	number = {3},
	journal = {SIAM Journal on Optimization},
	publisher = {Society for Industrial & Applied Mathematics (SIAM)},
	author = {Ang,  Andersen and De Sterck,  Hans and Vavasis,  Stephen},
	year = {2024},
	month = aug,
	pages = {2788–2820}
}

@book{Angrist2008,
  title = {Mostly Harmless Econometrics: An Empiricist’s Companion},
  ISBN = {9780691120348},
  url = {http://dx.doi.org/10.2307/j.ctvcm4j72},
  DOI = {10.2307/j.ctvcm4j72},
  publisher = {Princeton University Press},
  address = {Princeton, USA},
  author = {Angrist,  Joshua D. and Pischke,  J\"{o}rn-Steffen},
  year = {2008},
  month = dec 
}

@article{Attouch2016,
	title = {Fast convergence of inertial dynamics and algorithms with asymptotic vanishing viscosity},
	volume = {168},
	ISSN = {1436-4646},
	url = {http://dx.doi.org/10.1007/s10107-016-0992-8},
	DOI = {10.1007/s10107-016-0992-8},
	number = {1–2},
	journal = {Mathematical Programming},
	publisher = {Springer Science and Business Media LLC},
	author = {Attouch,  Hedy and Chbani,  Zaki and Peypouquet,  Juan and Redont,  Patrick},
	year = {2016},
	month = mar,
	pages = {123–175}
}

@book{Bauschke2023,
	title = {An Introduction to Convexity,  Optimization,  and Algorithms},
	ISBN = {9781611977806},
	url = {http://dx.doi.org/10.1137/1.9781611977806},
	DOI = {10.1137/1.9781611977806},
	publisher = {Society for Industrial and Applied Mathematics},
        address = {Philadelphia, USA},
	author = {Bauschke,  Heinz H. and Moursi,  Walaa M.},
	year = {2023},
	month = jan 
}

@article{Beck2009,
	title = {A Fast Iterative Shrinkage-Thresholding Algorithm for Linear Inverse Problems},
	volume = {2},
	ISSN = {1936-4954},
	url = {http://dx.doi.org/10.1137/080716542},
	DOI = {10.1137/080716542},
	number = {1},
	journal = {SIAM Journal on Imaging Sciences},
	publisher = {Society for Industrial & Applied Mathematics (SIAM)},
	author = {Beck,  Amir and Teboulle,  Marc},
	year = {2009},
	month = jan,
	pages = {183–202}
}

@article{Brandt1983,
	title = {Multigrid Algorithms for the Solution of Linear Complementarity Problems Arising from Free Boundary Problems},
	volume = {4},
	ISSN = {2168-3417},
	url = {http://dx.doi.org/10.1137/0904046},
	DOI = {10.1137/0904046},
	number = {4},
	journal = {SIAM Journal on Scientific and Statistical Computing},
	publisher = {Society for Industrial & Applied Mathematics (SIAM)},
	author = {Brandt,  Achi and Cryer,  Colin W.},
	year = {1983},
	month = dec,
	pages = {655–684}
}

@article{Caffarelli1998,
	title = {The obstacle problem revisited},
	volume = {4},
	ISSN = {1531-5851},
	url = {http://dx.doi.org/10.1007/BF02498216},
	DOI = {10.1007/bf02498216},
	number = {4–5},
	journal = {The Journal of Fourier Analysis and Applications},
	publisher = {Springer Science and Business Media LLC},
	author = {Caffarelli,  L. A.},
	year = {1998},
	month = jul,
	pages = {383–402}
}

@article{Chambolle2016,
  title = {An introduction to continuous optimization for imaging},
  volume = {25},
  ISSN = {1474-0508},
  url = {http://dx.doi.org/10.1017/S096249291600009X},
  DOI = {10.1017/s096249291600009x},
  journal = {Acta Numerica},
  publisher = {Cambridge University Press (CUP)},
  author = {Chambolle,  Antonin and Pock,  Thomas},
  year = {2016},
  month = may,
  pages = {161–319}
}

@article{DeSterck2021,
	title = {Optimizing multigrid reduction‐in‐time and Parareal coarse‐grid operators for linear advection},
	volume = {28},
	ISSN = {1099-1506},
	url = {http://dx.doi.org/10.1002/nla.2367},
	DOI = {10.1002/nla.2367},
	number = {4},
	journal = {Numerical Linear Algebra with Applications},
	publisher = {Wiley},
	author = {De Sterck,  Hans and Falgout,  Robert D. and Friedhoff,  Stephanie and Krzysik,  Oliver A. and MacLachlan,  Scott P.},
	year = {2021},
	month = mar 
}

@article{DeSterck2023_2,
	title = {Fast Multigrid Reduction-in-Time for Advection via Modified Semi-{L}agrangian Coarse-Grid Operators},
	volume = {45},
	ISSN = {1095-7197},
	url = {http://dx.doi.org/10.1137/22M1486522},
	DOI = {10.1137/22m1486522},
	number = {4},
	journal = {SIAM Journal on Scientific Computing},
	publisher = {Society for Industrial & Applied Mathematics (SIAM)},
	author = {De Sterck,  Hans and Falgout,  Robert D. and Krzysik,  Oliver A.},
	year = {2023},
	month = jul,
	pages = {A1890–A1916}
}

@article{DeSterck2023,
	title = {Efficient Multigrid Reduction-in-Time for Method-of-Lines Discretizations of Linear Advection},
	volume = {96},
	ISSN = {1573-7691},
	url = {http://dx.doi.org/10.1007/s10915-023-02223-4},
	DOI = {10.1007/s10915-023-02223-4},
	number = {1},
	journal = {Journal of Scientific Computing},
	publisher = {Springer Science and Business Media LLC},
	author = {De Sterck,  H. and Falgout,  R. D. and Krzysik,  O. A. and Schroder,  J. B.},
	year = {2023},
	month = may 
}

@article{DeSterck2019,
	title = {Convergence analysis for parallel‐in‐time solution of hyperbolic systems},
	volume = {27},
	ISSN = {1099-1506},
	url = {http://dx.doi.org/10.1002/nla.2271},
	DOI = {10.1002/nla.2271},
	number = {1},
	journal = {Numerical Linear Algebra with Applications},
	publisher = {Wiley},
	author = {De Sterck,  Hans and Friedhoff,  Stephanie and Howse,  Alexander J. M. and MacLachlan,  Scott P.},
	year = {2019},
	month = nov 
}

@article{DeSterck2024,
	title = {Multigrid Reduction‐In‐Time Convergence for Advection Problems: A {F}ourier Analysis Perspective},
	volume = {32},
	ISSN = {1099-1506},
	url = {http://dx.doi.org/10.1002/nla.2593},
	DOI = {10.1002/nla.2593},
	number = {1},
	journal = {Numerical Linear Algebra with Applications},
	publisher = {Wiley},
	author = {De Sterck,  H. and Friedhoff,  S. and Krzysik,  O. A. and MacLachlan,  S. P.},
	year = {2024},
	month = oct 
}

@article{Dobrev2017,
	title = {Two-Level Convergence Theory for Multigrid Reduction in Time ({MGRIT})},
	volume = {39},
	ISSN = {1095-7197},
	url = {http://dx.doi.org/10.1137/16M1074096},
	DOI = {10.1137/16m1074096},
	number = {5},
	journal = {SIAM Journal on Scientific Computing},
	publisher = {Society for Industrial \& Applied Mathematics (SIAM)},
	author = {Dobrev,  V. A. and Kolev,  Tz. and Petersson,  N. A. and Schroder,  J. B.},
	year = {2017},
	month = jan,
	pages = {S501–S527}
}

@article{Emmett2012,
  title = {Toward an efficient parallel in time method for partial differential equations},
  volume = {7},
  ISSN = {1559-3940},
  url = {http://dx.doi.org/10.2140/camcos.2012.7.105},
  DOI = {10.2140/camcos.2012.7.105},
  number = {1},
  journal = {Communications in Applied Mathematics and Computational Science},
  publisher = {Mathematical Sciences Publishers},
  author = {Emmett,  Matthew and Minion,  Michael},
  year = {2012},
  month = mar,
  pages = {105–132}
}

@article{Falgout2014,
	title = {Parallel Time Integration with Multigrid},
	volume = {36},
	ISSN = {1095-7197},
	url = {http://dx.doi.org/10.1137/130944230},
	DOI = {10.1137/130944230},
	number = {6},
	journal = {SIAM Journal on Scientific Computing},
	publisher = {Society for Industrial \& Applied Mathematics (SIAM)},
	author = {Falgout,  R. D. and Friedhoff,  S. and Kolev,  Tz. V. and MacLachlan,  S. P. and Schroder,  J. B.},
	year = {2014},
	month = jan,
	pages = {C635–C661}
}

@article{Falgout2017,
  title = {Multigrid Reduction in Time for Nonlinear Parabolic Problems: A Case Study},
  volume = {39},
  ISSN = {1095-7197},
  url = {http://dx.doi.org/10.1137/16M1082330},
  DOI = {10.1137/16m1082330},
  number = {5},
  journal = {SIAM Journal on Scientific Computing},
  publisher = {Society for Industrial & Applied Mathematics (SIAM)},
  author = {Falgout,  R. D. and Manteuffel,  T. A. and O’Neill,  B. and Schroder,  J. B.},
  year = {2017},
  month = jan,
  pages = {S298–S322}
}

@article{Friedlander2008,
  title = {Exact Regularization of Convex Programs},
  volume = {18},
  ISSN = {1095-7189},
  url = {http://dx.doi.org/10.1137/060675320},
  DOI = {10.1137/060675320},
  number = {4},
  journal = {SIAM Journal on Optimization},
  publisher = {Society for Industrial & Applied Mathematics (SIAM)},
  author = {Friedlander,  Michael P. and Tseng,  Paul},
  year = {2008},
  month = jan,
  pages = {1326–1350}
}

@inbook{Gander2015,
	title = {50 Years of Time Parallel Time Integration},
	ISBN = {9783319233215},
	ISSN = {2191-3048},
	url = {http://dx.doi.org/10.1007/978-3-319-23321-5_3},
	DOI = {10.1007/978-3-319-23321-5_3},
	booktitle = {Multiple Shooting and Time Domain Decomposition Methods},
	publisher = {Springer International Publishing},
        address = {Cham, Switzerland},
	author = {Gander,  Martin J.},
	year = {2015},
	pages = {69–113}
}

@article{Gander2020,
	title = {{PARAOPT}: A Parareal Algorithm for Optimality Systems},
	volume = {42},
	ISSN = {1095-7197},
	url = {http://dx.doi.org/10.1137/19M1292291},
	DOI = {10.1137/19m1292291},
	number = {5},
	journal = {SIAM Journal on Scientific Computing},
	publisher = {Society for Industrial & Applied Mathematics (SIAM)},
	author = {Gander,  Martin J. and Kwok,  Felix and Salomon,  Julien},
	year = {2020},
	month = jan,
	pages = {A2773–A2802}
}

@article{Gander2007,
	title = {Analysis of the Parareal Time‐Parallel Time‐Integration Method},
	volume = {29},
	ISSN = {1095-7197},
	url = {http://dx.doi.org/10.1137/05064607X},
	DOI = {10.1137/05064607x},
	number = {2},
	journal = {SIAM Journal on Scientific Computing},
	publisher = {Society for Industrial \& Applied Mathematics (SIAM)},
	author = {Gander,  Martin J. and Vandewalle,  Stefan},
	year = {2007},
	month = jan,
	pages = {556–578}
}

@misc{Gao2017,
Author = {Bolin Gao and Lacra Pavel},
Title = {On the Properties of the Softmax Function with Application in Game Theory and Reinforcement Learning},
Year = {2017},
Eprint = {arXiv:1704.00805},
note = { \\ \url{https://arxiv.org/abs/1704.00805}}
}

@article{Golub1999,
  title = {Tikhonov Regularization and Total Least Squares},
  volume = {21},
  ISSN = {1095-7162},
  url = {http://dx.doi.org/10.1137/S0895479897326432},
  DOI = {10.1137/s0895479897326432},
  number = {1},
  journal = {SIAM Journal on Matrix Analysis and Applications},
  publisher = {Society for Industrial & Applied Mathematics (SIAM)},
  author = {Golub,  Gene H. and Hansen,  Per Christian and O’Leary,  Dianne P.},
  year = {1999},
  month = jan,
  pages = {185–194}
}

@article{Hahne2023,
	title = {Parallel-in-time optimization of induction motors},
	volume = {13},
	ISSN = {2190-5983},
	url = {http://dx.doi.org/10.1186/s13362-023-00134-5},
	DOI = {10.1186/s13362-023-00134-5},
	number = {1},
	journal = {Journal of Mathematics in Industry},
	publisher = {Springer Science and Business Media LLC},
        address = {Cham, Switzerland},
	author = {Hahne,  Jens and Polenz,  Bj\"{o}rn and Kulchytska-Ruchka,  Iryna and Friedhoff,  Stephanie and Ulbrich,  Stefan and Sch\"{o}ps,  Sebastian},
	year = {2023},
	month = jun 
}

@article{Howse2019,
	title = {Parallel-In-Time Multigrid with Adaptive Spatial Coarsening for The Linear Advection and Inviscid Burgers Equations},
	volume = {41},
	ISSN = {1095-7197},
	url = {http://dx.doi.org/10.1137/17M1144982},
	DOI = {10.1137/17m1144982},
	number = {1},
	journal = {SIAM Journal on Scientific Computing},
	publisher = {Society for Industrial & Applied Mathematics (SIAM)},
	author = {Howse,  Alexander J. and Sterck,  Hans De and Falgout,  Robert D. and MacLachlan,  Scott and Schroder,  Jacob},
	year = {2019},
	month = jan,
	pages = {A538–A565}
}

@article{Krzysik2025,
author = {Krzysik, O. A. and De Sterck, H. and Falgout, R. D. and Schroder, J. B.},
title = {Parallel-in-Time Solution of Hyperbolic {PDE} Systems via Characteristic-Variable Block Preconditioning},
journal = {SIAM Journal on Scientific Computing},
volume = {0},
number = {0},
pages = {S337-S363},
year = {0},
url = {http://dx.doi.org/10.1137/24M1673310}
}

@article{Lions2001,
	title = {Resolution d'{EDP} par un schema en temps parareel},
	volume = {332},
	ISSN = {0764-4442},
	url = {http://dx.doi.org/10.1016/S0764-4442(00)01793-6},
	DOI = {10.1016/s0764-4442(00)01793-6},
	number = {7},
	journal = {Comptes Rendus de l’Académie des Sciences - Series I - Mathematics},
	publisher = {Elsevier BV},
	author = {Lions,  Jacques-Louis and Maday,  Yvon and Turinici,  Gabriel},
	year = {2001},
	month = apr,
	pages = {661–668}
}

@article{Lions1979,
	ISSN = {00361429},
	URL = {http://www.jstor.org/stable/2156649},
	author = {P. L. Lions and B. Mercier},
	journal = {SIAM Journal on Numerical Analysis},
	number = {6},
	pages = {964--979},
	publisher = {Society for Industrial and Applied Mathematics},
	title = {Splitting Algorithms for the Sum of Two Nonlinear Operators},
	urldate = {2024-04-06},
	volume = {16},
	year = {1979}
}

@inbook{Maday2013,
	title = {Parareal in Time Intermediate Targets Methods for Optimal Control Problems},
	ISBN = {9783034806312},
	ISSN = {2296-6072},
	url = {http://dx.doi.org/10.1007/978-3-0348-0631-2_5},
	DOI = {10.1007/978-3-0348-0631-2_5},
	booktitle = {Control and Optimization with PDE Constraints},
	publisher = {Springer Basel},
        address = {Basel, Switzerland},
	author = {Maday,  Yvon and Riahi,  Mohamed-Kamel and Salomon,  Julien},
	year = {2013},
	pages = {79–92}
}

@article{Mandel1984,
	title = {A multilevel iterative method for symmetric,  positive definite linear complementarity problems},
	volume = {11},
	ISSN = {1432-0606},
	url = {http://dx.doi.org/10.1007/BF01442171},
	DOI = {10.1007/bf01442171},
	number = {1},
	journal = {Applied Mathematics \& Optimization},
	publisher = {Springer Science and Business Media LLC},
	author = {Mandel,  Jan},
	year = {1984},
	month = feb,
	pages = {77–95}
}

@article{Mangasarian1985,
  title = {Sufficiency of Exact Penalty Minimization},
  volume = {23},
  ISSN = {1095-7138},
  url = {http://dx.doi.org/10.1137/0323003},
  DOI = {10.1137/0323003},
  number = {1},
  journal = {SIAM Journal on Control and Optimization},
  publisher = {Society for Industrial & Applied Mathematics (SIAM)},
  author = {Mangasarian,  O. L.},
  year = {1985},
  month = jan,
  pages = {30–37}
}

@article{Minion2010,
	author = {Minion,  Michael},
	title = {A hybrid parareal spectral deferred corrections method},
	volume = {5},
	ISSN = {1559-3940},
	url = {http://dx.doi.org/10.2140/camcos.2010.5.265},
	DOI = {10.2140/camcos.2010.5.265},
	number = {2},
	journal = {Communications in Applied Mathematics and Computational Science},
	publisher = {Mathematical Sciences Publishers},
	year = {2010},
	month = dec,
	pages = {265–301}
}

@article{Nesterov1983,
	title={A method for solving the convex programming problem with convergence rate $O(1/k^2)$},
	author={Nesterov,  Yurii},
	journal={Proceedings of the USSR Academy of Sciences},
	year={1983},
	volume={269},
	pages={543-547},
	url={https://api.semanticscholar.org/CorpusID:145918791}
}

@book{Nesterov2004,
	title = {Introductory Lectures on Convex Optimization},
	ISBN = {9781441988539},
	ISSN = {1384-6485},
	url = {http://dx.doi.org/10.1007/978-1-4419-8853-9},
	DOI = {10.1007/978-1-4419-8853-9},
	journal = {Applied Optimization},
	publisher = {Springer US},
	author = {Nesterov,  Yurii},
	year = {2004}
}

@article{Ong2020,
	title = {Applications of time parallelization},
	volume = {23},
	ISSN = {1433-0369},
	url = {http://dx.doi.org/10.1007/s00791-020-00331-4},
	DOI = {10.1007/s00791-020-00331-4},
	number = {1–4},
	journal = {Computing and Visualization in Science},
	publisher = {Springer Science and Business Media LLC},
	author = {Ong,  Benjamin W. and Schroder,  Jacob B.},
	year = {2020},
	month = sep 
}

@article{ODonoghue2013,
	title = {Adaptive Restart for Accelerated Gradient Schemes},
	volume = {15},
	ISSN = {1615-3383},
	url = {http://dx.doi.org/10.1007/s10208-013-9150-3},
	DOI = {10.1007/s10208-013-9150-3},
	number = {3},
	journal = {Foundations of Computational Mathematics},
	publisher = {Springer Science and Business Media LLC},
	author = {O’Donoghue,  Brendan and Candès,  Emmanuel},
	year = {2013},
	month = jul,
	pages = {715–732}
}

@article{Parikh2014,
  title = {Proximal Algorithms},
  volume = {1},
  ISSN = {2167-3918},
  url = {http://dx.doi.org/10.1561/2400000003},
  DOI = {10.1561/2400000003},
  number = {3},
  journal = {Foundations and Trends{\textregistered} in Optimization},
  publisher = {Now Publishers},
  author = {Parikh,  Neal},
  year = {2014},
  pages = {127–239}
}

@article{Ries1983,
	title = {A note on {MGR} methods},
	volume = {49},
	ISSN = {0024-3795},
	url = {http://dx.doi.org/10.1016/0024-3795(83)90091-5},
	DOI = {10.1016/0024-3795(83)90091-5},
	journal = {Linear Algebra and its Applications},
	publisher = {Elsevier BV},
	author = {Ries,  Manfred and Trottenberg,  Ulrich and Winter,  Gerd},
	year = {1983},
	month = feb,
	pages = {1–26}
}

@BOOK{Rodrigues2014,
	title     = "Obstacle problems in mathematical physics",
	author    = "Rodrigues, J-F",
	publisher = "North-Holland",
	series    = "North-Holland Mathematics Studies",
	month     =  may,
	year      =  2014
}

@inproceedings{Romero2022,
  title = {ODE Discretization Schemes as Optimization Algorithms},
  url = {http://dx.doi.org/10.1109/CDC51059.2022.9992691},
  DOI = {10.1109/cdc51059.2022.9992691},
  booktitle = {2022 IEEE 61st Conference on Decision and Control (CDC)},
  publisher = {IEEE},
  author = {Romero,  Orlando and Benosman,  Mouhacine and Pappas,  George J.},
  year = {2022},
  month = dec,
  pages = {6318–6325}
}

@article{Ruprecht2012,
	title = {Explicit parallel-in-time integration of a linear acoustic-advection system},
	volume = {59},
	ISSN = {0045-7930},
	url = {http://dx.doi.org/10.1016/j.compfluid.2012.02.015},
	DOI = {10.1016/j.compfluid.2012.02.015},
	journal = {Computers \& Fluids},
	publisher = {Elsevier BV},
	author = {Ruprecht,  D. and Krause,  R.},
	year = {2012},
	month = apr,
	pages = {72–83}
}

@article{Southworth2019,
	title = {Necessary Conditions and Tight Two-level Convergence Bounds for Parareal and Multigrid Reduction in Time},
	volume = {40},
	ISSN = {1095-7162},
	url = {http://dx.doi.org/10.1137/18M1226208},
	DOI = {10.1137/18m1226208},
	number = {2},
	journal = {SIAM Journal on Matrix Analysis and Applications},
	publisher = {Society for Industrial & Applied Mathematics (SIAM)},
	author = {Southworth,  Ben S.},
	year = {2019},
	month = jan,
	pages = {564–608}
}

@article{Su2016,
	author  = {Weijie Su and Stephen Boyd and Emmanuel J. Cand{{\`e}}s},
	title   = {A Differential Equation for Modeling {N}esterov's Accelerated Gradient Method: Theory and Insights},
	journal = {Journal of Machine Learning Research},
	year    = {2016},
	volume  = {17},
	number  = {153},
	pages   = {1--43},
	url     = {http://jmlr.org/papers/v17/15-084.html}
}

@article{Suykens1999,
  title = {Least Squares Support Vector Machine Classifiers},
  volume = {9},
  ISSN = {1573-773X},
  url = {http://dx.doi.org/10.1023/A:1018628609742},
  DOI = {10.1023/a:1018628609742},
  number = {3},
  journal = {Neural Processing Letters},
  publisher = {Springer Science and Business Media LLC},
  author = {Suykens,  J.A.K. and Vandewalle,  J.},
  year = {1999},
  month = jun,
  pages = {293–300}
}

@misc{Tai2021,
  doi = {10.48550/ARXIV.2104.12647},
  note = {\\ \url{https://arxiv.org/abs/2104.12647}},
  author = {Tai,  Yunpeng},
  title = {A Survey Of Regression Algorithms And Connections With Deep Learning},
  publisher = {arXiv},
  year = {2021},
  copyright = {Creative Commons Attribution 4.0 International}
}

@article{Tran2015,
  title = {An $L^1$ Penalty Method for General Obstacle Problems},
  volume = {75},
  ISSN = {1095-712X},
  url = {http://dx.doi.org/10.1137/140963303},
  DOI = {10.1137/140963303},
  number = {4},
  journal = {SIAM Journal on Applied Mathematics},
  publisher = {Society for Industrial & Applied Mathematics (SIAM)},
  author = {Tran,  Giang and Schaeffer,  Hayden and Feldman,  William M. and Osher,  Stanley J.},
  year = {2015},
  month = jan,
  pages = {1424–1444}
}

@BOOK{Trottenberg2000,
	title     = "Multigrid",
	author    = "Trottenberg, Ulrich and Oosterlee, Cornelius W and Schuller,
	Anton",
	publisher = "Academic Press",
	month     =  nov,
	year      =  2000,
	address   = "San Diego, CA",
	language  = "en"
}

@inbook{Ulbrich2015,
	title = {Preconditioners Based on “Parareal” Time-Domain Decomposition for Time-Dependent PDE-Constrained Optimization},
	ISBN = {9783319233215},
	ISSN = {2191-3048},
	url = {http://dx.doi.org/10.1007/978-3-319-23321-5_8},
	DOI = {10.1007/978-3-319-23321-5_8},
	booktitle = {Multiple Shooting and Time Domain Decomposition Methods},
	publisher = {Springer International Publishing},
	author = {Ulbrich,  Stefan},
	year = {2015},
	pages = {203–232}
}

@misc{Vuchkov2024,
	doi = {10.48550/ARXIV.2405.04808},
	note = {\\ \url{https://arxiv.org/abs/2405.04808}},
	author = {Vuchkov,  Radoslav and Cyr,  Eric C. and Ridzal,  Denis},
	title = {Multigrid-in-time preconditioners for {KKT} systems},
	publisher = {arXiv},
	year = {2024},
	copyright = {Creative Commons Attribution 4.0 International}
}

@article{Wu2015,
	title = {Multigrid Methods with {N}ewton-{G}auss-{S}eidel Smoothing and Constraint Preserving Interpolation for Obstacle Problems},
	volume = {8},
	ISSN = {2079-7338},
	url = {http://dx.doi.org/10.4208/nmtma.2015.w08si},
	DOI = {10.4208/nmtma.2015.w08si},
	number = {2},
	journal = {Numerical Mathematics: Theory,  Methods and Applications},
	publisher = {Global Science Press},
	author = {Wu,  Chunxiao and Wan,  Justin W.L.},
	year = {2015},
	month = may,
	pages = {199–219}
}

\end{document}